\documentclass[11pt,a4paper]{article}
\usepackage[left=2cm, top=2.5cm,bottom=2.5cm,right=2cm]{geometry}
\usepackage[T1]{fontenc}
\usepackage[utf8]{inputenc}
\usepackage[british]{babel}
\usepackage{lmodern,mathtools,amssymb,amsthm,mathrsfs,bm,xcolor}
\usepackage{microtype}
\usepackage[hidelinks]{hyperref}
\usepackage{cleveref}
\usepackage{graphicx}
\hypersetup{
  pdftitle={Poisson hyperplanes in constant curvature spaces},
  pdfauthor={Anna Gusakova, Zakhar Kabluchko, Christoph Thale},
  pdfkeywords={Poisson hyperplane tessellation, constant curvature, cell intensity}
}
\usepackage{enumitem}

\numberwithin{equation}{section}

\theoremstyle{plain}
\newtheorem{theorem}{Theorem}[section]
\newtheorem{lemma}[theorem]{Lemma}
\newtheorem{corollary}[theorem]{Corollary}
\newtheorem{proposition}[theorem]{Proposition}

\theoremstyle{definition}
\newtheorem{definition}[theorem]{Definition}
\newtheorem{example}[theorem]{Example}

\theoremstyle{remark}
\newtheorem{remark}[theorem]{Remark}

\newcommand{\dint}{\textup{d}}

\def\EE{\mathbb{E}}

\def\HH{\mathbb{H}}

\def\MM{\mathbb{M}}
\def\NN{\mathbb{N}}

\def\PP{\mathbb{P}}

\def\RR{\mathbb{R}}
\def\SS{\mathbb{S}}

\def\cA{\mathcal{A}}

\def\cC{\mathcal{C}}
\def\cD{\mathcal{D}}
\def\cE{\mathcal{E}}
\def\cF{\mathcal{F}}

\def\cH{\mathcal{H}}

\newcommand{\arccosh}{\mathop{\mathrm{arccosh}}\nolimits}

\newcommand{\proj}{\operatorname{proj}}

\newcommand{\PHT}{\operatorname{\kappa-PHT}}

\newcommand{\cs}{{\rm cs}}
\newcommand{\sn}{{\rm sn}}

\makeatletter
\let\@fnsymbol\@alph
\makeatother

\begin{document}

\title{\bfseries Poisson hyperplanes in constant curvature spaces}

\author{Anna Gusakova\footnotemark[1],\; Zakhar Kabluchko\footnotemark[2],\; Christoph Th\"ale\footnotemark[3]}

\date{}
\renewcommand{\thefootnote}{\fnsymbol{footnote}}
\footnotetext[1]{University of M\"unster, Germany. Email: gusakova@uni-muenster.de}
\footnotetext[2]{University of M\"unster, Germany. Email: zakhar.kabluchko@uni-muenster.de}
\footnotetext[3]{Ruhr University Bochum, Germany. Email: christoph.thaele@rub.de}

\maketitle

\begin{abstract}
	We study isotropic Poisson hyperplane tessellations in spaces of
	constant sectional curvature. Our main result is an explicit formula
	for the cell intensity in arbitrary dimensions, valid in Euclidean,
	elliptic and hyperbolic space. In the hyperbolic case, the formula
	holds throughout the regime when all cells are bounded, including the critical
	intensity. The proof combines the Crofton formula, a mass transport
	principle and the polyhedral Gauss--Bonnet formula, providing a common
	approach to all three curvature regimes. We also derive explicit
	mean-value formulas for the intrinsic volumes, Crofton functionals
	and face numbers of the typical cell. For the zero cell, we determine
	its mean volume and establish identities relating its face numbers,
	intrinsic volumes and total face volumes. In particular, at the
	hyperbolic critical intensity the typical cell has finite mean volume
	but infinite second volume moment. Finally, we identify typical and
	volume-weighted typical lower-dimensional faces with cells of
	intrinsic sectional Poisson hyperplane tessellations, yielding formulas for their
	intensities and mean volumes.
	
\end{abstract}

\begingroup
\small
\setlength{\parskip}{0pt}
\tableofcontents
\endgroup

\section{Introduction and overview of main results}
\label{sec:introduction}

Poisson hyperplane tessellations are classical objects of
stochastic geometry. In Euclidean space, their cell intensities
and many mean geometric characteristics have been known for
decades, see \cite{SW,SKM,HugSchneiderBook}. On the sphere,
the corresponding cell counts are closely related to
Schl\"afli's formula and the combinatorics of random great
hyperspheres, see
\cite{CE67,Miles,ArbeiterZahle,HS16,KT21,HR17,SchneiderConvexCones}.
The situation in hyperbolic space is substantially different.
Its exponential volume growth makes the boundary of a large
observation window comparable in size to its interior,
preventing a direct transfer of the usual Euclidean
large-window arguments.

This raises the question of whether the cell intensity of a
hyperbolic Poisson hyperplane tessellation admits an explicit
formula, and whether such a formula can be derived within a
common framework for spaces of constant curvature. In Euclidean
space, a classical approach uses large observation windows,
where boundary contributions become negligible compared with
volume, see \cite[Section~10.3]{SW}. We develop a different
argument based on Crofton formulas, mass transport and the
polyhedral Gauss--Bonnet formula. It provides an alternative
derivation even in Euclidean space and applies uniformly to all
three curvature regimes. In particular, it allows us to determine the
hyperbolic cell intensity in the regime when all cells are bounded,
including the critical intensity.

\begin{figure}[htbp]
\centering
\includegraphics[width=\textwidth]{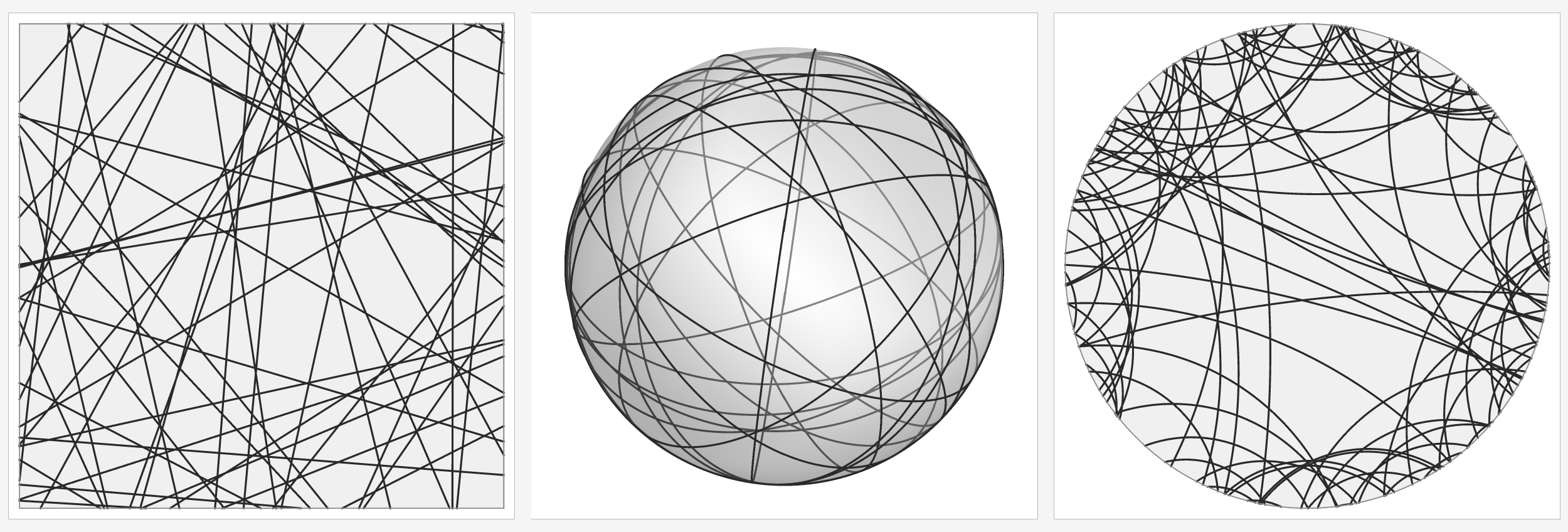}
\caption{Realizations of Poisson hyperplane tessellations
in dimension $d=2$. From left to right: the Euclidean
plane, the sphere, and the hyperbolic plane in the
Poincar\'e disk model. The elliptic plane is obtained
from the sphere by identifying antipodal points.}
\label{fig:PoissonTessellations}
\end{figure}

\paragraph{Cell intensity.}
For $d\geq2$ and $\kappa\in\RR$, let $\MM_\kappa^d$ denote the
$d$-dimensional space of constant sectional curvature $\kappa$,
namely Euclidean space for $\kappa=0$, hyperbolic space for
$\kappa<0$ and elliptic space for $\kappa>0$. Let $\eta_t^\kappa$ be an isometry-invariant Poisson
process of totally geodesic hyperplanes with intensity measure
$t\widetilde\mu_{d-1}^\kappa$, where $t>0$ and
$\widetilde\mu_{d-1}^\kappa$ is
the invariant measure on the space of hyperplanes, whose normalization is specified in
Section~\ref{subsec:space-hyperplanes}.
The resulting random tessellation is called the
\emph{Poisson hyperplane tessellation} and is denoted
by $\PHT_t$. We write $\gamma_{d,t,\kappa}$ for its
cell intensity.
Informally, the cell intensity is the mean number of cells per unit volume.
The precise definition is given in Section~\ref{subsec:typical-cell}.
In positive curvature, cells are understood modulo the antipodal
identification, and if the hyperplane process is empty, we adopt
the convention that the tessellation has no cells. Figure~\ref{fig:PoissonTessellations} illustrates
realizations of Poisson hyperplane tessellations in
the three curvature regimes.

The main question is whether the cell intensity can be determined
explicitly in hyperbolic space, and whether its computation can be
placed within a common framework for all three curvature regimes.
Our first result answers both questions.

\begin{theorem}[Cell intensity of the Poisson hyperplane tessellation]
	\label{thm:cell-intensity}
	Let $d\geq2$, $\kappa\in\RR$ and $t>0$. If $\kappa<0$, assume
	in addition that
	$t\geq t_{\rm crit}(d,\kappa):=\pi(d-1)\sqrt{-\kappa}$. Then
	\begin{equation}\label{eq:IntroCellIntensity}
		\gamma_{d,t,\kappa}
		=
		\frac{2}{\omega_{d+1}}
		\sum_{j=0}^{\lfloor d/2\rfloor}
		\kappa^j\frac{t^{d-2j}}{(d-2j)!}
		-
		\begin{cases}
			\displaystyle
			\frac{2\kappa^{d/2}}{\omega_{d+1}}e^{-t/\sqrt{\kappa}},
			&\kappa>0\text{ and }d\text{ even},\\[0.7em]
			0,&\text{otherwise}.
		\end{cases}
	\end{equation}
	Here $\omega_m=2\pi^{m/2}/\Gamma(m/2)$ is the surface area of the
	unit sphere in $\RR^m$, and the coefficient with $j=0$ is one
	also when $\kappa=0$.
\end{theorem}

For $\kappa=0$, \eqref{eq:IntroCellIntensity} reduces to the
classical Euclidean formula
\[
\gamma_{d,t,0}
=\frac{2t^d}{\omega_{d+1}d!},
\]
see \cite[Equation~(10.46)]{SW}, with the corresponding
normalization of the invariant hyperplane measure.

For $\kappa>0$, the total number $N$ of hyperplanes is Poisson
distributed with mean $t/\sqrt{\kappa}$. These hyperplanes lift
to great hyperspheres on the spherical cover of $\MM_\kappa^d$.
The closures of the connected components of their complement
are called spherical chambers. They occur in antipodal pairs,
each of which defines one cell of the elliptic tessellation.
Conditional on $N=n\geq1$, Schl\"afli's formula therefore gives
$\sum_{i=0}^d\binom{n-1}{i}$ elliptic cells almost surely.
Averaging this count yields \eqref{eq:IntroCellIntensity}.
The corresponding calculation appears in Hug and
Reichenbacher \cite[Theorem~9.5 and its proof]{HR17}.
The exponential correction in~\eqref{eq:IntroCellIntensity} depends on the treatment of the empty
process: counting the whole elliptic space as one cell would
remove the correction in even dimensions and introduce a positive
correction in odd dimensions. It disappears in every dimension
under the Euler-characteristic convention, which assigns the
empty arrangement the weight $\chi(\mathbb{RP}^d)$, see the end of Subsection~\ref{subsec:further-consequences} for details.

For $\kappa<0$, the parameter $t_{\rm crit}(d,\kappa)$ is the
threshold for boundedness of the cells. Above this threshold,
all cells are almost surely bounded according to
\cite[Lemma~5.3.3]{HeroldDiss} or
\cite[Theorem~3.3]{BetaStarPolytopes}. This remains true at the
critical intensity by
\cite[Theorem~1.1 and Corollary~1.2]{BuehlerGusakovaRecke}.
Thus \eqref{eq:IntroCellIntensity} gives the cell intensity
throughout the bounded-cell regime, including its critical
endpoint.

In Euclidean space, a classical approach to computing cell
intensities uses expanding observation windows, where boundary
contributions become negligible compared with volume, see
\cite[Section~10.3]{SW}. This argument does not extend directly
to hyperbolic space, where the boundary area of a large geodesic
ball grows at the same exponential rate as its volume.
In elliptic space, Schl\"afli's formula permits a global count
of chambers, but this approach has no hyperbolic counterpart.
The main feature of our proof is that neither large-window limits
nor global chamber counts are needed. Instead, we determine the
cell intensity from the geometry of the lower-dimensional
skeletons. Crofton and intersection formulas first give their
mean volume densities. The partition of exterior angles,
combined with a mass transport argument, then relates these
densities to expected $\kappa$-intrinsic volumes of the typical cell.
Finally, the polyhedral Gauss--Bonnet formula expresses the
cell intensity in terms of the skeleton densities, with
coefficients determined by the curvature. In Euclidean space,
this reduces to the equality between cell and vertex intensities.
For nonzero curvature, additional skeleton densities contribute.
The resulting argument applies uniformly to all three curvature
regimes and provides an alternative derivation even in the
classical Euclidean case. Mean-value relations for topological and geometric
characteristics of random cell complexes have been studied
in Euclidean and spherical spaces, see
\cite{Zaehle1987,Zaehle1988,LeistritzZaehle1992,ArbeiterZahle}. In particular, the spherical theory
includes relations between curvature-measure intensities
and formulas for sections and great-hypersphere mosaics.
In the hyperbolic plane ($d=2$), a related combination of
Gauss--Bonnet and the mass transport principle appears in
\cite[Theorems~5.2 and~5.5]{BenjaminiSchramm}.
Our argument combines the higher-dimensional polyhedral
Gauss--Bonnet formula with the invariant mass transport principle
to determine the cell intensity uniformly in all three
curvature regimes.

\paragraph{Mean geometric characteristics of the typical cell.}
The \emph{typical cell} $C_{d,t,\kappa}^{\rm typ}$ can be thought
of as a randomly selected cell of the tessellation, with all
cells receiving equal weight, regardless of their size or shape. Its precise definition uses the
Palm distribution of the cell-centre process, as explained in
Section~\ref{subsec:typical-cell}. In particular, this
construction does not require a limiting procedure involving
expanding observation windows. The mean volume of the typical
cell is directly related to the cell intensity, as explained in Proposition~\ref{prop:TypicalCellVolumeIntensity}.

We derive formulas for the expected volume, surface area and
face numbers of the typical cell in all three curvature regimes.
Here we present the particularly simple formulas for
nonpositive curvature. For every integer $r\geq0$, put
\[
P_r(t,\kappa)
:=
\sum_{j=0}^{\lfloor r/2\rfloor}
\kappa^j\frac{t^{r-2j}}{(r-2j)!}.
\]
Write $\cH_\kappa^j$ for $j$-dimensional Hausdorff measure
and $f_k(c)$ for the number of $k$-dimensional faces of a
polytope $c$ in $\MM_\kappa^d$. For $\kappa\leq0$, under
the assumptions of Theorem~\ref{thm:cell-intensity}, we obtain
\begin{equation}\label{eq:IntroTypicalCellVolumeSurface}
	\begin{aligned}
		\EE\cH_\kappa^d(C_{d,t,\kappa}^{\rm typ})
		&=\frac{\omega_{d+1}}{2P_d(t,\kappa)},\qquad
		\EE\cH_\kappa^{d-1}
		(\partial C_{d,t,\kappa}^{\rm typ})
		=\frac{\omega_d t}{P_d(t,\kappa)},
	\end{aligned}
\end{equation}
and
\begin{equation}\label{eq:IntroTypicalCellFaceVector}
	\EE f_k(C_{d,t,\kappa}^{\rm typ})
	=
	2^{d-k}\frac{t^{d-k}}{(d-k)!}
	\frac{P_k(t,\kappa)}{P_d(t,\kappa)},
	\qquad 0\leq k<d.
\end{equation}
The proofs are given in
Sections~\ref{sec:CellIntensityTypicalCell} and
\ref{subsec:typical-cell-face-vector}. All these expectations are finite for every $t>0$ when
$\kappa=0$, and for $t\geq t_{\rm crit}(d,\kappa)$ when
$\kappa<0$, including at the critical intensity.

For $\kappa=0$, we have
$P_r(t,0)=t^r/r!$, and
\eqref{eq:IntroTypicalCellFaceVector} reduces to
$\EE f_k(C_{d,t,0}^{\rm typ})=2^{d-k}\binom{d}{k}$,
the face vector ($f$-vector) of the $d$-dimensional cube.
This classical Euclidean identity also holds for
stationary nondegenerate Poisson hyperplane tessellations
without an isotropy assumption, see
\cite[Theorems~10.3.1 and~10.3.2]{SW}.

In the hyperbolic plane $\MM_\kappa^2$, $\kappa<0$,
let $A$, $L$ and $N$ denote the area, perimeter and
number of sides of the typical cell. Equations
\eqref{eq:IntroTypicalCellVolumeSurface} and
\eqref{eq:IntroTypicalCellFaceVector} yield
\begin{equation}\label{eq:IntroPlanarTypicalCell}
	\EE A=\frac{4\pi}{t^2+2\kappa},
	\qquad
	\EE L=\frac{4\pi t}{t^2+2\kappa},
	\qquad
	\EE N=\frac{4t^2}{t^2+2\kappa},
\end{equation}
for $t\geq\pi\sqrt{-\kappa}$. These expressions agree
with \cite[Equation~(7.11)]{SantaloYanez}, whose
line-intensity parameter is $\lambda=t/\pi$.
That work also derives all six second-order moments
of area, perimeter and side number, see
\cite[Equations~(7.12)--(7.17)]{SantaloYanez}.
However, its definition of an unweighted random polygon
uses an assumed limiting sampling procedure among polygons
meeting expanding observation domains. The derivation of
the corresponding mean values also relies on a boundary-effect
relation that was left unproved in hyperbolic space,
see \cite[Equation~(7.9)]{SantaloYanez}.
Our Palm construction establishes
\eqref{eq:IntroPlanarTypicalCell} in the bounded-cell
regime without these assumptions.
A detailed comparison is given in
Remark~\ref{rem:SantaloYanezTypicalCell}.

The preceding identities belong to a more general family
of mean-value formulas involving curvature-dependent
intrinsic volumes. For a compact full-dimensional convex
polytope $c\subset\MM_\kappa^d$ with $\kappa\leq0$,
we define
\[
\upsilon_k^\kappa(c)
:=
\frac{1}{\omega_{k+1}}
\sum_{F\in\cF_k(c)}
\alpha(F,c)\cH_\kappa^k(F),
\qquad 0\leq k\leq d,
\]
where $\alpha(F,c)$ denotes the normalized exterior angle
at the face $F$. In top dimension, we set
$\alpha(c,c)=1$. These functionals coincide, up to
normalization, with the usual intrinsic volumes when
$\kappa=0$. Their precise definition for all three
curvature regimes is given in
Section~\ref{subsec:MeanIntrinsicVolumes}.

For $\kappa=0$ and $t>0$, or for $\kappa<0$ and
$t\geq t_{\rm crit}(d,\kappa)$,
Theorem~\ref{thm:MeanIntrinsicVolumes} gives
\[
\EE\upsilon_k^\kappa(C_{d,t,\kappa}^{\rm typ})
=
\frac{t^{d-k}}{2(d-k)!\,P_d(t,\kappa)},
\qquad 0\leq k\leq d.
\]
The theorem also provides the corresponding formulas for
positive curvature. We further determine the mean Crofton
functionals, or quermassintegrals, which measure the sets of
totally geodesic flats meeting a cell. Their expectations
are given explicitly in
Theorem~\ref{thm:MeanTypicalCellCrofton}.
The relation between these functionals and the
$\kappa$-intrinsic volumes is established in
Proposition~\ref{prop:CroftonIntrinsicTransform}.

\paragraph{Zero cells and lower-dimensional faces.}
A second natural way to select a cell is to fix a reference
point $p\in\MM_\kappa^d$ and consider the cell
$C_{d,t,\kappa}^p$ containing it, called the
\emph{zero cell}. Unlike the typical cell, this selection
favours larger cells. More precisely, up to isometry,
the distribution of the zero cell coincides with the volume-weighted typical-cell distribution. In positive curvature, this identity holds
after conditioning on the hyperplane process being nonempty.
On the empty-process event, we define the zero cell to be
the whole elliptic space, although this realization
contributes no cell to the tessellation. The definitions
and the volume-weighting identity are given in
Section~\ref{subsec:zero-cell-definition} and
Proposition~\ref{prop:ZeroCellVolumeBias}.

The zero cell can be studied for every $t>0$, even when
the hyperbolic tessellation contains unbounded cells.
Using the void probability formula for Poisson processes, we obtain a common
integral representation for its mean volume in all three
curvature regimes and evaluate it explicitly in terms of
finite products and Gamma functions, see
Theorems~\ref{thm:ZeroCellMeanVolume} and~\ref{thm:ZeroCellVolumeClosedForm}.
For hyperbolic space of curvature $-1$, the mean-volume
formula and its finiteness criterion were previously
obtained in
\cite[Theorem~7.1]{BuehlerHugThaeleVisibility}.
The planar result already appears in
\cite[Equation~(6.4)]{SantaloYanez}.

The volume-weighting identity also relates the zero cell
to the second volume moment of the typical cell.
Write
$V^{\rm typ}=\cH_\kappa^d(C_{d,t,\kappa}^{\rm typ})$
and $V^0=\cH_\kappa^d(C_{d,t,\kappa}^p)$.
For $\kappa\leq0$, under the assumptions of
Theorem~\ref{thm:cell-intensity}, we have
\[
\EE[(V^{\rm typ})^2]
=
\EE V^{\rm typ}\,\EE V^0,
\]
as an identity in $[0,\infty]$, see
Proposition~\ref{prop:TypicalCellVolumeSecondMoment}.
In hyperbolic space, the mean zero-cell volume is finite
precisely when
$t>t_{\rm crit}(d,\kappa)$.
Consequently, at the critical intensity,
\[
0<\EE V^{\rm typ}<\infty,
\qquad
\EE V^0=\EE[(V^{\rm typ})^2]=+\infty.
\]
Thus the typical cell has finite mean volume but infinite
variance, while the zero cell has infinite mean volume,
even though both cells are almost surely bounded.

Beyond volumes, we establish an Efron-type identity
relating the mean Crofton functionals of the zero cell
to its expected face numbers, see
Theorem~\ref{thm:ZeroCellEfron}. Related face-number
results are available for Euclidean Poisson zero cells,
weighted spherical cells generated by a fixed number of
great hyperspheres and hyperbolic Poisson zero cells
in \cite{KabluchkoZeroFVector,KT21,BetaStarPolytopes},
respectively. Combining the hyperbolic face-number formula
with our identity yields explicit expressions for the mean
proper $\kappa$-intrinsic volumes of the zero cell.
In particular, each such mean is obtained from its
Euclidean counterpart at the same intensity parameter
by multiplication with a factor independent of the
intrinsic-volume index, see
Theorem~\ref{prop:HyperbolicZeroCellIntrinsicVolumes} and Remark~\ref{rem:ZeroCellIntrinsicVolumeCounterparts}.
We also determine the mean total volume of faces of
the zero cell in
Corollary~\ref{cor:ZeroCellSkeletonVolume}.

The relation between typical and volume-weighted cells
has a counterpart for lower-dimensional faces.
We define the typical $k$-face by sampling from the
collection of distinct $k$-faces of the
tessellation, rather than from the faces of the typical
full-dimensional cell. For $1\leq k\leq d-1$, we show
that, up to isometry, its distribution agrees with that
of the typical cell of an intrinsic $k$-dimensional
Poisson hyperplane tessellation with the same intensity
parameter $t$. Moreover, the volume-weighted typical
$k$-face has the distribution of the intrinsic
$k$-dimensional zero cell. In positive curvature,
the intrinsic comparison tessellation counts an uncut
whole $k$-space as one cell when the induced hyperplane
process is empty. These representations are established
in Theorem~\ref{thm:TypicalKFaceSection}.
The corresponding Euclidean representation for
volume-weighted faces was obtained in
\cite{SchneiderWeightedFaces}, and related fixed-number
spherical results appear in \cite{KT21}.
Our sectional representations yield explicit face
intensities and show that the mean $k$-volume of a
typical $k$-face depends only on $k$, $t$ and $\kappa$,
but not on the ambient dimension $d$, under the stated
bounded-cell assumptions, see
Corollary~\ref{cor:MeanTypicalKFaceVolume}.

Finally, the explicit formulas allow us to study the
dependence of these geometric characteristics on the
curvature. For fixed $d$ and $t$, the cell intensity,
the mean volume and mean intrinsic and Crofton functionals
of the typical cell, and the mean volume of the zero cell
are $C^\infty$ functions of $\kappa$ in a neighbourhood
of zero, see Propositions~\ref{prop:CurvatureSmoothness}
and \ref{prop:ZeroCellVolumeSmoothness}.
This dependence need not be real analytic. Indeed,
the positive-curvature formulas contain corrections
involving $e^{-t/\sqrt{\kappa}}$ multiplied by powers
of $\kappa$. Extended by zero to nonpositive curvature,
these terms are smooth but have all derivatives vanishing
at the origin. Thus the explicit formulas exhibit a
smooth transition between the three curvature regimes
which is not, in general, described by a convergent
power series.

\medspace

The remainder of the paper is organized as follows.
Section~\ref{sec:random-cells} introduces the
constant-curvature spaces, Poisson hyperplane processes
and the relevant random cells.
Section~\ref{sec:CellIntensityTypicalCell} establishes
the cell-intensity formula and the mean geometric
characteristics of the typical cell.
Section~\ref{sec:ZeroCell} studies the zero cell,
including its volume moments, face numbers and intrinsic
volumes.
Section~\ref{sec:Kfaces} develops the sectional
representations of typical faces and derives the
corresponding face intensities and mean-value formulas.

\section{Poisson hyperplane tessellations and their cells}
\label{sec:random-cells}

\subsection{Constant-curvature spaces and Poisson hyperplane processes}
\label{subsec:space-hyperplanes}
\label{subsec:invariant-measures}

\paragraph{Constant-curvature spaces.}
Throughout, let $d\geq2$ and $\kappa\in\RR$. We first specify the
models used for the $d$-dimensional Euclidean, hyperbolic and elliptic
spaces. Let $g_{{\rm euc}}$ be the Euclidean metric on $\RR^d$.
Let
\[
\SS^d:=\{x\in\RR^{d+1}:x_0^2+x_1^2+\cdots+x_d^2=1\}
\]
be the unit $d$-dimensional sphere.
The $d$-dimensional elliptic space is the real projective
space
\[
\cE^d:=\SS^d/\{\pm{\rm id}\}=\mathbb{RP}^d
\]
equipped with the quotient Riemannian metric
$g_{{\rm ell}}$ induced by the standard spherical metric
on $\SS^d$.
For hyperbolic space we use the hyperboloid model
\[
\HH^d
:=
\{x\in\RR^{d+1}:-x_0^2+x_1^2+\cdots+x_d^2=-1,\ x_0>0\},
\]
with the Riemannian metric $g_{{\rm hyp}}$ induced by the Lorentzian
bilinear form
\[
b_{-1}(x,y):=-x_0y_0+x_1y_1+\cdots+x_dy_d,
\qquad x,y\in\RR^{d+1}.
\]
For arbitrary curvature $\kappa \in \RR$, we set
\begin{equation}\label{eq:Volumescale}
(\MM_\kappa^d,g_\kappa)
=
\begin{cases}
(\cE^d,\kappa^{-1}g_{{\rm ell}}),&\kappa>0,\\
(\RR^d,g_{{\rm euc}}),&\kappa=0,\\
(\HH^d,(-\kappa)^{-1}g_{{\rm hyp}}),&\kappa<0.
\end{cases}
\end{equation}
The Riemannian manifold $(\MM_\kappa^d,g_\kappa)$ has constant sectional
curvature $\kappa$, see \cite[p.~236]{Lang}.
We denote by $d_\kappa$ the distance induced by $g_\kappa$.
Explicitly,
\begin{equation*}
d_\kappa(x,y)
=
\begin{cases}
\kappa^{-1/2}\arccos|\langle\widehat x,\widehat y\rangle|,
&\kappa>0,\\[0.4em]
\|x-y\|,
&\kappa=0,\\[0.4em]
(-\kappa)^{-1/2}\arccosh\bigl(-b_{-1}(x,y)\bigr),
&\kappa<0,
\end{cases}
\end{equation*}
where, in the elliptic case, $\widehat x,\widehat y\in\SS^d$ are
arbitrary representatives of $x,y\in\cE^d$. The absolute value makes
the first expression independent of these choices.

For $0\leq j\leq d$, let $\cH_\kappa^j$ denote the $j$-dimensional
Hausdorff measure induced by $d_\kappa$. In particular,
$\cH_\kappa^d$ is Riemannian volume and $\cH_\kappa^0$ is counting
measure.
As in the
introduction, we write $\omega_m:=2\pi^{m/2}/\Gamma(m/2)$ for the surface area of the unit sphere in $\RR^m$. In positive curvature $\kappa>0$, we also have
\begin{equation}\label{eq:EllipticSpaceVolume}
\cH_\kappa^d(\MM_\kappa^d)
=
\frac{\omega_{d+1}}{2\kappa^{d/2}}.
\end{equation}
For what follows, we fix a reference point $p\in\MM_\kappa^d$, the ``origin'', and we write
$I(\MM_\kappa^d)$ for the full isometry group of $(\MM_\kappa^d,g_\kappa)$.

\paragraph{Grassmannians and invariant measures.}
For $0\leq k\leq d$, let $A_\kappa(d,k)$ be the space of complete
totally geodesic $k$-dimensional subspaces, called $k$-flats, of
$\MM_\kappa^d$. All flats are unoriented. We identify
$A_\kappa(d,0)$ with $\MM_\kappa^d$ and have
$A_\kappa(d,d)=\{\MM_\kappa^d\}$. Let
$G_\kappa(d,k)\subset A_\kappa(d,k)$ consist of those $k$-flats
containing the reference point $p$. The spaces $A_\kappa(d,k)$ and $G_\kappa(d,k)$ are equipped with their natural Borel
$\sigma$-algebras. We denote by $\nu_k^\kappa$ the unique probability
measure on $G_\kappa(d,k)$ invariant under all isometries of $\MM_\kappa^d$ fixing $p$.

For $0\leq k<d$, $L\in G_\kappa(d,d-k)$ and $x\in L$, let
$H(L,x)$ be the unique totally geodesic $k$-flat through $x$
orthogonal to $L$ at $x$. For every Borel set
$A\subset A_\kappa(d,k)$, define
\begin{equation}\label{eq:FlatMeasureGeneralCurv}
\mu_k^\kappa(A)
=
\int_{G_\kappa(d,d-k)}\int_L
\cs_\kappa\bigl(d_\kappa(x,p)\bigr)^k
{\bf1}\{H(L,x)\in A\}\,
\cH_\kappa^{d-k}(\dint x)\,
\nu_{d-k}^\kappa(\dint L),
\end{equation}
where
\[
\cs_\kappa(r)
=
\begin{cases}
\cos(\sqrt{\kappa}\,r),&\kappa>0,\\
1,&\kappa=0,\\
\cosh(\sqrt{-\kappa}\,r),&\kappa<0,
\end{cases}
\]
is the generalized cosine function.
When $k=0$, the factor involving $\cs_\kappa$ equals
$1$, and \eqref{eq:FlatMeasureGeneralCurv} gives
$\mu_0^\kappa=\cH_\kappa^d$.

The measure $\mu_k^\kappa$ is a Radon measure invariant under the
natural action of the isometry group $I(\MM_\kappa^d)$ on $A_\kappa(d,k)$ and does not
depend on the choice of $p$. Up to multiplication by a positive
constant, it is the unique nonzero isometry-invariant Radon measure
on $A_\kappa(d,k)$, see \cite[Equation~(17.41)]{SantaloBook} for the
construction. For $\kappa=-1$, the normalization in
\eqref{eq:FlatMeasureGeneralCurv} agrees with
\cite[Equation~(6)]{HeroldHugThaele} and its relation to the normalization
in \cite{SantaloBook} is explained in
\cite[Remark~7]{HeroldHugThaele}.

Throughout the paper we use the renormalized measures
\begin{equation}\label{eq:FlatMeasureNormalization}
\widetilde\mu_k^\kappa
:=
\frac{\omega_{k+1}}{\omega_{d+1}}\mu_k^\kappa,
\qquad 0\leq k<d, \qquad
\widetilde\mu_d^\kappa:=\delta_{\MM_\kappa^d}.
\end{equation}
The factors in \eqref{eq:FlatMeasureNormalization} simplify the
constants in the Crofton and intersection formulae used below.  In particular, under this normalization the intersection formula \eqref{eq:Formula_1} has no additional dimensional factor.
Another advantage of this normalization is that, for $\kappa=1$, the measure $\widetilde\mu_k^1$ is a probability
measure on $A_1(d,k)$. Indeed, for $0\leq k<d$, polar coordinates
in $L$ give
\[
\widetilde\mu_k^1(A_1(d,k))
=
\frac{\omega_{k+1}\omega_{d-k}}{\omega_{d+1}}
\int_0^{\pi/2}
(\cos r)^k(\sin r)^{d-k-1}\,\dint r
=1,
\]
and the assertion for $k=d$ follows from the definition.
Under the identifications of the underlying spaces in
\eqref{eq:Volumescale}, rescaling gives
\begin{equation}\label{eq:HyperplaneMeasureScale}
\widetilde\mu_k^\kappa
=
\begin{cases}
\kappa^{-(d-k)/2}\widetilde\mu_k^1,&\kappa>0,\\
(-\kappa)^{-(d-k)/2}\widetilde\mu_k^{-1},&\kappa<0,
\end{cases}
\qquad 0\leq k\leq d.
\end{equation}
Consequently,
\begin{equation}\label{eq:FlatMeasureScalePositive}
\widetilde\mu_k^\kappa(A_\kappa(d,k))
=
\kappa^{-(d-k)/2},
\qquad \kappa>0,\quad 0\leq k\leq d.
\end{equation}

\begin{remark}
The space $A_\kappa(d,d-1)$ consists of geodesic hyperplanes. For every
minimizing geodesic segment $[x,y]\subset\MM_\kappa^d$, with
$\cH_\kappa^1([x,y])<\pi/(2\sqrt{\kappa})$ when $\kappa>0$, we have
\begin{equation}\label{eq:HyperplaneSegmentCrofton}
\widetilde\mu_{d-1}^\kappa
\bigl(\{H\in A_\kappa(d,d-1):H\cap [x,y]\neq\emptyset\}\bigr)
=
\frac{\cH_\kappa^1([x,y])}{\pi}.
\end{equation}
This is the case $i=0$, $k=d-1$ of the Crofton formula
\eqref{eq:Formula_2} below. Thus \eqref{eq:HyperplaneSegmentCrofton}
provides a geometric interpretation of the normalization in
\eqref{eq:FlatMeasureNormalization}.
\end{remark}

\paragraph{Poisson hyperplane processes.}
For $t>0$, let $\eta_t^\kappa$ be a Poisson process on
$A_\kappa(d,d-1)$ with intensity measure
$t\widetilde\mu_{d-1}^\kappa$. Thus, for every Borel set
$A\subset A_\kappa(d,d-1)$ of finite invariant measure,
$\eta_t^\kappa(A)$ is Poisson distributed with mean
$t\widetilde\mu_{d-1}^\kappa(A)$, and the counts in pairwise disjoint
Borel sets are independent. The law of $\eta_t^\kappa$ is invariant
under the isometry group $I(\MM_\kappa^d)$. By \eqref{eq:HyperplaneSegmentCrofton},
the expected number of hyperplanes meeting a geodesic segment of
length $r$ is $tr/\pi$, subject to the same restriction on the
segment in positive curvature. Almost surely, only finitely many hyperplanes of $\eta_t^\kappa$
meet any compact set of $\MM_\kappa^d$, and the hyperplanes are in general position.
More precisely, every nonempty intersection of $m$ distinct
hyperplanes has dimension $d-m$ for $1\leq m\leq d$, and no $d+1$
distinct hyperplanes have a common point.

For $\kappa\leq0$, the invariant hyperplane measure
$\widetilde\mu_{d-1}^\kappa$ has infinite total mass, and hence
$\eta_t^\kappa$ contains infinitely many hyperplanes almost surely. For $\kappa>0$, \eqref{eq:FlatMeasureScalePositive} with $k=d-1$
shows that the total number of hyperplanes is a.s.\ finite and satisfies
\begin{equation}\label{eq:EllipticHyperplaneCount}
N:=\eta_t^\kappa\bigl(A_\kappa(d,d-1)\bigr)
\sim\operatorname{Poisson}\!\left(\frac{t}{\sqrt{\kappa}}\right).
\end{equation}
The probability that the Poisson hyperplane process $\eta_t^\kappa$ is non-empty will be denoted by
\begin{equation}\label{eq:NonemptyProbability}
p_\kappa(t)
:=
\PP(\eta_t^\kappa\neq\emptyset)
=
\begin{cases}
1-e^{-t/\sqrt{\kappa}},&\kappa>0,\\
1,&\kappa\leq0.
\end{cases}
\end{equation}

\subsection{Cells, faces and the bounded-cell regime}
\label{subsec:cells-faces}

The hyperplanes of the Poisson hyperplane process $\eta_t^\kappa$ determine a random collection of
cells, called the Poisson hyperplane tessellation and denoted by $\PHT_t$. We first describe the construction in
nonpositive curvature and then explain the convention used in
positive curvature.

\paragraph{Nonpositive curvature.}
Suppose that $\kappa\leq0$. The cells of the Poisson hyperplane tessellation $\PHT_t$ are the closures
of the connected components of the open set
\begin{equation}\label{eq:HyperplaneComplement}
\MM_\kappa^d\setminus\bigcup_{H\in\eta_t^\kappa}H.
\end{equation}
The cells  are closed geodesically convex sets with nonempty interior.
By local finiteness of the hyperplane process, each cell is locally
an intersection of finitely many closed geodesic halfspaces.

For $0\leq k<d$, the $k$-faces of a cell $c$ are the nonempty
$k$-dimensional intersections
$
c\cap H_1\cap\cdots\cap H_{d-k},
$
where $H_1,\ldots,H_{d-k}$ are distinct hyperplanes of
$\eta_t^\kappa$. We denote the collection of $k$-faces of $c$  by $\cF_k(c)$ and put
$f_k(c):=\#\cF_k(c)$. The collection of all distinct $k$-faces of
the tessellation is denoted by $\cF_k(\PHT_t)$ with the convention that a face shared by
several cells is included only once.
In Euclidean space, all cells are almost surely bounded for every
$t>0$, see \cite[Theorem~10.3.2]{SW}. In hyperbolic space this is not the case. Define
\begin{equation}\label{eq:HyperbolicCriticalIntensity}
t_{\rm crit}(d,\kappa)
:=
\pi(d-1)\sqrt{-\kappa},
\qquad \kappa<0.
\end{equation}
Almost surely, all cells are  bounded when
$t\geq t_{\rm crit}(d,\kappa)$. For
$t>t_{\rm crit}(d,\kappa)$, this follows from
\cite[Lemma~5.3.3]{HeroldDiss} or
\cite[Theorem~3.3]{BetaStarPolytopes}. Boundedness at
$t=t_{\rm crit}(d,\kappa)$ follows from
\cite[Theorem~1.1 and Corollary~1.2]{BuehlerGusakovaRecke}.
In these cases the cells are compact, and local finiteness of the
hyperplane process implies that they are convex polytopes with
finitely many faces.

\begin{remark}
\label{rem:HyperbolicIntensityNormalization}
For curvature $-1$, the Poisson hyperplane process in
\cite[Section~3.2, Equation~(3.9)]{BetaStarPolytopes}
has intensity measure $\frac{\lambda}{2}\mu_{d-1}^{-1}$,
where $\mu_{d-1}^{-1}$ is defined by
\eqref{eq:FlatMeasureGeneralCurv}. Consequently,
\eqref{eq:FlatMeasureNormalization} gives the parameter conversion
\begin{equation}\label{eq:HyperbolicParameterConversion}
\lambda
=
2\frac{\omega_d}{\omega_{d+1}}\,t.
\end{equation}
The critical value in \cite[Theorem~3.3]{BetaStarPolytopes} is
\[
\lambda_d^{\rm crit}
=
\frac{(d-1)^2\sqrt{\pi}\,
      \Gamma\bigl((d-1)/2\bigr)}
     {\Gamma(d/2)}
=
2\frac{\omega_d}{\omega_{d+1}}\,\pi(d-1).
\]
Thus \eqref{eq:HyperbolicParameterConversion} yields
$t_{\rm crit}(d,-1)=\pi(d-1)$. The expression for arbitrary
$\kappa<0$ in \eqref{eq:HyperbolicCriticalIntensity} follows by
rescaling according to \eqref{eq:HyperplaneMeasureScale}.
\end{remark}

\paragraph{Positive curvature.}
Suppose now that $\kappa>0$. Let $\SS_\kappa^d$ denote the unit
sphere $\SS^d$ equipped with its standard spherical metric
multiplied by $\kappa^{-1}$, and let
$
q_\kappa:\SS_\kappa^d\longrightarrow\MM_\kappa^d
$
be the quotient map identifying antipodal points. We use
$\cH_\kappa^j$ also for the Hausdorff measures on this spherical
cover.

Each hyperplane of $\eta_t^\kappa$ lifts to a great hypersphere
of $\SS_\kappa^d$. On the event $\eta_t^\kappa\neq\emptyset$,
the closures of the connected components of the complement
of the union of these great hyperspheres are called
\emph{spherical chambers}. They are
compact spherical polyhedra and occur in antipodal pairs.
An \emph{elliptic cell}, also called an \emph{elliptic chamber},
is an antipodal pair
$c=\{\widehat c,-\widehat c\}$ of spherical chambers.
Its \emph{support} is
\begin{equation*}
|c|
:=
q_\kappa(\widehat c)
=
q_\kappa(-\widehat c)
\subset\MM_\kappa^d.
\end{equation*}
We define $\PHT_t$ to be the collection of these elliptic cells.
On the event $\eta_t^\kappa=\emptyset$, we set
$\PHT_t=\emptyset$.
The volume of an elliptic cell is
\begin{equation*}
\cH_\kappa^d(c)
:=
\cH_\kappa^d(|c|)
=
\cH_\kappa^d(\widehat c).
\end{equation*}
The last equality holds because $q_\kappa$ is a local isometry
and is one-to-one on the interior of $\widehat c$. Any
identifications within $\widehat c$ occur on its boundary,
which has zero $d$-dimensional volume.

For $0\leq k<d$, choose one spherical representative
$\widehat c$ of an elliptic cell $c$. The spherical
$k$-faces of $\widehat c$ are the $k$-dimensional intersections of
$\widehat c$ with collections of the great hyperspheres
bounding its defining hemispheres. The $k$-faces of the elliptic cell $c$
are the distinct images of these spherical faces under
$q_\kappa$.  Thus, each $k$-face is regarded as a subset of $\MM_\kappa^d$. We write $\cF_k(c)$ for the collection of $k$-faces
of $c$ and put $f_k(c):=\#\cF_k(c)$.
The collection of all $k$-faces of the cells of
$\PHT_t$ is denoted by
\[
\cF_k(\PHT_t):=\bigcup_{c\in\PHT_t}\cF_k(c).
\]
A face shared by several elliptic
cells is included only once in $\cF_k(\PHT_t)$.   If $\eta_t^\kappa=\emptyset$, we set
$\cF_k(\PHT_t)=\emptyset$.

A spherical $k$-face may be a whole great $k$-sphere.
Its projective image is then a complete elliptic
$k$-flat, which is counted as one face. For example,
if $\eta_t^\kappa$ consists of a single hyperplane $H$,
its spherical lift bounds two hemispheres forming one
elliptic cell $c$. In this case, $|c|=\MM_\kappa^d$,
although $H$ is a facet of $c$. Thus an elliptic cell
should not be identified with its support.

\paragraph{Common notation.}
For every curvature $\kappa\in\RR$, we set
\[
\cF_d(\PHT_t):=\PHT_t,
\qquad
\cF_d(c):=\{c\},
\qquad
f_d(c):=1.
\]
For $\kappa\leq0$ we put $|c|:=c$. When no confusion can arise,
we also write $c$ for $|c|$ in set-theoretic expressions and
volume integrals. The same convention applies to lower-dimensional faces:
we write $|F|$ for the geometric support of a face $F$
and use $\cH_\kappa^k(F)$ as shorthand for
$\cH_\kappa^k(|F|)$.

\paragraph{The space of chambers.}
We denote by $\cC_\kappa^d$ the measurable space of compact
full-dimensional polyhedral chambers. For $\kappa\leq0$, its
elements are compact geodesically convex polytopes in
$\MM_\kappa^d$ with nonempty interior. We equip this space
with the Borel $\sigma$-algebra induced by the Hausdorff
metric associated with $d_\kappa$.

For $\kappa>0$, the elements of $\cC_\kappa^d$ are
elliptic chambers, that is, antipodal pairs of spherical
polyhedra obtained as intersections of finite nonempty
collections of closed hemispheres and having nonempty
spherical interior. In particular, the whole sphere $\SS_\kappa^d$ is excluded
from the space of spherical representatives used to define
$\cC_\kappa^d$.  We equip $\cC_\kappa^d$ with
the quotient Borel structure induced by the Hausdorff metric
on the spherical representatives under
$\widehat c\mapsto-\widehat c$.

\paragraph{Default intensity regime.}
Unless stated otherwise, throughout the remainder of the paper
we assume
\begin{equation}\label{eq:AdmissibleRegime}
\begin{aligned}
&t>0 &&\text{for all }\kappa\in\RR,\\
&t\geq t_{\rm crit}(d,\kappa) &&\text{if }\kappa<0.
\end{aligned}
\end{equation}
Under \eqref{eq:AdmissibleRegime}, all cells of $\PHT_t$
belong to $\cC_\kappa^d$ almost surely.

\subsection{Cell intensity and the typical cell}
\label{subsec:typical-cell}

Throughout this subsection we assume \eqref{eq:AdmissibleRegime}.
Recall that $\cC_\kappa^d$ is the measurable space of compact
full-dimensional polyhedral chambers introduced in
Section~\ref{subsec:cells-faces}. For $\kappa\leq0$, its elements
are compact geodesically convex polytopes. For $\kappa>0$, each chamber is an antipodal pair of spherical
polyhedra.

Cell intensity may be understood as an average number of cells per unit volume. To count cells per unit volume, we assign a centre to each chamber.
This is done using a centre map.

\begin{definition}\label{def:CentreMap}
A \emph{centre map} is a measurable map
$
z:\cC_\kappa^d\longrightarrow\MM_\kappa^d
$
such that
\begin{equation}\label{eq:CentreMapProperties}
z(c)\in|c|,
\qquad
z(\phi c)=\phi z(c),
\qquad
c\in\cC_\kappa^d,\quad \phi\in I(\MM_\kappa^d).
\end{equation}
Here $\phi c$ denotes the image of the chamber $c$ under
the isometry $\phi$.
\end{definition}

\begin{lemma}\label{lem:CentreMapExists}
A centre map exists for every $\kappa\in\RR$.
\end{lemma}

\begin{proof}
For $\kappa\leq0$, take the circumcentre,
\[
z(c)
=
\operatorname*{arg\,min}_{x\in\MM_\kappa^d}
\max_{y\in c}d_\kappa(x,y).
\]
In Euclidean and hyperbolic space this minimizer exists, is unique,
and belongs to the compact convex set $c$, see also
\cite[Lemma~5.2.13]{HeroldDiss} for the hyperbolic case.
Uniqueness implies covariance under isometries. The circumcentre
depends continuously on $c$ in the Hausdorff metric: under Hausdorff
convergence, the functions being minimized converge uniformly on
bounded sets, the minimizers remain bounded, and every limit point
is the unique minimizer for the limiting set. In particular, the
map is measurable.

For $\kappa>0$, recall that
$c=\{\widehat c,-\widehat c\}$ is an antipodal pair of spherical
polyhedra and that
$q_\kappa:\SS_\kappa^d\to\MM_\kappa^d$ is the quotient map.
Regarding $\widehat c$ as a subset of the unit sphere in
$\RR^{d+1}$, put
\begin{equation}\label{eq:SphericalCentre}
b(\widehat c)
:=
\int_{\widehat c}u\,\cH_\kappa^d(\dint u),
\qquad
z(c)
:=
q_\kappa\!\left(
\frac{b(\widehat c)}{\|b(\widehat c)\|}
\right).
\end{equation}
Let $n$ be the inward unit normal of any hemisphere defining
$\widehat c$. Since $\widehat c$ has nonempty spherical interior,
$\langle n,u\rangle>0$ for almost every $u\in\widehat c$.
Consequently,
$\langle n,b(\widehat c)\rangle>0$. Applying this to all defining
hemispheres shows that $b(\widehat c)\neq0$ and that its normalized
direction belongs to $\widehat c$.
Replacing $\widehat c$ by $-\widehat c$ changes the sign of
$b(\widehat c)$ and therefore leaves $z(c)$ unchanged. Moreover,
the integral in \eqref{eq:SphericalCentre} commutes with orthogonal
transformations, which proves covariance under elliptic isometries.
It depends continuously on the spherical representative in the
Hausdorff metric: the corresponding indicator functions converge
away from the boundary of the limiting spherical polyhedron, which
has spherical volume zero. Thus \eqref{eq:SphericalCentre} defines
a measurable map on the quotient chamber space.
\end{proof}

Fix a centre map $z:\cC_\kappa^d\longrightarrow\MM_\kappa^d$. The associated \emph{cell-centre process} is
the random counting measure
\begin{equation*}
\xi_t^{\kappa,z}
:=
\sum_{c\in\PHT_t}\delta_{z(c)}.
\end{equation*}
If several cells have the same centre, they contribute separate
atoms to this sum. By \eqref{eq:CentreMapProperties} and the
isometry invariance of $\PHT_t$, the law of $\xi_t^{\kappa,z}$ is
isometry invariant.

\begin{lemma}\label{lem:CellIntensityFinite}
Under \eqref{eq:AdmissibleRegime}, for every centre map $z$ there
is a constant $\gamma_{d,t,\kappa}^{\,z}\in(0,\infty)$ such that
\[
\EE\xi_t^{\kappa,z}(B)
=
\gamma_{d,t,\kappa}^{\,z}\cH_\kappa^d(B)
\]
for every Borel set $B\subset\MM_\kappa^d$.
\end{lemma}

\begin{proof}
First suppose that $\kappa\leq0$. Given a compact set $B$, choose
a compact geodesically convex window $W$ whose interior contains
$B$, and let $N_W$ be the number of hyperplanes meeting $W$.
Every full-dimensional cell whose support meets $B$ has an interior
point in $W$. The $N_W$ hyperplanes meeting $W$ produce at most
$2^{N_W}$ different cell interiors in $W$. Since the centre of a
cell belongs to its support,
\[
\xi_t^{\kappa,z}(B)\leq2^{N_W}.
\]
The random variable $N_W$ is Poisson distributed with finite mean,
and hence $\EE2^{N_W}<\infty$. On the other hand, if $\kappa>0$, let $N$ be the total number of hyperplanes, as in
\eqref{eq:EllipticHyperplaneCount}. The total number of elliptic
cells is at most $2^N$, including the empty-process realization,
on which there are no cells. Again, $\EE2^N<\infty$.

Thus $\EE\xi_t^{\kappa,z}$ is a locally finite isometry-invariant
measure. Uniqueness of Riemannian volume up to normalization gives
the asserted representation with a finite nonnegative constant.
If this constant were zero, the centre process would almost surely
be empty on every member of a countable cover by bounded balls.
It would therefore be empty almost surely. This contradicts
$p_\kappa(t)>0$ in \eqref{eq:NonemptyProbability}, since every
nonempty tessellation has at least one cell and hence at least
one centre.
\end{proof}

We call $\gamma_{d,t,\kappa}^{\,z}$ the \emph{cell intensity}
associated with $z$. Equivalently,
\begin{equation}\label{eq:CellIntensityDefinition}
\gamma_{d,t,\kappa}^{\,z}
=
\frac{1}{\cH_\kappa^d(B)}
\EE\sum_{c\in\PHT_t}{\bf1}\{z(c)\in B\},
\qquad
0<\cH_\kappa^d(B)<\infty.
\end{equation}
Thus it is the mean number of cell centres per unit volume, and
the quotient in \eqref{eq:CellIntensityDefinition} does not depend
on the choice of $B$.

The \emph{typical cell} can be thought of as a randomly selected
cell of the tessellation, with each cell receiving equal weight,
regardless of its size or shape.
The Palm construction provides a rigorous formulation of this
sampling procedure without requiring a limiting sequence of
observation windows. We shall define the typical cell as a random chamber
$C_{d,t,\kappa}^{\mathrm{typ},z}$ satisfying
\begin{equation}\label{eq:TypicalCellInvariantExpectations}
\gamma_{d,t,\kappa}^{\,z}\,
\EE h(C_{d,t,\kappa}^{\mathrm{typ},z})
=
\frac{1}{\cH_\kappa^d(B)}
\EE\sum_{c\in\PHT_t}
{\bf1}\{z(c)\in B\}h(c)
\end{equation}
for every nonnegative measurable isometry-invariant functional
$h:\cC_\kappa^d\to[0,\infty]$ and every Borel set $B$ of positive
finite volume. Here isometry invariance means that
$h(\phi c)=h(c)$ for all $c\in\cC_\kappa^d$ and
$\phi\in I(\MM_\kappa^d)$. Such functionals depend on a chamber
only up to isometry. To construct a typical cell centred at the ``origin'' $p$, we use the
following recentering procedure.

Let $I(\MM_\kappa^d,p)$ be the subgroup of isometries fixing $p$.
It is compact and carries a unique Haar probability measure,
denoted by $\widetilde\lambda_\kappa(p,\cdot)$ and viewed as a
measure on $I(\MM_\kappa^d)$. For $x\in\MM_\kappa^d$, choose
an isometry $\phi_x$ with $\phi_x(p)=x$ and define
\begin{equation}\label{eq:RecenteringKernel}
\widetilde\lambda_\kappa(x,A)
:=
\int_{I(\MM_\kappa^d,p)}
{\bf1}\{\phi_x\psi\in A\}\,
\widetilde\lambda_\kappa(p,\dint\psi),
\qquad A\subset I(\MM_\kappa^d)\ \text{Borel}.
\end{equation}
Invariance makes \eqref{eq:RecenteringKernel} independent
of the choice of $\phi_x$. The resulting family is a measurable
probability kernel, and $\widetilde\lambda_\kappa(x,\cdot)$ is
supported on the isometries sending $p$ to $x$. Its inverse
isometries therefore move $x$ to $p$, with the remaining freedom
averaged over the isotropy group.

\begin{definition}\label{def:TypicalCell}
The \emph{typical cell} associated with $z$ is a
$\cC_\kappa^d$-valued random chamber
$C_{d,t,\kappa}^{\mathrm{typ},z}$ with distribution
\begin{equation}\label{eq:TypicalCellDefinition}
\PP(C_{d,t,\kappa}^{\mathrm{typ},z}\in A)
:=
\frac{1}{\gamma_{d,t,\kappa}^{\,z}\cH_\kappa^d(B)}
\EE\sum_{c\in\PHT_t}{\bf1}\{z(c)\in B\}
\int_{I(\MM_\kappa^d)}
{\bf1}\{\phi^{-1}c\in A\}\,
\widetilde\lambda_\kappa(z(c),\dint\phi),
\end{equation}
where $A\subset\cC_\kappa^d$ is Borel and
$B\subset\MM_\kappa^d$ is any Borel set of positive finite volume.
\end{definition}

\begin{remark}
For $\kappa<0$ and $t<t_{\rm crit}(d,\kappa)$, the tessellation
$\PHT_t$ contains unbounded cells of infinite volume almost surely.
By \cite[Lemma~2.2]{BuehlerGusakovaRecke} and the discussion
following it, there is no measurable isometry-covariant centre
assignment for all cells of this tessellation. Hence, our
definition of cell intensity and typical cell does not extend
beyond the regime \eqref{eq:AdmissibleRegime}.
\end{remark}

The right-hand side of \eqref{eq:TypicalCellDefinition} does not
depend on $B$. Indeed, for fixed $A$, its numerator defines an
isometry-invariant measure in $B$, dominated by
$\EE\xi_t^{\kappa,z}$. It is therefore a constant multiple of
$\cH_\kappa^d$. Taking $A=\cC_\kappa^d$ shows that
\eqref{eq:TypicalCellDefinition} is a probability measure.
This is the Palm mark distribution of the cell-centre process,
see \cite{RZ90,Last10} for background. By construction, we have
\[
z(C_{d,t,\kappa}^{\mathrm{typ},z})=p
\qquad\text{almost surely}.
\]
Also, the distribution of $C_{d,t,\kappa}^{\mathrm{typ},z}$ is invariant under $I(\MM_\kappa^d,p)$.
For an isometry-invariant functional $h$, we have
$h(\phi^{-1}c)=h(c)$ for every isometry $\phi$.
Since $\widetilde\lambda_\kappa(z(c),\cdot)$ is a probability
measure, \eqref{eq:TypicalCellDefinition} immediately
implies \eqref{eq:TypicalCellInvariantExpectations}.

The cell intensity and the distribution of the isometry class
of the typical cell do not depend on the choice of centre map,
as the following proposition shows.

\begin{proposition}
\label{prop:CentreIndependence}
Let $z_1$ and $z_2$ be centre maps. Under
\eqref{eq:AdmissibleRegime},
$
\gamma_{d,t,\kappa}^{\,z_1}
=
\gamma_{d,t,\kappa}^{\,z_2}.
$
Moreover, for every nonnegative measurable isometry-invariant
functional $h:\cC_\kappa^d\to[0,\infty]$,
\[
\EE h(C_{d,t,\kappa}^{\mathrm{typ},z_1})
=
\EE h(C_{d,t,\kappa}^{\mathrm{typ},z_2}).
\]
\end{proposition}

Proposition~\ref{prop:CentreIndependence} will be proved in
Section~\ref{sec:CellIntensityTypicalCell}, using the mass transport
principle of Lemma~\ref{lm:MassTransport}. Henceforth we fix one
centre map and suppress it from the notation, writing
$\xi_t^\kappa$, $\gamma_{d,t,\kappa}$ and
$C_{d,t,\kappa}^{\rm typ}$. The centred distribution of
$C_{d,t,\kappa}^{\rm typ}$ may still depend on this fixed choice.
Only its distribution up to isometry is independent of it.

We conclude with the relation between cell intensity and the mean
volume of the typical cell. Recall from
\eqref{eq:NonemptyProbability} that
\[
p_\kappa(t)
=
\begin{cases}
1-e^{-t/\sqrt{\kappa}},&\kappa>0,\\
1,&\kappa\leq0.
\end{cases}
\]

\begin{proposition}
\label{prop:TypicalCellVolumeIntensity}
Under \eqref{eq:AdmissibleRegime},
\begin{equation}\label{eq:VolumeViaIntensity}
\gamma_{d,t,\kappa}\,
\EE\cH_\kappa^d(C_{d,t,\kappa}^{\rm typ})
=
p_\kappa(t).
\end{equation}
In particular, $\EE\cH_\kappa^d(C_{d,t,\kappa}^{\rm typ})<\infty$.
\end{proposition}

\begin{proof}
Let $\lambda_\kappa$ be Haar measure on $I(\MM_\kappa^d)$
normalized by
\begin{equation}\label{eq:IsometryMeasureNormalization}
\int_{I(\MM_\kappa^d)}
{\bf1}\{\phi(p)\in A\}\,\lambda_\kappa(\dint\phi)
=
\cH_\kappa^d(A),
\qquad A\subset\MM_\kappa^d\ \text{Borel}.
\end{equation}
The isometry groups considered here are unimodular, so
$\lambda_\kappa$ is also invariant under inversion.
For every nonnegative measurable function
$f:I(\MM_\kappa^d)\times\cC_\kappa^d\to[0,\infty]$, the refined
Campbell formula states that
\begin{equation}\label{eq:Campbell1}
\EE\sum_{c\in\PHT_t}
\int_{I(\MM_\kappa^d)}
f(\phi,\phi^{-1}c)\,
\widetilde\lambda_\kappa(z(c),\dint\phi)
=
\gamma_{d,t,\kappa}\,
\EE\int_{I(\MM_\kappa^d)}
f(\phi,C_{d,t,\kappa}^{\rm typ})\,
\lambda_\kappa(\dint\phi),
\end{equation}
see \cite{RZ90,Last10}.

Apply \eqref{eq:Campbell1} with $f(\phi,c):={\bf1}\{\phi^{-1}(p)\in|c|\}$.
On its left-hand side, $f(\phi,\phi^{-1}c)={\bf1}\{p\in|c|\}$.
Almost surely, $p$ lies on no hyperplane of the process and
belongs to exactly one cell support when the tessellation is
nonempty. Hence the left-hand side of \eqref{eq:Campbell1} equals
$p_\kappa(t)$.

For each fixed chamber $c$, invariance of $\lambda_\kappa$ under
inversion and \eqref{eq:IsometryMeasureNormalization} give
\[
\int_{I(\MM_\kappa^d)}
{\bf1}\{\phi^{-1}(p)\in|c|\}\,
\lambda_\kappa(\dint\phi)
=
\cH_\kappa^d(c).
\]
Substitution into the right-hand side of \eqref{eq:Campbell1}
proves \eqref{eq:VolumeViaIntensity}.
Positivity and finiteness of the expected volume now follow from
Lemma~\ref{lem:CellIntensityFinite} and $0<p_\kappa(t)\leq1$.
\end{proof}

For $\kappa\leq0$, \eqref{eq:VolumeViaIntensity} says that the
mean volume of the typical cell is the reciprocal of the cell
intensity. For $\kappa>0$, the numerator is instead
$1-e^{-t/\sqrt{\kappa}}$, because the empty hyperplane process
contributes no cell.

\subsection{The zero cell and volume weighting}
\label{subsec:zero-cell-definition}

While the typical cell represents sampling with equal weight among cells,
the zero cell is defined by fixing a reference point and selecting the
cell containing it. This provides a different sampling principle, since
cells occupying more volume are more likely to contain a given point.
As we shall see below, the resulting distribution is closely related
to the volume-weighted typical-cell distribution.
Unlike the typical cell, the zero cell can be defined without assuming
that all cells are bounded. We therefore initially allow every $t>0$,
including intensities below the hyperbolic threshold in
\eqref{eq:AdmissibleRegime}. Almost surely, the fixed reference point
$p$ lies on no hyperplane of $\eta_t^\kappa$.

\begin{definition}\label{def:ZeroCell}
For $\kappa\leq0$, the \emph{zero cell} $C_{d,t,\kappa}^p$ is the
closure of the connected component of
\eqref{eq:HyperplaneComplement} containing $p$. Equivalently, it
is the intersection of the closed halfspaces bounded by the
hyperplanes of $\eta_t^\kappa$ and containing $p$.

For $\kappa>0$, on the event $\eta_t^\kappa\neq\emptyset$, the
zero cell is the unique elliptic cell whose support contains $p$,
with its inherited chamber structure. On the event
$\eta_t^\kappa=\emptyset$, we put $C_{d,t,\kappa}^p:=\MM_\kappa^d$,
regarded as a boundaryless whole-space object with no proper faces.
This extension does not add a cell to $\PHT_t$ or a centre to
$\xi_t^\kappa$.
\end{definition}

To specify the state space of the zero cell, we enlarge the space
$\cC_\kappa^d$ of compact chambers. For $\kappa\leq0$, let
$\cD_\kappa^d$ be the space of closed geodesically convex subsets
of $\MM_\kappa^d$ with nonempty interior, allowing unbounded sets.
We equip it with the $\sigma$-algebra generated by the sets
\[
\{D\in\cD_\kappa^d:D\cap U\neq\emptyset\},
\qquad U\subset\MM_\kappa^d\ \text{open}.
\]
For $\kappa>0$, put
\[
\cD_\kappa^d:=\cC_\kappa^d\sqcup\{\MM_\kappa^d\},
\]
with the disjoint-union $\sigma$-algebra. Here the added whole-space
object is distinct from every chamber in $\cC_\kappa^d$, even if
their supports coincide. Thus $C_{d,t,\kappa}^p$ is a
$\cD_\kappa^d$-valued random object for every $t>0$. Under
\eqref{eq:AdmissibleRegime}, it belongs to $\cC_\kappa^d$ almost
surely on the event $\eta_t^\kappa\neq\emptyset$.

For the remainder of this subsection, assume
\eqref{eq:AdmissibleRegime}. Recall that
\eqref{eq:VolumeViaIntensity} gives
\[
\EE\cH_\kappa^d(C_{d,t,\kappa}^{\rm typ})
=
\frac{p_\kappa(t)}{\gamma_{d,t,\kappa}}
\in(0,\infty).
\]
We may therefore weight the distribution of the typical cell by
its volume.

\begin{definition}
\label{def:VolumeWeightedTypicalCell}
The \emph{volume-weighted typical cell}
$C_{d,t,\kappa}^{\rm vol}$ is a $\cC_\kappa^d$-valued random
chamber with distribution
\begin{equation}\label{eq:def_vol_typ}
\PP(C_{d,t,\kappa}^{\rm vol}\in A)
:=
\frac{
\EE\!\left[
\cH_\kappa^d(C_{d,t,\kappa}^{\rm typ})
{\bf1}\{C_{d,t,\kappa}^{\rm typ}\in A\}
\right]
}{
\EE\cH_\kappa^d(C_{d,t,\kappa}^{\rm typ})
},
\qquad A\subset\cC_\kappa^d\ \text{Borel}.
\end{equation}
\end{definition}

The next proposition makes precise the relation between volume
weighting and sampling the cell containing a fixed point. The
conditioning in positive curvature is necessary because
$\PHT_t$ has no cells on the empty-process event, whereas the
zero cell was extended to that event in Definition~\ref{def:ZeroCell}.

\begin{proposition}
\label{prop:ZeroCellVolumeBias}
Under \eqref{eq:AdmissibleRegime}, for every nonnegative measurable
isometry-invariant functional
$h:\cD_\kappa^d\to[0,\infty]$,
\begin{equation}\label{eq:ZeroCellVolumeBias}
\EE h(C_{d,t,\kappa}^{\rm vol})
=
\EE\!\left[
h(C_{d,t,\kappa}^p)
\,\middle|\,
\eta_t^\kappa\neq\emptyset
\right].
\end{equation}
For $\kappa\leq0$, the conditioning can be omitted, since
$p_\kappa(t)=1$ by \eqref{eq:NonemptyProbability}.
\end{proposition}

\begin{proof}
Apply \eqref{eq:Campbell1} to the nonnegative measurable function
\[
f(\phi,c)
:=
{\bf1}\{\phi^{-1}(p)\in|c|\}h(c),
\qquad
(\phi,c)\in I(\MM_\kappa^d)\times\cC_\kappa^d.
\]
Since $h$ is isometry invariant,
\[
f(\phi,\phi^{-1}c)
=
{\bf1}\{p\in|c|\}h(c).
\]
Almost surely, exactly one cell support contains $p$ on the event
$\eta_t^\kappa\neq\emptyset$, and the sum on the left-hand side
of \eqref{eq:Campbell1} is zero otherwise. Its expectation is
therefore
\[
p_\kappa(t)\,
\EE\!\left[
h(C_{d,t,\kappa}^p)
\,\middle|\,
\eta_t^\kappa\neq\emptyset
\right].
\]
For each fixed chamber $c$, invariance of Haar measure under
inversion and \eqref{eq:IsometryMeasureNormalization} give
\[
\int_{I(\MM_\kappa^d)}
f(\phi,c)\,\lambda_\kappa(\dint\phi)
=
\cH_\kappa^d(c)h(c).
\]
Consequently, \eqref{eq:Campbell1} yields
\begin{equation*}
p_\kappa(t)\,
\EE\!\left[
h(C_{d,t,\kappa}^p)
\,\middle|\,
\eta_t^\kappa\neq\emptyset
\right]
=
\gamma_{d,t,\kappa}\,
\EE\!\left[
\cH_\kappa^d(C_{d,t,\kappa}^{\rm typ})
h(C_{d,t,\kappa}^{\rm typ})
\right].
\end{equation*}
Divide this expression by $p_\kappa(t)>0$
and use \eqref{eq:VolumeViaIntensity} and
\eqref{eq:def_vol_typ} to obtain \eqref{eq:ZeroCellVolumeBias}.
All integrands are nonnegative, so the argument also applies when
the expectations in \eqref{eq:ZeroCellVolumeBias} are infinite.
\end{proof}

Thus, conditional on the hyperplane process being nonempty,
observing the cell containing $p$ amounts, up to isometry, to
weighting the typical-cell distribution by volume. Equation
\eqref{eq:ZeroCellVolumeBias} does not identify the centred laws:
$C_{d,t,\kappa}^{\rm vol}$ has centre $p$ under the fixed centre
map, whereas $C_{d,t,\kappa}^p$ is selected by containing $p$,
not by having its centre there.

\section{Cell intensity and mean geometric characteristics of the typical cell}
\label{sec:CellIntensityTypicalCell}

\subsection{Mean \texorpdfstring{$\kappa$}{kappa}-intrinsic volumes of the typical cell}
\label{subsec:MeanIntrinsicVolumes}

Recall from Section~\ref{subsec:cells-faces} that $\cC_\kappa^d$
is the space of compact full-dimensional polyhedral chambers.
For $c\in\cC_\kappa^d$, the collection $\cF_k(c)$ consists of its
distinct $k$-faces. We first define the exterior angles
that enter the $\kappa$-intrinsic volumes of $c$.

\paragraph{Exterior angles: Nonpositive curvature.}
Suppose first that $\kappa\leq0$, let $F\in\cF_k(c)$ with
$0\leq k<d$, and choose $x\in\operatorname{relint}F$.
The Riemannian metric makes the tangent space $T_x\MM_\kappa^d$ a Euclidean vector
space. The normal space of $F$ at $x$ is
\[
N_xF:=(T_xF)^\perp\subset T_x\MM_\kappa^d.
\]
Let $n_1,\ldots,n_m\in N_xF$ be the outward unit normals at $x$
to the facets of $c$ containing $F$. The outer normal cone is
\begin{equation}\label{eq:ExteriorNormalCone}
N(c,F;x)
:=
\left\{\sum_{i=1}^m a_i n_i:a_i\geq0\right\}
\subset N_xF.
\end{equation}
The \emph{exterior angle} of $c$ at $F$ is its normalized solid
angle,
\begin{equation}\label{eq:ExteriorAngleDefinition}
\alpha(F,c)
:=
\frac{
\cH^{d-k-1}\bigl(N(c,F;x)\cap\SS(N_xF)\bigr)
}{\omega_{d-k}},
\end{equation}
where $\SS(N_xF)$ is the unit sphere in the $(d-k)$-dimensional
Euclidean space $N_xF$, and $\cH^{d-k-1}$ denotes its usual
spherical Hausdorff measure. For $k=d-1$, this measure is counting
measure on the two-point unit sphere. Parallel transport along
$\operatorname{relint}F$ identifies the normal cones isometrically,
since the supporting faces and facets are totally geodesic.
Thus \eqref{eq:ExteriorAngleDefinition} is independent of $x$.

\paragraph{Exterior angles: Positive curvature.}
For $\kappa>0$, let $c=\{\widehat c,-\widehat c\}$
be an elliptic cell and $F\in\cF_k(c)$. Choose a
spherical $k$-face $\widehat F$ of $\widehat c$
whose projective image is $F$.

We define the \emph{local incidences} of $F$ with $c$
as follows. If $\widehat F$ is not a whole great
$k$-sphere, then $F$ is said to have one local
incidence with $c$. If $\widehat F$ is a whole great
$k$-sphere, then $F$ is said to have two local
incidences with $c$, corresponding to the two
antipodal points of $\widehat F$ above each point
of $F$.

Let $\alpha(\widehat F,\widehat c)$ denote the
normalized exterior angle of $\widehat c$ at
$\widehat F$, computed on the spherical cover
using \eqref{eq:ExteriorNormalCone} and
\eqref{eq:ExteriorAngleDefinition}. Here the
normal cones are taken in the tangent spaces
of $\SS_\kappa^d$.

We assign this angle to each local incidence $\iota$,
denote it by $\alpha(F,c;\iota)$, and define
\begin{equation}\label{eq:EllipticExteriorAngle}
\alpha(F,c)
:=
\sum_\iota\alpha(F,c;\iota)
=
\begin{cases}
\alpha(\widehat F,\widehat c),
&\widehat F\text{ is not a whole great $k$-sphere},\\[0.4em]
2\alpha(\widehat F,\widehat c),
&\widehat F\text{ is a whole great $k$-sphere}.
\end{cases}
\end{equation}

\begin{example}[Local incidences]
Consider an arrangement of $N\geq1$ elliptic hyperplanes
$H_1,\ldots,H_N$ in general position in $\MM_\kappa^d$,
where $\kappa>0$. If $1\leq N\leq d$, these hyperplanes
have a common intersection
$
F_*:=H_1\cap\cdots\cap H_N,
$
which is a complete elliptic $(d-N)$-flat.
There are $2^{N-1}$ elliptic cells, and $F_*$ is
a face of each of them, with two local incidences.
Every other face has exactly one local incidence
with each elliptic cell containing it.
In particular, when $N=d$, the exceptional face
$F_*$ is a single point in elliptic space, whose
spherical lift consists of two antipodal points.
If $N\geq d+1$, every face has exactly one local
incidence with each elliptic cell containing it.
\end{example}

By summing exterior angles over local incidences,
we ensure that, for every face of the tessellation,
the exterior angles of all elliptic cells containing
that face sum to one, as shown in
Lemma~\ref{lm:ExteriorAnglePartition} below.

\paragraph{$\kappa$-intrinsic volumes.}
We now define the $\kappa$-intrinsic volumes of compact
full-dimensional polyhedral chambers in $\MM_\kappa^d$.

\begin{definition}
\label{def:KappaIntrinsicVolumes}
For $\kappa \in \RR$ and  $0\leq k\leq d$, the \emph{$k$-th $\kappa$-intrinsic volume}
is the functional $\upsilon_k^\kappa:\cC_\kappa^d\to[0,\infty)$
defined by
\begin{equation}\label{eq:ExteriorAngleFunctionals}
\upsilon_k^\kappa(c)
:=
\frac{1}{\omega_{k+1}}
\sum_{F\in\cF_k(c)}
\alpha(F,c)\cH_\kappa^k(F),
\end{equation}
where $\cF_d(c)=\{c\}$ and $\alpha(c,c):=1$.
\end{definition}

The functionals in \eqref{eq:ExteriorAngleFunctionals} are
nonnegative, measurable and isometry invariant. In particular,
\begin{equation}\label{eq:TopIntrinsicVolume}
\upsilon_d^\kappa(c)
=
\frac{\cH_\kappa^d(c)}{\omega_{d+1}}.
\end{equation}

\begin{remark}[$\kappa = 0$]
For Euclidean polytopes, the face-angle representation of the
usual intrinsic volumes $V_0,\ldots,V_d$ gives
\begin{equation}\label{eq:EuclideanIntrinsicVolumes}
\upsilon_k^0(c)=\frac{V_k(c)}{\omega_{k+1}},
\qquad 0\leq k\leq d.
\end{equation}
\end{remark}

\begin{remark}[Agreement of the \(\kappa\)-intrinsic volumes with the spherical face-angle sums]
For $\kappa>0$, the face-angle sum in
\eqref{eq:ExteriorAngleFunctionals} agrees with the
corresponding sum computed on either spherical representative
$\widehat c$. Indeed, let $\widehat F$ be a spherical
$k$-face of $\widehat c$ and let $F=q_\kappa(\widehat F)$.
If $\widehat F$ is not a whole great $k$-sphere, then
$\cH_\kappa^k(F)=\cH_\kappa^k(\widehat F)$ and  $\alpha(F,c)=\alpha(\widehat F,\widehat c)$.
If $\widehat F$ is a whole great $k$-sphere, then
$\cH_\kappa^k(F)=\frac12\cH_\kappa^k(\widehat F)$ and  $\alpha(F,c)=2\alpha(\widehat F,\widehat c)$.
Thus, in both cases,
\[
\alpha(F,c)\cH_\kappa^k(F)
=
\alpha(\widehat F,\widehat c)
\cH_\kappa^k(\widehat F).
\]
Summing over the faces proves the assertion.
\end{remark}

\begin{remark}[Steiner formula for $\kappa\leq0$]
\label{rem:SteinerFormula}
For $\kappa\leq0$, the functionals in
\eqref{eq:ExteriorAngleFunctionals} are also the normalized
coefficients in the volume formula for parallel sets. To state it,
define the generalized sine function
\begin{equation}\label{eq:GeneralizedSine}
\sn_\kappa(s)
:=
\begin{cases}
\sin(\sqrt{\kappa}\,s)/\sqrt{\kappa},&\kappa>0,\\
s,&\kappa=0,\\
\sinh(\sqrt{-\kappa}\,s)/\sqrt{-\kappa},&\kappa<0.
\end{cases}
\end{equation}
Recall that the generalized cosine $\cs_\kappa$ was defined in
Section~\ref{subsec:space-hyperplanes}. For $c\in\cC_\kappa^d$,
$\kappa\leq0$, and $r\geq0$, put
\[
c_r:=\{x\in\MM_\kappa^d:d_\kappa(x,c)\leq r\},
\qquad
d_\kappa(x,c):=\inf_{y\in c}d_\kappa(x,y).
\]
Then
\begin{equation}\label{eq:SteinerFormula}
\cH_\kappa^d(c_r)
=\cH_\kappa^d(c)
+\sum_{k=0}^{d-1}
\omega_{k+1}\omega_{d-k}\upsilon_k^\kappa(c)
\int_0^r
\sn_\kappa(s)^{d-k-1}\cs_\kappa(s)^k\,\dint s.
\end{equation}
Indeed, the normal exponential map over a $k$-face has Jacobian
$\sn_\kappa(s)^{d-k-1}\cs_\kappa(s)^k$. Integration over its unit
normal cone contributes $\omega_{d-k}\alpha(F,c)$. Every point
outside $c$ has a unique nearest point in $c$, so summing over
the faces proves \eqref{eq:SteinerFormula}, see \cite[Equation (1.11)]{Solanes}.

For $\kappa=0$, the integral in \eqref{eq:SteinerFormula} is
$r^{d-k}/(d-k)$. Together with
\eqref{eq:EuclideanIntrinsicVolumes}, this recovers the usual
Euclidean Steiner formula. For comparison with the smooth
formulation in terms of mean curvature integrals in constant
curvature, see \cite[Section~1.2, Equation~(1.11)]{Solanes}.
Equation~\eqref{eq:SteinerFormula} is asserted here only for
$\kappa\leq0$. In particular, it is not a parallel-volume formula
for the projective supports of arbitrary elliptic chambers.
\end{remark}

\paragraph{Mean $\kappa$-intrinsic volumes of the typical cell.}
We now state the formula for the expected $\kappa$-intrinsic volumes of the typical cell of $\PHT_t$. Recall from
Section~\ref{subsec:typical-cell} that $\gamma_{d,t,\kappa}$ is
the cell intensity defined by \eqref{eq:CellIntensityDefinition},
and that $C_{d,t,\kappa}^{\rm typ}$ denotes the typical cell
specified by \eqref{eq:TypicalCellDefinition}.

\begin{theorem}[Mean $\kappa$-intrinsic volumes of the typical cell]
\label{thm:MeanIntrinsicVolumes}
Under \eqref{eq:AdmissibleRegime}, all the following expectations
are finite, and for $0\leq k\leq d$,
\begin{equation}\label{eq:MeanExteriorAngleFunctionals}
\gamma_{d,t,\kappa}\,
\EE\upsilon_k^\kappa(C_{d,t,\kappa}^{\rm typ})
=
\begin{cases}
\displaystyle
\frac{1}{\omega_{d+1}}\frac{t^{d-k}}{(d-k)!},
&0\leq k<d,\\[0.8em]
\displaystyle
\frac{p_\kappa(t)}{\omega_{d+1}},&k=d,
\end{cases}
\end{equation}
where $p_\kappa(t)=\PP(\eta_t^\kappa\neq\emptyset)$ was defined
in \eqref{eq:NonemptyProbability}.
\end{theorem}

By \eqref{eq:TopIntrinsicVolume}, the case $k=d$ of
\eqref{eq:MeanExteriorAngleFunctionals} is precisely
\eqref{eq:VolumeViaIntensity}. The proof of
Theorem~\ref{thm:MeanIntrinsicVolumes} is given in the next two
subsections. We first compute the volume densities of the
$k$-skeletons and then use the exterior angles to distribute this
volume among incident cells by the mass transport principle.

\subsection{Skeleton densities}
\label{subsec:skeleton-densities}

To prove Theorem~\ref{thm:MeanIntrinsicVolumes}, we first determine
the expected $k$-dimensional volume of the tessellation's skeleton
in an observation window. Recall that $\cF_k(\PHT_t)$ consists
of the distinct $k$-faces of $\PHT_t$, with
$\cF_d(\PHT_t)=\PHT_t$. For a Borel set
$B\subset\MM_\kappa^d$, define the random measures
\begin{equation}\label{eq:SkeletonMeasures}
S_k(B)
:=
\sum_{F\in\cF_k(\PHT_t)}\cH_\kappa^k(F\cap B),
\qquad 0\leq k\leq d.
\end{equation}
In positive curvature, the intersections and Hausdorff measures in
\eqref{eq:SkeletonMeasures} refer to the geometric supports of the
faces. Each face is counted once.  In particular, $S_0$ counts vertices. For $k=d$, the cell
supports partition the space up to a set of volume zero whenever
the process is nonempty, and hence
\begin{equation}\label{eq:TopSkeletonMeasure}
S_d(B)
=
\begin{cases}
\cH_\kappa^d(B),&\eta_t^\kappa\neq\emptyset,\\
0,&\eta_t^\kappa=\emptyset.
\end{cases}
\end{equation}

For $k<d$, there is an equivalent description in terms of
intersections of hyperplanes. By general position, each $k$-face
almost surely has a unique set of $d-k$ supporting hyperplanes. Conversely,
the faces supported by a fixed $k$-dimensional intersection of
$d-k$ hyperplanes partition that flat up to a set of
$\cH_\kappa^k$-measure zero. Therefore, almost surely,
\begin{equation}\label{eq:SkeletonIntersectionRepresentation}
S_k(B)
=
\frac{1}{(d-k)!}
\sum_{\substack{
(H_1,\ldots,H_{d-k})\in(\eta_t^\kappa)_{\neq}^{d-k}\\
\dim(H_1\cap\cdots\cap H_{d-k})=k}}
\cH_\kappa^k(B\cap H_1\cap\cdots\cap H_{d-k}).
\end{equation}
Here $(\eta_t^\kappa)_{\neq}^{d-k}$ denotes the set of ordered
tuples of distinct hyperplanes, so the factor $1/(d-k)!$ removes
the ordering. Equation~\eqref{eq:SkeletonIntersectionRepresentation}
identifies $S_k(B)$ with the total $k$-volume in $B$ of the
intersection process of order $d-k$. 

The next lemma gives an explicit formula for
$\EE S_k(B)$, the expected $k$-dimensional volume
of the tessellation's $k$-skeleton in a Borel set $B$.
For $k=d$, the formula follows immediately from
\eqref{eq:TopSkeletonMeasure}. For $k<d$, the formula
in hyperbolic space of curvature $-1$ follows from
\cite[Theorem~1]{HeroldHugThaele}, after converting
the normalization of the hyperplane measure using
\eqref{eq:FlatMeasureNormalization}, see also
\cite[Remark~1]{HeroldHugThaele}.
Related integral-geometric calculations in spaces
of constant curvature can be found in
\cite{BetkenHugThaele}.
We include a proof in our notation to make explicit
the constants needed in the proof of
Theorem~\ref{thm:MeanIntrinsicVolumes}.

\begin{lemma}[$k$-skeleton densities]
	\label{lm:SkeletonIntensities}
Let $t>0$. For every Borel set
$B\subset\MM_\kappa^d$ and every $0\leq k\leq d$,
\begin{equation}\label{eq:SkeletonIntensities}
\begin{gathered}
\EE S_k(B)
=
\sigma_{d,k,t,\kappa}\,\cH_\kappa^d(B),\qquad
\sigma_{d,k,t,\kappa}
=
\begin{cases}
\displaystyle
\frac{\omega_{k+1}}{\omega_{d+1}}
\frac{t^{d-k}}{(d-k)!},&0\leq k<d,\\[0.8em]
p_\kappa(t),&k=d,
\end{cases}
\end{gathered}
\end{equation}
where $p_\kappa(t)=\PP(\eta_t^\kappa\neq\emptyset)$ was defined
in \eqref{eq:NonemptyProbability}.
\end{lemma}

The constant $\sigma_{d,k,t,\kappa}$ is called the
\emph{$k$-skeleton density}. It is the expected $k$-dimensional
face volume per unit ambient volume. For $k=0$, it is the vertex
intensity.  In particular, \eqref{eq:SkeletonIntensities}
asserts that each measure $\EE S_k$ is locally finite.

To compute $\EE S_k(B)$ using
\eqref{eq:SkeletonIntersectionRepresentation}, we need two
integral-geometric identities.
Recall that $A_\kappa(d,k)$ is
the space of totally geodesic $k$-flats and that
$\widetilde\mu_k^\kappa$ is the invariant measure normalized in
\eqref{eq:FlatMeasureNormalization}.
Lemma~\ref{lm:CroftonTypeFormula} gives Crofton's formula
for mean section volumes, while
Lemma~\ref{lm:IntersectionFormula} allows us to replace
integration over hyperplane tuples by integration over
their intersection flats.

\begin{lemma}[Crofton formula]\label{lm:CroftonTypeFormula}
Let $0\leq i\leq k\leq d-1$, and let
$W\subset\MM_\kappa^d$ be a countably
$(d+i-k)$-rectifiable Borel set. (Recall that a set is countably $m$-rectifiable if,
up to a set of zero $m$-dimensional Hausdorff measure,
it can be covered by countably many Lipschitz images
of subsets of $\RR^m$.) Then, as an identity in
$[0,\infty]$,
\begin{equation}\label{eq:Formula_2}
\int_{A_\kappa(d,k)}
\cH_\kappa^i(W\cap E)\,
\widetilde\mu_k^\kappa(\dint E)
=
\frac{\omega_{i+1}}{\omega_{d-k+i+1}}
\cH_\kappa^{d+i-k}(W).
\end{equation}
\end{lemma}

For the Crofton formula in standard spaces of constant curvature,
see \cite[Lemma~2.1]{BetkenHugThaele}. The hyperbolic case is also
stated in \cite[Lemma~2]{HeroldHugThaele}. Multiplying the
classically normalized measures there by the factor in
\eqref{eq:FlatMeasureNormalization} gives the ratio of sphere
areas in \eqref{eq:Formula_2}. In positive curvature, apply the
spherical formula to the antipodal preimage $q_\kappa^{-1}(W)$.
Under the identification of elliptic flats with great subspheres,
the normalized invariant flat measures agree, and both Hausdorff
measures in the formula are multiplied by two on passing to the
cover. Dividing by two yields \eqref{eq:Formula_2} on elliptic
space. Rescaling as in \eqref{eq:HyperplaneMeasureScale} gives
the formula for arbitrary nonzero curvature. The two sides acquire
the same factor $|\kappa|^{-(d+i-k)/2}$.
In particular, taking $i=k$ in \eqref{eq:Formula_2} gives
\begin{equation}\label{eq:CroftonSpecialCase}
\int_{A_\kappa(d,k)}
\cH_\kappa^k(B\cap E)\,
\widetilde\mu_k^\kappa(\dint E)
=
\frac{\omega_{k+1}}{\omega_{d+1}}
\cH_\kappa^d(B),
\qquad 0\leq k<d,
\end{equation}
for every Borel set $B\subset\MM_\kappa^d$. The choice $i=0$,
$k=d-1$ in \eqref{eq:Formula_2} also recovers the segment formula
\eqref{eq:HyperplaneSegmentCrofton}, with constant
$\omega_1/\omega_2=1/\pi$.

\begin{lemma}[Intersection formula]\label{lm:IntersectionFormula}
Let $0\leq k<d$. For every nonnegative measurable function
$f:A_\kappa(d,k)\to[0,\infty]$,
\begin{multline}\label{eq:Formula_1}
\int_{A_\kappa(d,d-1)^{d-k}}
f(H_1\cap\cdots\cap H_{d-k})
\, {\bf1}\{\dim(H_1\cap\cdots\cap H_{d-k})=k\}
\prod_{j=1}^{d-k}\widetilde\mu_{d-1}^\kappa(\dint H_j)
\\
=\int_{A_\kappa(d,k)}f(E)\,
\widetilde\mu_k^\kappa(\dint E).
\end{multline}
The integrand on the left is defined to be zero whenever the
intersection is empty or has dimension different from $k$.
\end{lemma}

The normalization in \eqref{eq:FlatMeasureNormalization} removes
any additional dimensional constant from
\eqref{eq:Formula_1}. For comparison, the hyperbolic formula with
$\mu_k^{-1}$ in place of $\widetilde\mu_k^{-1}$ is given in
\cite[Lemma~4]{HeroldHugThaele}.

\begin{proof}[Proof of Lemma \ref{lm:IntersectionFormula}]
Let $\Theta$ be the image of
$(\widetilde\mu_{d-1}^\kappa)^{\otimes(d-k)}$,
restricted to tuples with $k$-dimensional intersection,
under the map
\[
(H_1,\ldots,H_{d-k})
\longmapsto H_1\cap\cdots\cap H_{d-k}.
\]
Then $\Theta$ is an isometry-invariant Borel measure on
$A_\kappa(d,k)$. It is also locally finite. Namely, if an intersection meets
a compact ball $B$, every hyperplane in its defining tuple meets
$B$. Consequently,
\[
\Theta\bigl(\{E\in A_\kappa(d,k):E\cap B\neq\emptyset\}\bigr)
\leq
\left[
\widetilde\mu_{d-1}^\kappa
\bigl(\{H\in A_\kappa(d,d-1):H\cap B\neq\emptyset\}\bigr)
\right]^{d-k}
<\infty.
\]
Given $E\in A_\kappa(d,k)$, choose a point $x\in E$
and a compact ball $B$ centred at $x$ with positive radius.
Every flat sufficiently close to $E$ also meets $B$.
The preceding estimate therefore shows that $\Theta$
is locally finite. By uniqueness of the isometry-invariant
Radon measure, we have
$\Theta=a\widetilde\mu_k^\kappa$ for some $a\geq0$.

To determine $a$, choose a small geodesic ball $W$ with positive
finite volume. Successively apply \eqref{eq:Formula_2} to the
hyperplanes in the defining tuple, using the section dimension
at each step. Nonempty nontransverse intersections form a null
set, and all integrands are nonnegative. Thus Tonelli's theorem
gives
\begin{equation}\label{eq:IntersectionNormalization}
\int_{A_\kappa(d,k)}\cH_\kappa^k(W\cap E)\,\Theta(\dint E)
=
\left(\prod_{\ell=k+1}^{d}
\frac{\omega_\ell}{\omega_{\ell+1}}\right)
\cH_\kappa^d(W)
=
\frac{\omega_{k+1}}{\omega_{d+1}}
\cH_\kappa^d(W).
\end{equation}
By \eqref{eq:CroftonSpecialCase}, the same integral with
$\widetilde\mu_k^\kappa$ in place of $\Theta$ has the value
in \eqref{eq:IntersectionNormalization}. Since this value is
positive and finite, we conclude that $a=1$, proving \eqref{eq:Formula_1}.
\end{proof}

\begin{proof}[Proof of Lemma~\ref{lm:SkeletonIntensities}]
For $k=d$, taking expectations in \eqref{eq:TopSkeletonMeasure}
gives \eqref{eq:SkeletonIntensities}.

Let $0\leq k<d$. Apply the multivariate Mecke formula
\cite[Theorem~4.4]{LPbook} to
\eqref{eq:SkeletonIntersectionRepresentation}. The intensity
measure of $\eta_t^\kappa$ is
$t\widetilde\mu_{d-1}^\kappa$. Using
\eqref{eq:Formula_1} to integrate over the intersection flat,
and then \eqref{eq:CroftonSpecialCase}, yields
\begin{align*}
\EE S_k(B)
=
\frac{t^{d-k}}{(d-k)!}
\int_{A_\kappa(d,k)}\cH_\kappa^k(B\cap E)\,
\widetilde\mu_k^\kappa(\dint E)
=
\frac{\omega_{k+1}}{\omega_{d+1}}
\frac{t^{d-k}}{(d-k)!}\,\cH_\kappa^d(B).
\end{align*}
All interchanges are justified by nonnegativity.
\end{proof}

For $k<d$, the constants $\sigma_{d,k,t,\kappa}$ in
\eqref{eq:SkeletonIntensities} are independent of the curvature.
Only the top-dimensional density $\sigma_{d,d,t,\kappa}=p_\kappa(t)$
contains the positive-curvature empty-process correction.

\subsection{From skeleton densities to mean \texorpdfstring{$\kappa$}{kappa}-intrinsic volumes}
\label{subsec:intrinsic-volumes-proof}

We now prove Theorem~\ref{thm:MeanIntrinsicVolumes}
by using exterior angles to distribute the $k$-skeleton
measures $S_k$ among incident cells and applying the
mass transport principle. Recall from
\eqref{eq:SkeletonMeasures} and \eqref{eq:SkeletonIntensities}
that $S_k$ is the $k$-volume measure on the $k$-face
skeleton and that
\[
\EE S_k(B)
=
\sigma_{d,k,t,\kappa}\cH_\kappa^d(B),
\qquad 0\leq k\leq d.
\]

\paragraph{The angle partition.}
Recall that $\alpha(F,c)$ denotes the normalized exterior angle
in \eqref{eq:ExteriorAngleDefinition}, with the convention~\eqref{eq:EllipticExteriorAngle} when $\kappa>0$.
The following lemma shows that, for every $k$-face of the
tessellation, the exterior angles of all cells containing
that face sum to one.

\begin{lemma}\label{lm:ExteriorAnglePartition}
Almost surely, for every $0\leq k<d$ and every
$F\in\cF_k(\PHT_t)$,
\begin{equation}\label{eq:ExteriorAnglePartition}
\sum_{\substack{c\in\PHT_t\\F\in\cF_k(c)}}\alpha(F,c)=1.
\end{equation}
\end{lemma}

\begin{proof}
Choose a point $x\in F$ lying on no hyperplanes of
$\eta_t^\kappa$ other than those supporting $F$.
By general position, exactly $d-k$ hyperplanes of
$\eta_t^\kappa$ pass through $x$, and their unit normals
$n_1,\ldots,n_{d-k}$ form a basis of the normal space
$N_xF=(T_xF)^\perp$.
All $2^{d-k}$ choices of sides occur
as local chamber sectors at $x$. The corresponding outer normal
cones are
\[
\left\{\sum_{i=1}^{d-k} a_i\varepsilon_i n_i:a_i\geq0\right\},
\qquad
(\varepsilon_1,\ldots,\varepsilon_{d-k})\in\{-1,1\}^{d-k}.
\]
These cones are the images of the coordinate orthants under the
linear isomorphism sending the standard orthonormal basis of $\RR^{d-k}$ to
$n_1,\ldots,n_{d-k}$. Moreover, they cover $N_xF$ with disjoint interiors,
and their boundaries have spherical measure zero. Their
normalized solid angles therefore sum to one.

In positive curvature, two of these local contributions may
come from the same elliptic cell. They are then combined
in $\alpha(F,c)$ by \eqref{eq:EllipticExteriorAngle}.
Thus all $2^{d-k}$ contributions are counted exactly once,
proving \eqref{eq:ExteriorAnglePartition}.
\end{proof}

\paragraph{Mass transport.}
The mass transport principle, stated in the next lemma,
allows us to relate different mean geometric characteristics
of the tessellation. Roughly speaking, it asserts that, for
an isometry-invariant mass transport, the mean mass sent
from a region equals the mean mass received by it.
An early application of the mass transport principle to invariant
percolation on regular trees appears in \cite{Haggstrom97}. The
method was developed systematically for group-invariant percolation
on graphs in \cite{BLPS}. A continuous version for the hyperbolic
plane is given in \cite[Theorem~5.2]{BenjaminiSchramm}, see also
\cite[Theorem~2.1]{BuehlerGusakovaRecke} for hyperbolic space of
arbitrary dimension. We use the homogeneous-space formulation of
\cite[Corollary~6.2]{GentnerLast}.

\begin{lemma}[Mass transport principle]
	\label{lm:MassTransport}
Let $M$ be a nonnegative Borel measure on
$\MM_\kappa^d\times\MM_\kappa^d$, invariant under the diagonal
action of $I(\MM_\kappa^d)$, meaning that $M(\phi A\times\phi B)=M(A\times B)$
for all Borel sets $A,B\subset\MM_\kappa^d$ and all
$\phi\in I(\MM_\kappa^d)$. Suppose that the first marginal
$A\mapsto M(A\times\MM_\kappa^d)$ is locally finite.
Then there exists a number $a\in[0,\infty)$ such that
\begin{equation}\label{eq:MassTransport}
M(B\times\MM_\kappa^d)
=
M(\MM_\kappa^d\times B)
=
a\cH_\kappa^d(B)
\end{equation}
for every Borel set $B\subset\MM_\kappa^d$.
\end{lemma}

\begin{proof}
The isometry group $I(\MM_\kappa^d)$ acts transitively and properly on
$\MM_\kappa^d$. It is unimodular in each of the three
geometries. In positive curvature this follows from compactness.
Local finiteness of the first marginal implies that $M$ is
$\sigma$-finite.
Hence \cite[Corollary~6.2]{GentnerLast}, with the modular function
equal to one, gives equality of the two marginals. The first
marginal is a locally finite isometry-invariant measure and is
therefore $a\cH_\kappa^d$ for some finite $a\geq0$.
This proves \eqref{eq:MassTransport}.
\end{proof}

\begin{remark}
For $\kappa>0$, the conclusion can also be seen directly.
Indeed, compactness makes $M$ a finite measure, both marginals are invariant
multiples of $\cH_\kappa^d$, and their total masses coincide.
In particular, passing from the sphere to elliptic space causes
no difficulty for Lemma~\ref{lm:MassTransport}.
\end{remark}

We first apply Lemma~\ref{lm:MassTransport} to establish the
independence of the centre map stated in
Proposition~\ref{prop:CentreIndependence}. We recall its statement
before giving the proof.

\begingroup
\theoremstyle{plain}
\newtheorem*{centreindependencerecall}
  {Proposition~\ref{prop:CentreIndependence}}
\begin{centreindependencerecall}
Let $z_1,z_2:\cC_\kappa^d\to\MM_\kappa^d$ be centre maps.
Under \eqref{eq:AdmissibleRegime}, we have $\gamma_{d,t,\kappa}^{\,z_1}=\gamma_{d,t,\kappa}^{\,z_2}$.
Moreover, for every nonnegative measurable isometry-invariant
functional $h:\cC_\kappa^d\to[0,\infty]$,
\[
\EE h(C_{d,t,\kappa}^{\mathrm{typ},z_1})
=
\EE h(C_{d,t,\kappa}^{\mathrm{typ},z_2}).
\]
\end{centreindependencerecall}
\endgroup

\begin{proof}
Transport one unit from $z_1(c)$ to $z_2(c)$ for every cell $c$.
The expected transported mass is the nonnegative measure defined by
\[
M(A\times B)
:=
\EE\sum_{c\in\PHT_t}
{\bf1}\{z_1(c)\in A\}{\bf1}\{z_2(c)\in B\}
\]
for Borel sets $A,B\subset\MM_\kappa^d$.
Covariance of the centre maps makes $M$ diagonally invariant.
Its first marginal is
$\gamma_{d,t,\kappa}^{\,z_1}\cH_\kappa^d$, which is locally
finite by Lemma~\ref{lem:CellIntensityFinite}. Its second
marginal is $\gamma_{d,t,\kappa}^{\,z_2}\cH_\kappa^d$.
Equation~\eqref{eq:MassTransport} therefore gives equality of
the two intensities.

For the assertion about typical cells, fix $h$ as in the
statement and replace the unit mass above by
$h_n(c):=\min\{h(c),n\}$, where $n\in\NN$. This transport
is again diagonally invariant, and its first marginal is bounded
by $n\gamma_{d,t,\kappa}^{\,z_1}\cH_\kappa^d$.
Combining \eqref{eq:MassTransport} with
\eqref{eq:TypicalCellInvariantExpectations} yields
\[
\gamma_{d,t,\kappa}^{\,z_1}
\EE h_n(C_{d,t,\kappa}^{\mathrm{typ},z_1})
=
\gamma_{d,t,\kappa}^{\,z_2}
\EE h_n(C_{d,t,\kappa}^{\mathrm{typ},z_2}).
\]
The intensities are equal and positive. Dividing by their common
value, letting $n\to\infty$ and using monotone convergence,
proves the assertion.
\end{proof}

\paragraph{Proof of the mean intrinsic-volume formula.}
We now have all the ingredients needed to compute the mean
$\kappa$-intrinsic volumes of the typical cell in terms of
the cell intensity. We fix a centre map $z$ and use the
notation $\gamma_{d,t,\kappa}$ and
$C_{d,t,\kappa}^{\rm typ}$ from
Section~\ref{subsec:typical-cell}.

\begin{proof}[Proof of Theorem~\ref{thm:MeanIntrinsicVolumes}]
First, let $0\leq k<d$, and define a nonnegative Borel measure $M_k$ on
$\MM_\kappa^d\times\MM_\kappa^d$ by
\begin{equation}\label{eq:AngleMassTransport}
M_k(A\times B)
:=
\frac{1}{\omega_{k+1}}
\EE\sum_{c\in\PHT_t}\sum_{F\in\cF_k(c)}
\alpha(F,c)\cH_\kappa^k(F\cap A)
{\bf1}\{z(c)\in B\}.
\end{equation}
Thus mass is sent from points of a face to the centres of its
incident cells, with the exterior angles as weights. For elliptic
faces, $F\cap A$ refers to their projective supports, while
$\alpha(F,c)$ includes all local incidences as in
\eqref{eq:EllipticExteriorAngle}. The measure $M_k$ is invariant under the diagonal
action of $I(\MM_\kappa^d)$ because the Poisson hyperplane tessellation is isometry invariant and the
centre map is covariant.

By \eqref{eq:ExteriorAnglePartition} and
\eqref{eq:SkeletonIntensities}, the first marginal of $M_k$  is
\begin{equation}\label{eq:AngleTransportFirstMarginal}
M_k(A\times\MM_\kappa^d)
=
\frac{1}{\omega_{k+1}}\EE S_k(A)
=
\frac{\sigma_{d,k,t,\kappa}}{\omega_{k+1}}
\cH_\kappa^d(A).
\end{equation}
This measure is locally finite on $\MM_\kappa^d$, so Lemma~\ref{lm:MassTransport} applies.
On the other hand, \eqref{eq:ExteriorAngleFunctionals} and
\eqref{eq:TypicalCellInvariantExpectations} give, for every
Borel set $B\subset\MM_\kappa^d$ of positive finite volume,
\begin{equation}\label{eq:AngleTransportSecondMarginal}
M_k(\MM_\kappa^d\times B)
=
\EE\sum_{c\in\PHT_t}{\bf1}\{z(c)\in B\}
\upsilon_k^\kappa(c)
=
\gamma_{d,t,\kappa}\cH_\kappa^d(B)
\EE\upsilon_k^\kappa(C_{d,t,\kappa}^{\rm typ}).
\end{equation}
Comparing \eqref{eq:AngleTransportFirstMarginal} and
\eqref{eq:AngleTransportSecondMarginal} by
\eqref{eq:MassTransport}, and dividing by
$\cH_\kappa^d(B)$, yields
\begin{equation}\label{eq:IntrinsicVolumesViaSkeletons}
\gamma_{d,t,\kappa}
\EE\upsilon_k^\kappa(C_{d,t,\kappa}^{\rm typ})
=
\frac{\sigma_{d,k,t,\kappa}}{\omega_{k+1}}.
\end{equation}
All sums and integrals used above are nonnegative.
The expectation in \eqref{eq:IntrinsicVolumesViaSkeletons}
is finite because $\sigma_{d,k,t,\kappa}<\infty$ and
$\gamma_{d,t,\kappa}>0$.

For $k=d$, \eqref{eq:TopIntrinsicVolume},
\eqref{eq:VolumeViaIntensity} and
\eqref{eq:SkeletonIntensities} give directly
\[
\gamma_{d,t,\kappa}
\EE\upsilon_d^\kappa(C_{d,t,\kappa}^{\rm typ})
=
\frac{p_\kappa(t)}{\omega_{d+1}}
=
\frac{\sigma_{d,d,t,\kappa}}{\omega_{d+1}}.
\]
Thus \eqref{eq:IntrinsicVolumesViaSkeletons} holds also for
$k=d$. Substituting the constants from \eqref{eq:SkeletonIntensities}
into \eqref{eq:IntrinsicVolumesViaSkeletons} proves
\eqref{eq:MeanExteriorAngleFunctionals}, including finiteness
of all the expectations.
\end{proof}

\paragraph{Extension to general invariant hyperplane tessellations.}
The relation between skeleton densities and mean
$\kappa$-intrinsic volumes of the typical cell is not specific
to Poisson hyperplane tessellations.

Let $\mathcal H$ be an isometry-invariant, locally finite
random collection of totally geodesic hyperplanes in
$\MM_\kappa^d$, and let $\mathcal T$ be the hyperplane
tessellation induced by $\mathcal H$. Assume that almost
surely all cells are compact and that the hyperplanes are
in general position, in the following sense: for every
$1\leq m\leq d$, every nonempty intersection of $m$
distinct hyperplanes has dimension $d-m$, and no $d+1$
distinct hyperplanes have a common point. In positive
curvature, we use the conventions for elliptic chambers
and their faces from Subsection~\ref{subsec:cells-faces}
and the definition of local incidences from
Subsection~\ref{subsec:MeanIntrinsicVolumes}.

\begin{proposition}[Mean $\kappa$-intrinsic volumes through skeleton densities]
\label{prop:GeneralHyperplaneIntrinsicVolumes}
Under the preceding assumptions, suppose that $\mathcal T$
has a finite positive cell intensity $\gamma$, and fix
$0\leq k<d$. If its $k$-skeleton has finite volume density
$\sigma_k$, that is, if
\[
\EE\sum_{F\in\cF_k(\mathcal T)}
\cH_\kappa^k(F\cap B)
=
\sigma_k\,\cH_\kappa^d(B)
\]
for every Borel set $B\subset\MM_\kappa^d$ of finite
volume, then the typical cell $C^{\rm typ}$ of
$\mathcal T$ satisfies
\begin{equation}\label{eq:GeneralHyperplaneIntrinsicVolumes}
\gamma\,
\EE\upsilon_k^\kappa(C^{\rm typ})
=
\frac{\sigma_k}{\omega_{k+1}}.
\end{equation}
\end{proposition}

\begin{proof}
By general position, every $k$-face is contained in
exactly $d-k$ hyperplanes. Their normal vectors are
linearly independent, and the exterior normal cones
of the incident local chambers partition the normal
space. Consequently, the exterior angles of all
cells containing a given $k$-face sum to one, with
multiple local incidences of an elliptic cell combined
according to \eqref{eq:EllipticExteriorAngle}.

We can therefore apply the mass transport argument
from the proof of Theorem~\ref{thm:MeanIntrinsicVolumes}.
The first marginal of the corresponding transport
has density $\sigma_k/\omega_{k+1}$, while its second
marginal has density
$\gamma\EE\upsilon_k^\kappa(C^{\rm typ})$.
Lemma~\ref{lm:MassTransport} gives
\eqref{eq:GeneralHyperplaneIntrinsicVolumes}.
The expectation is finite because $\sigma_k<\infty$
and $\gamma>0$.
\end{proof}

Thus the Poisson assumption is not needed for
\eqref{eq:GeneralHyperplaneIntrinsicVolumes}.
In the present paper, it is used to evaluate the
skeleton densities explicitly in
\eqref{eq:SkeletonIntensities}.

\paragraph{Internal-angle functionals.}
We now introduce analogues of the $\kappa$-intrinsic volumes
in which exterior angles are replaced by internal angles.
The argument used in the proof of
Theorem~\ref{thm:MeanIntrinsicVolumes} applies to these
functionals as well.
Let $c\in\cC_\kappa^d$ be a compact full-dimensional
polyhedral chamber, let $F\in\cF_k(c)$ with
$0\leq k<d$, and choose $x\in\operatorname{relint}F$.
For one local incidence $\iota$ of $F$ with $c$, let
$n_1,\ldots,n_m\in N_xF$ be the outward unit normals
to all facets containing this incidence. The corresponding
transverse tangent cone is
\[
D(c,F;x;\iota)
:=
\{v\in N_xF:\langle v,n_j\rangle\leq0,\
1\leq j\leq m\}.
\]
This is the tangent cone of that local chamber sector after
removing the directions tangent to $F$. Its normalized solid
angle is the \emph{internal angle}
\begin{equation*}
\beta(F,c;\iota)
:=
\frac{\cH^{d-k-1}
\bigl(D(c,F;x;\iota)\cap\SS(N_xF)\bigr)}{\omega_{d-k}}.
\end{equation*}
As for \eqref{eq:ExteriorAngleDefinition}, parallel transport
shows that this value is independent of $x$. For $\kappa\leq0$
there is only one local incidence. For $\kappa>0$ we sum over
local incidences using one spherical representative, exactly as
in \eqref{eq:EllipticExteriorAngle}. In both cases, write
$\beta(F,c)$ for the resulting angle and set $\beta(c,c):=1$.

Define the nonnegative measurable isometry-invariant functionals
$\iota_k^\kappa:\cC_\kappa^d\to[0,\infty)$ by
\begin{equation*}
\iota_k^\kappa(c)
:=
\frac{1}{\omega_{k+1}}
\sum_{F\in\cF_k(c)}\beta(F,c)\cH_\kappa^k(F),
\qquad 0\leq k\leq d.
\end{equation*}
In particular,
$\iota_d^\kappa(c)=\cH_\kappa^d(c)/\omega_{d+1}$.
For $k=0$, the functional $\iota_0^\kappa(c)$ is one half
of the sum of the normalized internal solid angles at
the vertices of $c$.

\begin{proposition}
\label{prop:MeanInternalAngleFunctionals}
Under \eqref{eq:AdmissibleRegime}, all the following expectations
are finite, and for $0\leq k\leq d$,
\begin{equation}\label{eq:MeanInternalAngleFunctionals}
\gamma_{d,t,\kappa}
\EE\iota_k^\kappa(C_{d,t,\kappa}^{\rm typ})
=
\frac{\sigma_{d,k,t,\kappa}}{\omega_{k+1}}
=
\begin{cases}
\displaystyle
\frac{1}{\omega_{d+1}}\frac{t^{d-k}}{(d-k)!},&0\leq k<d,\\[0.8em]
\displaystyle
\frac{p_\kappa(t)}{\omega_{d+1}},&k=d,
\end{cases}
\end{equation}
where $p_\kappa(t)$ was defined in
\eqref{eq:NonemptyProbability}. Consequently,
\begin{equation}\label{eq:EqualInternalExteriorMeans}
\EE\iota_k^\kappa(C_{d,t,\kappa}^{\rm typ})
=
\EE\upsilon_k^\kappa(C_{d,t,\kappa}^{\rm typ}),
\qquad 0\leq k\leq d.
\end{equation}
\end{proposition}

\begin{proof}
At a relative interior point of a $k$-face, the transverse
tangent cones of all local chamber sectors partition $N_xF$
up to their boundaries. Indeed, these cones are obtained by
choosing the signs of the $d-k$ linear inequalities given by
the supporting hyperplanes. Their normalized solid angles
therefore sum to one.
Thus, we obtain the following analogue of
Lemma~\ref{lm:ExteriorAnglePartition} for internal angles:
\begin{equation}\label{eq:InternalAnglePartition}
\sum_{\substack{c\in\PHT_t\\F\in\cF_k(c)}}\beta(F,c)=1,
\qquad F\in\cF_k(\PHT_t),\quad 0\leq k<d.
\end{equation}
The proof of Theorem~\ref{thm:MeanIntrinsicVolumes}
now applies upon replacing $\alpha(F,c)$ by $\beta(F,c)$
in \eqref{eq:AngleMassTransport}.
By \eqref{eq:InternalAnglePartition}, the first marginal remains
\eqref{eq:AngleTransportFirstMarginal}, while the second is
\eqref{eq:AngleTransportSecondMarginal} with $\iota_k^\kappa$
in place of $\upsilon_k^\kappa$. Lemma~\ref{lm:MassTransport}
therefore proves \eqref{eq:MeanInternalAngleFunctionals} for
$k<d$. For $k=d$, the two functionals both equal
$\cH_\kappa^d(c)/\omega_{d+1}$, and
\eqref{eq:VolumeViaIntensity} applies. Comparing
\eqref{eq:MeanInternalAngleFunctionals} with
\eqref{eq:MeanExteriorAngleFunctionals} proves
\eqref{eq:EqualInternalExteriorMeans}.
\end{proof}

\subsection{Gauss--Bonnet formula and cell intensity}
\label{subsec:cell-intensity-proof}

We now prove Theorem~\ref{thm:cell-intensity}. Recall that
$\gamma_{d,t,\kappa}$ is the cell intensity defined in
\eqref{eq:CellIntensityDefinition}, and that
$C_{d,t,\kappa}^{\rm typ}$ is the typical cell. The preceding
subsection established \eqref{eq:MeanExteriorAngleFunctionals},
which determines the products
$\gamma_{d,t,\kappa}\EE\upsilon_k^\kappa(C_{d,t,\kappa}^{\rm typ})$
for all $0\leq k\leq d$. To determine the cell intensity itself,
we use the following polyhedral \textit{Gauss--Bonnet formula}.

\begin{lemma}[Polyhedral Gauss--Bonnet formula]
\label{lm:PolyhedralGaussBonnet}
Let $c\in\cC_\kappa^d$ be a compact full-dimensional
polyhedral chamber. Then
\begin{equation}\label{eq:PolyhedralGaussBonnet}
2\sum_{j=0}^{\lfloor d/2\rfloor}
\kappa^j\upsilon_{2j}^\kappa(c)=1.
\end{equation}
The coefficient with $j=0$ is understood to be one
also when $\kappa=0$.
\end{lemma}

\begin{proof}
For a compact $d$-dimensional polyhedron $P$ in Euclidean,
hyperbolic or spherical space of curvature $\kappa$, the
Gauss--Bonnet theorem for polyhedral manifolds with totally
geodesic faces, \cite[Theorem~1.1]{Moroianu}, gives
\begin{equation}\label{eq:GaussBonnetFaceSum}
\frac{\chi(P)}{2}
=
\sum_{j=0}^{\lfloor d/2\rfloor}
\frac{\kappa^j}{\omega_{2j+1}}
\sum_{F\in\cF_{2j}(P)}
\alpha(F,P)\cH_\kappa^{2j}(F).
\end{equation}
Here $\chi(P)$ is the Euler characteristic of $P$, the faces and angles
are taken in $P$ itself, and the exterior angles are normalized
as in \eqref{eq:ExteriorAngleDefinition}. The top-dimensional
face, when present in the sum, has angle $\alpha(P,P)=1$.

For $\kappa\leq0$, apply \eqref{eq:GaussBonnetFaceSum} to
$P=c$. The compact convex polytope $c$ is contractible, so
$\chi(c)=1$. Equation~\eqref{eq:ExteriorAngleFunctionals}
therefore gives \eqref{eq:PolyhedralGaussBonnet}.
For $\kappa>0$, choose one spherical representative $\widehat c$
of $c$. It is a compact polyhedral manifold with totally geodesic
faces. To see that $\chi(\widehat c)=1$, regard $\widehat c$ as
a subset of the unit sphere and choose $u$ in its spherical
interior. Since $\widehat c$ is defined by a nonempty collection
of hemisphere inequalities, $-u\notin\widehat c$. The maps
\[
y\longmapsto
\frac{(1-s)y+su}{\|(1-s)y+su\|},
\qquad y\in\widehat c,\quad 0\leq s\leq1,
\]
are therefore well defined and give a contraction of $\widehat c$
to $u$. Indeed, each hemisphere inequality is preserved by the convex
combination and by normalization. Now, apply \eqref{eq:GaussBonnetFaceSum} to $P=\widehat c$.
By the incidence convention in \eqref{eq:EllipticExteriorAngle},
the face-angle sums on $\widehat c$ agree with those defining
$\upsilon_k^\kappa(c)$ in \eqref{eq:ExteriorAngleFunctionals}.
Thus \eqref{eq:PolyhedralGaussBonnet} follows from
$\chi(\widehat c)=1$. In this argument the Gauss--Bonnet formula is applied
to the spherical representative $\widehat c$, not to the projective support
$|c|$.
\end{proof}

\begin{proof}[Proof of Theorem~\ref{thm:cell-intensity}]
The idea is to apply the polyhedral Gauss--Bonnet formula
\eqref{eq:PolyhedralGaussBonnet} to the typical cell and
take expectations. Theorem~\ref{thm:MeanIntrinsicVolumes}
then allows us to evaluate the resulting terms.
Under \eqref{eq:AdmissibleRegime}, the typical cell belongs to
$\cC_\kappa^d$ almost surely. By
Theorem~\ref{thm:MeanIntrinsicVolumes},
\[
\EE\left[
\sum_{j=0}^{\lfloor d/2\rfloor}
2|\kappa|^j\upsilon_{2j}^\kappa
(C_{d,t,\kappa}^{\rm typ})
\right]<\infty.
\]
We may therefore apply \eqref{eq:PolyhedralGaussBonnet} to the
typical cell and take expectations term by term, including the case
$\kappa<0$ and $t=t_{\rm crit}(d,\kappa)$. Multiplying by
$\gamma_{d,t,\kappa}$ and using
\eqref{eq:MeanExteriorAngleFunctionals} together with
\eqref{eq:SkeletonIntensities} yields
\begin{equation}\label{eq:StructuralCellIntensity}
\gamma_{d,t,\kappa}
=
2\sum_{j=0}^{\lfloor d/2\rfloor}
\kappa^j\gamma_{d,t,\kappa}
\EE\upsilon_{2j}^\kappa(C_{d,t,\kappa}^{\rm typ})
=
2\sum_{j=0}^{\lfloor d/2\rfloor}
\frac{\kappa^j}{\omega_{2j+1}}
\sigma_{d,2j,t,\kappa}.
\end{equation}
For $2j<d$, \eqref{eq:SkeletonIntensities} gives
\[
\frac{\sigma_{d,2j,t,\kappa}}{\omega_{2j+1}}
=
\frac{1}{\omega_{d+1}}
\frac{t^{d-2j}}{(d-2j)!}.
\]
If $d$ is even, \eqref{eq:StructuralCellIntensity} also contains
the term with $2j=d$, which is
\[
\frac{2\kappa^{d/2}}{\omega_{d+1}}
\sigma_{d,d,t,\kappa}
=
\frac{2\kappa^{d/2}}{\omega_{d+1}}p_\kappa(t).
\]
Recall from \eqref{eq:NonemptyProbability} that
$p_\kappa(t)=1$ for $\kappa\leq0$ and
$p_\kappa(t)=1-e^{-t/\sqrt{\kappa}}$ for $\kappa>0$.
Substitution into \eqref{eq:StructuralCellIntensity} therefore
gives
\begin{equation}\label{eq:CellIntensityAllCurvatures}
\gamma_{d,t,\kappa}
=
\frac{2}{\omega_{d+1}}
\sum_{j=0}^{\lfloor d/2\rfloor}
\kappa^j\frac{t^{d-2j}}{(d-2j)!}
-
\begin{cases}
\displaystyle
\frac{2\kappa^{d/2}}{\omega_{d+1}}
e^{-t/\sqrt{\kappa}},
&\kappa>0\text{ and }d\text{ is even},\\[0.8em]
0,&\text{otherwise}.
\end{cases}
\end{equation}
This proves Theorem~\ref{thm:cell-intensity}. The exponential
correction occurs only in even dimensions because only then does
the sum in \eqref{eq:StructuralCellIntensity} contain the
top-dimensional density $\sigma_{d,d,t,\kappa}=p_\kappa(t)$.
\end{proof}

\begin{remark}[Mean volume of the typical cell]
Combining \eqref{eq:CellIntensityAllCurvatures} with
\eqref{eq:VolumeViaIntensity} also determines the expected volume:
\[
\EE\cH_\kappa^d(C_{d,t,\kappa}^{\rm typ})
=
\frac{p_\kappa(t)}{\gamma_{d,t,\kappa}}.
\]
In particular, for $\kappa<0$ this is the reciprocal of the
intensity in \eqref{eq:CellIntensityAllCurvatures}, including at
$t=t_{\rm crit}(d,\kappa)$. Its positivity and finiteness were
already established in
Proposition~\ref{prop:TypicalCellVolumeIntensity}.
\end{remark}

\subsection{Mean Crofton functionals of the typical cell}
\label{subsec:crofton-functionals}

We next determine the mean Crofton functionals of the typical cell.
These functionals measure the sets of totally geodesic flats meeting
a chamber. We first state their mean values in terms of
lower-dimensional cell intensities, and then prove the relation to
the $\kappa$-intrinsic volumes that yields these values.

\begin{definition}\label{def:CroftonFunctionals}
For a nonempty closed set $K\subset\MM_\kappa^d$ and
$0\leq m\leq d$, the \emph{$m$-th Crofton functional}
of $K$ is defined by
\begin{equation}\label{eq:CroftonFunctionalDefinition}
U_m^\kappa(K)
:=
\frac12\int_{A_\kappa(d,d-m)}
{\bf1}\{K\cap E\neq\emptyset\}\,
\widetilde\mu_{d-m}^\kappa(\dint E).
\end{equation}
For a chamber $c$, we put $U_m^\kappa(c):=U_m^\kappa(|c|)$,
where $|c|$ denotes its support. We also set
$U_m^\kappa(\emptyset):=0$.
\end{definition}

The integral in \eqref{eq:CroftonFunctionalDefinition} may be
infinite for unbounded sets. On Euclidean convex bodies, Crofton functionals
are the quermassintegrals in our normalization. By
\eqref{eq:FlatMeasureNormalization},
$\widetilde\mu_0^\kappa=2\cH_\kappa^d/\omega_{d+1}$, where we
 also recall  that
$\widetilde\mu_d^\kappa=\delta_{\MM_\kappa^d}$. Consequently, for every nonempty closed set
$K\subset\MM_\kappa^d$,
\begin{equation}\label{eq:CroftonEndpoints}
	U_0^\kappa(K)=\frac12,
	\qquad
	U_d^\kappa(K)=\frac{\cH_\kappa^d(K)}{\omega_{d+1}}.
\end{equation}
For Euclidean convex bodies, the usual Crofton formula gives
\begin{equation}\label{eq:EuclideanCroftonFunctionals}
U_m^0(K)=\frac{V_m(K)}{\omega_{m+1}},
\qquad 0\leq m\leq d,
\end{equation}
see \cite{SW}. In particular, for Euclidean polytopes,
\eqref{eq:EuclideanIntrinsicVolumes} and
\eqref{eq:EuclideanCroftonFunctionals} give
$U_m^0=\upsilon_m^0$.

For $1\leq j<d$, let $\gamma_{j,t,\kappa}$ denote the cell
intensity of the Poisson hyperplane tessellation in
$\MM_\kappa^j$, with the same parameter $t$ and with the invariant
hyperplane measure normalized intrinsically in dimension $j$
as in Section~\ref{subsec:space-hyperplanes}. In positive
curvature we retain the convention that an empty hyperplane
process contributes no cell. In dimension one,
$\MM_\kappa^1$ is a geodesic line for $\kappa\leq0$, and for
$\kappa>0$ it is the elliptic line
$\SS_\kappa^1/\{\pm{\rm id}\}$, a circle of length
$\pi/\sqrt{\kappa}$. Its hyperplanes are points, and their
intensity measure is $t\cH_\kappa^1/\pi$. Thus
\begin{equation}\label{eq:OneDimensionalCellIntensity}
\gamma_{1,t,\kappa}=\frac{t}{\pi},
\qquad \kappa\in\RR.
\end{equation}
On a line the interval intensity equals the point intensity.
On the elliptic line the number of cells equals the number of
points, including the empty realization under our convention.
For $\kappa<0$, assumption \eqref{eq:AdmissibleRegime} implies
\[
t\geq\pi(d-1)\sqrt{-\kappa}
>\pi(j-1)\sqrt{-\kappa},
\qquad 1\leq j<d.
\]
Hence all lower-dimensional tessellations used here have bounded
cells almost surely.

Recall that $\gamma_{d,t,\kappa}$ is the cell intensity defined
in \eqref{eq:CellIntensityDefinition}, and that
$C_{d,t,\kappa}^{\rm typ}$ is the typical cell defined in
\eqref{eq:TypicalCellDefinition}.

\begin{theorem}[Mean Crofton functionals of the typical cell]
\label{thm:MeanTypicalCellCrofton}
Under \eqref{eq:AdmissibleRegime}, all the mean Crofton
functionals of the typical cell are finite. For
$1\leq j\leq d-1$,
\begin{equation}\label{eq:MeanTypicalCellCrofton}
2\EE U_{d-j}^\kappa(C_{d,t,\kappa}^{\rm typ})
=
\frac{\omega_{j+1}}{\omega_{d+1}}
\frac{\gamma_{j,t,\kappa}}{\gamma_{d,t,\kappa}}.
\end{equation}
The remaining two mean values are
\begin{equation}\label{eq:MeanCroftonEndpoints}
\EE U_0^\kappa(C_{d,t,\kappa}^{\rm typ})=\frac12,
\qquad
\EE U_d^\kappa(C_{d,t,\kappa}^{\rm typ})
=
\frac{p_\kappa(t)}{\omega_{d+1}\gamma_{d,t,\kappa}},
\end{equation}
where $p_\kappa(t)=\PP(\eta_t^\kappa\neq\emptyset)$ was defined
in \eqref{eq:NonemptyProbability}.
\end{theorem}

The intensities on the right-hand side of
\eqref{eq:MeanTypicalCellCrofton} are explicit by
\eqref{eq:CellIntensityAllCurvatures} and
\eqref{eq:OneDimensionalCellIntensity}.
To prove Theorem~\ref{thm:MeanTypicalCellCrofton},
we first establish the relation between Crofton functionals
and $\kappa$-intrinsic volumes stated in
Proposition~\ref{prop:CroftonIntrinsicTransform}.
We begin with a sectional identity needed for the case
of nonpositive curvature. For a related formula involving
mean curvature integrals, see
\cite[Proposition~2.2.4]{Solanes}.
We give a direct proof of the sectional identity to
establish the normalization in
\eqref{eq:SectionExteriorFunctionals}.

\begin{lemma}
\label{lem:SectionalCrofton}
Let $d\geq2$, $\kappa\leq0$, and let $c\subset\MM_\kappa^d$
be a compact geodesically convex polytope with nonempty interior.
For $1\leq m\leq d-1$ and $0\leq j\leq d-m$,
\begin{equation}\label{eq:SectionExteriorFunctionals}
\int_{A_\kappa(d,d-m)}
\upsilon_j^{\kappa,E}(c\cap E)\,
\widetilde\mu_{d-m}^\kappa(\dint E)
=
\upsilon_{j+m}^\kappa(c),
\end{equation}
where $\upsilon_j^{\kappa,E}$ is computed intrinsically in $E$
using \eqref{eq:ExteriorAngleFunctionals} in dimension $d-m$.
\end{lemma}
\begin{proof}
For $j=d-m$, the assertion follows directly from
\eqref{eq:CroftonSpecialCase} and
$\upsilon_{d-m}^{\kappa,E}(c\cap E)
=\cH_\kappa^{d-m}(c\cap E)/\omega_{d-m+1}$.
We prove the hyperplane case for $0\leq j\leq d-2$ and then proceed by iteration.

Fix $F\in\cF_{j+1}(c)$ and take
$x\in\operatorname{relint}F$. Write $T:=T_xF$, $V:=T^\perp$,
and let $N\subset V$ be the outer normal cone of $c$ at $F$.
For a unit vector $u\in T_x\MM_\kappa^d$, let $H_{x,u}$ be the
totally geodesic hyperplane through $x$ with normal $u$.

We shall use the following incidence formula. To state it, define the incidence space
\[
\mathcal I_\kappa
:=
\{(x,H)\in\MM_\kappa^d\times A_\kappa(d,d-1):x\in H\}.
\]
Then, for every nonnegative measurable function
$g:\mathcal I_\kappa\to[0,\infty]$,
\begin{equation}\label{eq:SectionIncidenceCoarea}
	\begin{split}
		&\int_{A_\kappa(d,d-1)}\int_{F\cap H}
		g(x,H)\,\cH_\kappa^j(\dint x)\,
		\widetilde\mu_{d-1}^\kappa(\dint H)=
		\frac{\omega_d}{\omega_{d+1}}
		\int_F
		\EE_u\!\left[
		\|\proj_{T_xF}u\|\,g(x,H_{x,u})
		\right]\cH_\kappa^{j+1}(\dint x),
	\end{split}
\end{equation}
where $\EE_u$ denotes averaging with respect to the uniform
probability measure on the unit sphere of $T_x\MM_\kappa^d$.
To explain the Jacobian factor, observe that a hyperplane with
unit normal $u$ is locally a level set of a smooth function whose
gradient at $x$ is $u$. The restriction of this function to $F$
has gradient $\proj_{T_xF}u$. The coarea formula on $F$ therefore
contributes the factor $\|\proj_{T_xF}u\|$. Finally, the constant
$\omega_d/\omega_{d+1}$ follows from the normalization
\eqref{eq:FlatMeasureGeneralCurv} and
\eqref{eq:FlatMeasureNormalization}, using the reference point
$x$ in the parametrization of hyperplanes.

For a unit vector $u\in T_x\MM_\kappa^d$ with
$\|\proj_Tu\|>0$, consider the linear map $A_u:=\proj_{u^\perp}|_V$.
Since $u$ has a nonzero component in $T$, this map is an
isomorphism from $V$ onto the normal space of
$F\cap H_{x,u}$ within $H_{x,u}$. The outer normal cone
of the section $c\cap H_{x,u}$ at $F\cap H_{x,u}$ is
therefore $A_uN$. As before, we write $\alpha$ for the
normalized solid angle of a cone in its linear span.
We claim that
\[
\EE_u\left[
\alpha(A_uN)\,\middle|\,\|\proj_Tu\|
\right]
=\alpha(N).
\]
To prove this, let $R$ be an orthogonal transformation of
$V$, extended by the identity on $T$. The orthogonal
projection commutes with this transformation in the sense
that $A_{Ru}N=R A_uR^{-1}N$.
Since normalized solid angles are invariant under
orthogonal transformations, it follows that $\alpha(A_{Ru}N)=\alpha(A_uR^{-1}N)$.
Moreover, the conditional distribution of $u$ given
$\|\proj_Tu\|$ is invariant under $u\mapsto Ru$.

Now fix $u$ with $\|\proj_Tu\|>0$ and average over a
Haar-uniform orthogonal transformation $R$ of $V$.
For every unit vector $w\in A_uV$, the direction of
$RA_u^{-1}w$ is uniform on the unit sphere of $V$.
Consequently, averaging the normalized solid angle of
$A_uR^{-1}N$ over $R$ gives $\alpha(N)$.
Fubini's theorem and the preceding invariance of the
conditional distribution establish the claim. Finally,
$\|\proj_Tu\|=0$ occurs with probability zero.

Integration with respect to the uniform distribution on the unit sphere
of $T_x\MM_\kappa^d$ gives
\[
\begin{split}
	\frac{\omega_d}{\omega_{d+1}}
	\EE_u\|\proj_Tu\|
	&=
	\frac{\omega_d}{\omega_{d+1}}
	\frac{\Gamma((j+2)/2)\Gamma(d/2)}
	{\Gamma((j+1)/2)\Gamma((d+1)/2)}
	=\frac{\omega_{j+1}}{\omega_{j+2}}.
\end{split}
\]
Indeed, since $\dim T=j+1$, the squared length of the orthogonal
projection of a uniform unit vector onto $T$ has a beta distribution
with parameters $(j+1)/2$ and $(d-j-1)/2$, see \cite[Section~2.8]{KabluchkoSteigenbergerThaele2026}.

Now, we apply \eqref{eq:SectionIncidenceCoarea} with
$g(x,H)=\alpha(F\cap H,c\cap H)$, setting the integrand to zero
on nontransverse incidences. The preceding two identities yield
\begin{equation}\label{eq:FacewiseSectionCrofton}
\begin{split}
&\int_{A_\kappa(d,d-1)}
\alpha(F\cap H,c\cap H)\,
\cH_\kappa^j(F\cap H)\,
\widetilde\mu_{d-1}^\kappa(\dint H)
=
\frac{\omega_{j+1}}{\omega_{j+2}}
\alpha(F,c)\cH_\kappa^{j+1}(F).
\end{split}
\end{equation}
For almost every $H$, the $j$-faces of $c\cap H$ are precisely
the nonempty transverse sections of the $(j+1)$-faces of $c$.
Summing \eqref{eq:FacewiseSectionCrofton} over $F$ and dividing
by $\omega_{j+1}$ therefore gives
\begin{equation}\label{eq:OneSectionExterior}
\int_{A_\kappa(d,d-1)}
\upsilon_j^{\kappa,H}(c\cap H)\,
\widetilde\mu_{d-1}^\kappa(\dint H)
=
\upsilon_{j+1}^\kappa(c).
\end{equation}
Together with the top-dimensional case already considered, this
holds for every $0\leq j\leq d-1$.

For general $m$, iterate \eqref{eq:OneSectionExterior} inside
successive intersection flats. The iterated integral
can then be written over $m$ ambient hyperplanes. Equation
\eqref{eq:Formula_1} then replaces their product measure by
$\widetilde\mu_{d-m}^\kappa$, giving
\eqref{eq:SectionExteriorFunctionals}. All integrands are
nonnegative, so Tonelli's theorem justifies the iteration.
\end{proof}

\begin{proposition}[Crofton functionals and $\kappa$-intrinsic volumes]
\label{prop:CroftonIntrinsicTransform}
For every $c\in\cC_\kappa^d$ and $0\leq m\leq d$,
\begin{equation}\label{eq:CroftonIntrinsicTransform}
U_m^\kappa(c)
=
\sum_{r=0}^{\lfloor(d-m)/2\rfloor}
\kappa^r\upsilon_{m+2r}^\kappa(c).
\end{equation}
The coefficient with $r=0$ is understood to be one also when
$\kappa=0$.
\end{proposition}

\begin{proof}
For $m=0$, the assertion follows from
\eqref{eq:CroftonEndpoints} and
\eqref{eq:PolyhedralGaussBonnet}. For $m=d$, it follows from
\eqref{eq:CroftonEndpoints} and \eqref{eq:TopIntrinsicVolume}.
We therefore assume $1\leq m\leq d-1$.

Suppose first that $\kappa\leq0$. For almost every
$E\in A_\kappa(d,d-m)$, the section $c\cap E$ is empty or is
a compact convex polytope with nonempty relative interior in $E$.
Applying \eqref{eq:PolyhedralGaussBonnet} intrinsically in $E$
gives
\[
\frac12{\bf1}\{c\cap E\neq\emptyset\}
=
\sum_{r=0}^{\lfloor(d-m)/2\rfloor}
\kappa^r\upsilon_{2r}^{\kappa,E}(c\cap E),
\]
with empty sections contributing zero. When $d-m=1$, the same
identity follows directly from
$\upsilon_0^{\kappa,E}(c\cap E)=1/2$ for a nondegenerate segment.
Each term is integrable by Lemma~\ref{lem:SectionalCrofton}, so
we may integrate the finite signed sum term by term.
Equations \eqref{eq:CroftonFunctionalDefinition} and
\eqref{eq:SectionExteriorFunctionals} then give
\[
U_m^\kappa(c)
=
\sum_{r=0}^{\lfloor(d-m)/2\rfloor}
\kappa^r\upsilon_{m+2r}^\kappa(c),
\]
which proves \eqref{eq:CroftonIntrinsicTransform} for
$\kappa\leq0$.

Now suppose that $\kappa>0$. Choose one spherical representative
$\widehat c$, regarded as a subset of the unit sphere in
$\RR^{d+1}$, and let
\[
Q:=\{ru:r\geq0,\ u\in\widehat c\}
\]
be the cone over it. This is a full-dimensional polyhedral cone
contained in a closed halfspace, and therefore it is not a linear
subspace. Write $v_i(Q)$, $0\leq i\leq d+1$, for its conic
intrinsic volumes. Equivalently, $v_i(Q)$ is the probability that
the Euclidean nearest-point projection of a standard Gaussian
vector onto $Q$ belongs to the relative interior of an
$i$-dimensional face, see \cite[Section~2.2]{AmelunxenLotz}.
The face-angle representation of these quantities gives
\begin{equation}\label{eq:ConicExteriorAngles}
v_{k+1}(Q)=\kappa^{k/2}\upsilon_k^\kappa(c),
\qquad 0\leq k\leq d.
\end{equation}
Indeed, the contribution of a $(k+1)$-dimensional cone face is
the volume of its trace on the unit sphere divided by
$\omega_{k+1}$, multiplied by its normalized outer normal angle.
The factor $\kappa^{k/2}$ converts the volume of a face
on $\SS_\kappa^d$ to its volume on the unit sphere. The incidence
convention in \eqref{eq:EllipticExteriorAngle} identifies this
spherical face-angle sum with
\eqref{eq:ExteriorAngleFunctionals}. This includes linear cone
faces, whose spherical traces are whole subspheres with their
full Hausdorff measures.

Let $L$ be a uniformly distributed $(d-m+1)$-dimensional linear
subspace of $\RR^{d+1}$. The conic Crofton formula
\cite[Corollary~5.2]{AmelunxenLotz} gives
\begin{equation}\label{eq:ConicCroftonFormula}
\PP(Q\cap L\neq\{0\})
=
2\sum_{r=0}^{\lfloor(d-m)/2\rfloor}v_{m+1+2r}(Q).
\end{equation}
The projective flat corresponding to $L$ meets $|c|$ exactly
when $Q\cap L\neq\{0\}$. Its distribution is the normalized
invariant probability measure on $A_\kappa(d,d-m)$.
By \eqref{eq:FlatMeasureScalePositive}, the total mass of
$\widetilde\mu_{d-m}^\kappa$ is $\kappa^{-m/2}$. Consequently,
\eqref{eq:CroftonFunctionalDefinition},
\eqref{eq:ConicCroftonFormula} and
\eqref{eq:ConicExteriorAngles} yield
\[
U_m^\kappa(c)
=\frac{\kappa^{-m/2}}{2}\PP(Q\cap L\neq\{0\})
=\kappa^{-m/2}
\sum_{r=0}^{\lfloor(d-m)/2\rfloor}v_{m+1+2r}(Q)
=\sum_{r=0}^{\lfloor(d-m)/2\rfloor}
\kappa^r\upsilon_{m+2r}^\kappa(c).
\]
The cited conic Crofton formula only excludes linear subspaces, not cones
containing lines. Thus this argument applies to every chamber in
$\cC_\kappa^d$.
\end{proof}

\begin{proof}[Proof of Theorem~\ref{thm:MeanTypicalCellCrofton}]
Theorem~\ref{thm:MeanIntrinsicVolumes} gives finiteness of the
expectation of each term in
\eqref{eq:CroftonIntrinsicTransform}, also after replacing
$\kappa^r$ by $|\kappa|^r$. We may therefore take expectations
term by term, including at the hyperbolic critical intensity.
For $1\leq j\leq d-1$, use
\eqref{eq:CroftonIntrinsicTransform} with $m=d-j$ and substitute
\eqref{eq:MeanExteriorAngleFunctionals}. This gives
\begin{equation}\label{eq:MeanCroftonExpansion}
\omega_{d+1}\gamma_{d,t,\kappa}
\EE U_{d-j}^\kappa(C_{d,t,\kappa}^{\rm typ})
=
\sum_{r=0}^{\lfloor(j-1)/2\rfloor}
\kappa^r\frac{t^{j-2r}}{(j-2r)!}
+
\begin{cases}
\kappa^{j/2}p_\kappa(t),&j\text{ even},\\[0.3em]
0,&j\text{ odd}.
\end{cases}
\end{equation}
The separate term for even $j$ comes from the top-dimensional
intrinsic volume $\upsilon_d^\kappa$, whose mean involves
$p_\kappa(t)$ in \eqref{eq:MeanExteriorAngleFunctionals}.
By \eqref{eq:CellIntensityAllCurvatures} for $j\geq2$
and \eqref{eq:OneDimensionalCellIntensity} for $j=1$,
the right-hand side of \eqref{eq:MeanCroftonExpansion}
equals $\omega_{j+1}\gamma_{j,t,\kappa}/2$ for every
$1\leq j\leq d-1$. Dividing by
$\omega_{d+1}\gamma_{d,t,\kappa}/2$ proves
\eqref{eq:MeanTypicalCellCrofton}.
Finally, \eqref{eq:MeanCroftonEndpoints} follows from
\eqref{eq:CroftonEndpoints} and \eqref{eq:VolumeViaIntensity}.
\end{proof}

\begin{remark}
\label{rem:CroftonWholeEllipticSpace}
For $\kappa>0$, on the event $\eta_t^\kappa=\emptyset$,
the zero cell is defined to be the boundaryless whole
elliptic space $\MM_\kappa^d$, which does not belong to
$\cC_\kappa^d$. Hence,
Proposition~\ref{prop:CroftonIntrinsicTransform} does not
apply to this object. Its Crofton functionals are
nevertheless defined by \eqref{eq:CroftonFunctionalDefinition}.
Using \eqref{eq:FlatMeasureScalePositive} gives
\begin{equation}\label{eq:CroftonWholeEllipticSpace}
U_m^\kappa(\MM_\kappa^d)
=\frac{1}{2\kappa^{m/2}},
\qquad 0\leq m\leq d.
\end{equation}
Since Crofton functionals depend only on the support,
these values also agree with those of an elliptic
chamber represented by a pair of antipodal hemispheres.
\end{remark}

\subsection{Further consequences and examples}
\label{subsec:further-consequences}

We record several consequences of the cell-intensity formula
\eqref{eq:CellIntensityAllCurvatures} and the mean-value
formulas \eqref{eq:VolumeViaIntensity},
\eqref{eq:MeanExteriorAngleFunctionals} and
\eqref{eq:MeanTypicalCellCrofton}.
First we illustrate the cell-intensity formula and the
corresponding mean-volume identity in hyperbolic
spaces of dimensions $d=2,3,4$.

\paragraph{Low-dimensional hyperbolic examples.}
For $\kappa<0$, the first three cases of
\eqref{eq:CellIntensityAllCurvatures} are
\begin{equation}\label{eq:LowDimensionalHyperbolicIntensities}
\gamma_{2,t,\kappa}
=\frac{t^2+2\kappa}{4\pi},
\qquad
\gamma_{3,t,\kappa}
=\frac{t^3+6\kappa t}{6\pi^2},
\qquad
\gamma_{4,t,\kappa}
=\frac{t^4+12\kappa t^2+24\kappa^2}{32\pi^2}.
\end{equation}
Since $p_\kappa(t)=1$, \eqref{eq:VolumeViaIntensity} gives
\begin{equation*}
\EE\cH_\kappa^2(C_{2,t,\kappa}^{\rm typ})
=\frac{4\pi}{t^2+2\kappa},
\qquad
\EE\cH_\kappa^3(C_{3,t,\kappa}^{\rm typ})
=\frac{6\pi^2}{t^3+6\kappa t},
\qquad
\EE\cH_\kappa^4(C_{4,t,\kappa}^{\rm typ})
=\frac{32\pi^2}{t^4+12\kappa t^2+24\kappa^2}.
\end{equation*}
Each formula is valid for
$t\geq t_{\rm crit}(d,\kappa)=\pi(d-1)\sqrt{-\kappa}$ in its
respective dimension. In particular, all three expected volumes
are finite and positive at the corresponding critical intensity,
as also follows from
Proposition~\ref{prop:TypicalCellVolumeIntensity}.

\begin{remark}[Subcritical regime]
\label{rem:SubcriticalPolynomial}
The polynomial expressions in
\eqref{eq:LowDimensionalHyperbolicIntensities} cannot be interpreted
as cell intensities for all $t>0$. For example, writing
$\kappa=-a^2$ with $a>0$, the expression in dimension two is
\[
\frac{t^2-2a^2}{4\pi}<0,
\qquad 0<t<\sqrt{2}\,a,
\]
whereas an intensity is nonnegative. The sign-change value
$\sqrt{2}\,a$ is strictly below
$t_{\rm crit}(2,-a^2)=\pi a$.
\end{remark}

\paragraph{Dependence on curvature.}
We next examine the dependence of the cell intensity on
curvature, beginning with its behaviour under rescaling.

\begin{proposition}
\label{prop:IntensityScaling}
Under \eqref{eq:AdmissibleRegime}, for $\kappa\neq0$,
\begin{equation}\label{eq:CurvatureEffectIntensity}
\gamma_{d,t,\kappa}
=
|\kappa|^{d/2}
\gamma_{d,t/\sqrt{|\kappa|},\,\operatorname{sgn}(\kappa)}.
\end{equation}
\end{proposition}

\begin{proof}
Put $\kappa_0:=\operatorname{sgn}(\kappa)\in\{-1,1\}$.
Under the identification of the underlying spaces in
\eqref{eq:Volumescale}, we have
$g_\kappa=|\kappa|^{-1}g_{\kappa_0}$, and hence $\cH_\kappa^d=|\kappa|^{-d/2}\cH_{\kappa_0}^d$.
Moreover, \eqref{eq:HyperplaneMeasureScale} gives
\[
t\widetilde\mu_{d-1}^\kappa
=
\frac{t}{\sqrt{|\kappa|}}
\widetilde\mu_{d-1}^{\kappa_0}.
\]
Thus, under this identification, the hyperplane process
$\eta_t^\kappa$ has the same distribution as
$\eta_{t/\sqrt{|\kappa|}}^{\kappa_0}$.
Since rescaling preserves the isometry group and the cells,
we may use the same centre map for both tessellations.
Comparing their cell intensities using
\eqref{eq:CellIntensityDefinition} proves
\eqref{eq:CurvatureEffectIntensity}.
For $\kappa<0$, the rescaled parameter satisfies
$t/\sqrt{-\kappa}\geq\pi(d-1)$, so the intensity on
the right-hand side is defined under
\eqref{eq:AdmissibleRegime}.
\end{proof}

We now show that, for fixed $d$ and $t$, the cell intensity
and several mean geometric characteristics of the typical cell are $C^\infty$
functions of the curvature $\kappa$ in a neighbourhood of
zero. However, these functions need not be real analytic.

\begin{proposition}[Smooth dependence on curvature]
\label{prop:CurvatureSmoothness}
Fix $d\geq2$ and $t>0$. There exists $\varepsilon>0$ such that
\eqref{eq:AdmissibleRegime} holds for every
$\kappa\in(-\varepsilon,\varepsilon)$ and each of the functions
\[
\begin{aligned}
	\kappa &\longmapsto \gamma_{d,t,\kappa},
	&\qquad
	\kappa &\longmapsto
	\EE\cH_\kappa^d(C_{d,t,\kappa}^{\rm typ}),\\[0.6em]
	\kappa &\longmapsto
	\EE\upsilon_k^\kappa(C_{d,t,\kappa}^{\rm typ}),
	\quad 0\leq k\leq d,
	&\qquad
	\kappa &\longmapsto
	\EE U_m^\kappa(C_{d,t,\kappa}^{\rm typ}),
	\quad 0\leq m\leq d.
\end{aligned}
\]
is $C^\infty$ on $(-\varepsilon,\varepsilon)$.
\end{proposition}

\begin{proof}
Choose $\varepsilon>0$ with $\pi(d-1)\sqrt{\varepsilon}<t$.
Then \eqref{eq:AdmissibleRegime} holds for every
$\kappa\in(-\varepsilon,\varepsilon)$. The function
\[
\rho_t(\kappa)
:=
\begin{cases}
e^{-t/\sqrt{\kappa}},&\kappa>0,\\
0,&\kappa\leq0,
\end{cases}
\]
is $C^\infty$. Indeed, every derivative on the positive half-line is a
finite sum of powers of $\kappa$ multiplied by
$e^{-t/\sqrt{\kappa}}$, and each such derivative tends to zero
as $\kappa\downarrow0$. Consequently,
$p_\kappa(t)=1-\rho_t(\kappa)$ is $C^\infty$. The same observation
applied to the correction term in
\eqref{eq:CellIntensityAllCurvatures} shows that
$\gamma_{d,t,\kappa}$ is $C^\infty$ as well.
At $\kappa=0$,
\[
\gamma_{d,t,0}=\frac{2t^d}{\omega_{d+1}d!}>0.
\]
After decreasing $\varepsilon$ if necessary,
$1/\gamma_{d,t,\kappa}$ is therefore $C^\infty$ on the same
interval. Equation~\eqref{eq:VolumeViaIntensity} proves the
assertion for the mean volume, and
\eqref{eq:MeanExteriorAngleFunctionals} proves it for all mean
$\kappa$-intrinsic volumes.
For $2\leq j<d$, the same argument in dimension $j$ applies to
$\gamma_{j,t,\kappa}$, while
$\gamma_{1,t,\kappa}=t/\pi$ by
\eqref{eq:OneDimensionalCellIntensity}. Equations
\eqref{eq:MeanTypicalCellCrofton} and
\eqref{eq:MeanCroftonEndpoints} now prove the assertion for all
mean Crofton functionals.
\end{proof}

\paragraph{Poissonization of Schl\"afli's formula and the empty-process correction.}
We give an alternative derivation of the elliptic
cell-intensity formula by Poissonizing Schl\"afli's
formula. This approach was also used by Hug and
Reichenbacher \cite[Theorem~9.5 and its proof]{HR17},
who work in spherical rather than elliptic space.
We also show that the exponential correction in
\eqref{eq:CellIntensityAllCurvatures} disappears naturally
if the whole elliptic space is counted with its Euler
characteristic on the empty-process event.

Fix $\kappa>0$. For $n\geq1$ deterministic hyperplanes in general
position in $\MM_\kappa^d$, let $R_d(n)$ be the number of elliptic
chambers. Schl\"afli's formula gives
\begin{equation}\label{eq:SchlaefliEllipticCount}
R_d(n)
=
\sum_{i=0}^{d}\binom{n-1}{i}
=
\sum_{j=0}^{\lfloor d/2\rfloor}\binom{n}{d-2j},
\qquad n\geq1,
\end{equation}
where $\binom{n}{r}=0$ for $r>n$. The first sum is half the
spherical chamber count, the second equality follows from Pascal's
identity.
At $n=0$, the second sum in \eqref{eq:SchlaefliEllipticCount}
equals
\[
\sum_{j=0}^{\lfloor d/2\rfloor}\binom{0}{d-2j}
=
\chi(\MM_\kappa^d)
=
\chi(\mathbb{RP}^d)
=
\begin{cases}
1,&d\text{ even},\\
0,&d\text{ odd}.
\end{cases}
\]
This suggests an alternative count $R_d^\chi(n)$: for $n\geq1$
retain the ordinary chamber count, but at $n=0$ assign the whole
elliptic space its Euler characteristic. Then
\begin{equation}\label{eq:EulerCompletedChamberCount}
R_d^\chi(n)
=
\sum_{j=0}^{\lfloor d/2\rfloor}\binom{n}{d-2j},
\qquad n=0,1,2,\ldots.
\end{equation}
The alternative count $R_d^\chi(n)$ has a natural
topological interpretation: it is the sum of the
Euler characteristics of the connected components
of the complement of the hyperplanes in $\MM_\kappa^d$. Indeed,
for $n\geq1$, each such component is contractible
and hence has Euler characteristic $1$.
For $n=0$, the only component is the whole elliptic
space, whose Euler characteristic is
$\chi(\MM_\kappa^d)$.

Recall from \eqref{eq:EllipticHyperplaneCount} that the total number
of hyperplanes is
$N\sim\operatorname{Poisson}(t/\sqrt{\kappa})$. Its factorial
moments give
\begin{equation}\label{eq:PoissonBinomialMoments}
\EE\binom{N}{r}
=
\frac{(t/\sqrt{\kappa})^r}{r!},
\qquad r=0,1,2,\ldots,
\end{equation}
and combining \eqref{eq:EulerCompletedChamberCount},
\eqref{eq:PoissonBinomialMoments} and
\eqref{eq:EllipticSpaceVolume} yields
\begin{equation}\label{eq:EulerCompletedCellIntensity}
\frac{\EE R_d^\chi(N)}{\cH_\kappa^d(\MM_\kappa^d)}
=
\frac{2}{\omega_{d+1}}
\sum_{j=0}^{\lfloor d/2\rfloor}
\kappa^j\frac{t^{d-2j}}{(d-2j)!}.
\end{equation}
Thus counting the connected components of the
complement of the hyperplanes by their Euler characteristics gives precisely
the polynomial part of \eqref{eq:CellIntensityAllCurvatures}
in every dimension.

Under our convention, the whole elliptic space is not
counted as a cell when $N=0$, so $R_d(0)=0$.
Taking $B=\MM_\kappa^d$ in
\eqref{eq:CellIntensityDefinition}, we obtain
\begin{equation}\label{eq:SchlaefliEmptyCorrection}
\gamma_{d,t,\kappa}
=
\frac{\EE R_d^\chi(N)}{\cH_\kappa^d(\MM_\kappa^d)}
-
\frac{2\kappa^{d/2}}{\omega_{d+1}}
\chi(\MM_\kappa^d)e^{-t/\sqrt{\kappa}}.
\end{equation}
Equations~\eqref{eq:EulerCompletedCellIntensity} and
\eqref{eq:SchlaefliEmptyCorrection} recover
\eqref{eq:CellIntensityAllCurvatures} for $\kappa>0$.
Counting the whole elliptic space as one cell when $N=0$ would
instead add
$2\kappa^{d/2}e^{-t/\sqrt{\kappa}}/\omega_{d+1}$ to the
right-hand side of \eqref{eq:SchlaefliEmptyCorrection}. The
exponential correction would then vanish in even dimensions and
be positive in odd dimensions.

\section{Expected geometric characteristics of the zero cell}
\label{sec:ZeroCell}

\subsection{Mean volume of the zero cell}
\label{subsec:zero-cell-volume}

Recall from Section~\ref{subsec:zero-cell-definition} that
$C_{d,t,\kappa}^p$ is the zero cell of $\PHT_t$, selected as the almost surely unique cell containing the
fixed reference point $p\in\MM_\kappa^d$.
In positive curvature, when the hyperplane process is empty,
the zero cell is defined to be the whole elliptic space.
We first
compute the mean volume of the zero cell for every $t>0$, without assuming that
all cells are bounded. We then evaluate the resulting integral
and derive consequences for volume moments and their dependence
on curvature. Finally, we derive identities for the expected face numbers
and mean total face volumes of the zero cell.

\begin{theorem}[Mean volume of the zero cell]
\label{thm:ZeroCellMeanVolume}
For every $d\geq2$, $\kappa\in\RR$ and $t>0$,
\begin{equation}\label{eq:ZeroCellMeanVolumeIntegral}
\EE\cH_\kappa^d(C_{d,t,\kappa}^p)
=
\omega_d\int_0^{R_\kappa}
e^{-tr/\pi}\sn_\kappa(r)^{d-1}\,\dint r
-
\begin{cases}
\displaystyle
\frac{\omega_{d+1}}{2\kappa^{d/2}}e^{-t/\sqrt{\kappa}},
&\kappa>0,\\[0.6em]
0,&\kappa\leq0,
\end{cases}
\end{equation}
where
\[
R_\kappa
:=
\begin{cases}
\pi/\sqrt{\kappa},&\kappa>0,\\
+\infty,&\kappa\leq0,
\end{cases}
\]
and $\sn_\kappa$ is the generalised sine function from
\eqref{eq:GeneralizedSine}. The expectation in
\eqref{eq:ZeroCellMeanVolumeIntegral} is finite for every $t>0$
when $\kappa\geq0$. For $\kappa<0$, it is finite if and only if
\[
t>t_{\rm crit}(d,\kappa)=\pi(d-1)\sqrt{-\kappa},
\]
and equals $+\infty$ for $0<t\leq t_{\rm crit}(d,\kappa)$.
\end{theorem}

\begin{remark}[Related literature]\label{rem:ZeroCellVolumeLiterature}
The hyperbolic zero cell and its relation to beta-star polytopes
were studied in \cite[Section~3.2]{BetaStarPolytopes}.
For curvature $-1$, the mean volume of the hyperbolic zero cell
and its finiteness threshold were obtained in
\cite[Theorem~7.1]{BuehlerHugThaeleVisibility}. After converting
their normalization of the hyperplane intensity to ours, their
formula agrees with the case $\kappa=-1$ of
\eqref{eq:ZeroCellMeanVolumeIntegral}. The parameter
$\gamma_d^\ast$ in \cite{BuehlerHugThaeleVisibility} equals $t/\pi$ here.
Related results for weighted typical spherical cells and their
quermassintegrals appear in \cite{KT21}, for tessellations
generated by a fixed number of great hyperspheres.  For the second volume moment of the hyperbolic zero cell,
see \cite{BuehlerThaeleZeroCell}.
\end{remark}

\begin{remark}[Critical intensity]
For $\kappa<0$, Theorem~\ref{thm:ZeroCellMeanVolume} gives
\[
\EE\cH_\kappa^d(C_{d,t,\kappa}^p)<\infty
\quad\Longleftrightarrow\quad
t>t_{\rm crit}(d,\kappa).
\]
At $t=t_{\rm crit}(d,\kappa)$, all cells are almost surely
compact by the result recalled in
Section~\ref{subsec:cells-faces}, but the zero cell has infinite
mean volume. This differs from the typical cell, whose mean
volume is finite there by
Proposition~\ref{prop:TypicalCellVolumeIntensity}.
\end{remark}

\begin{proof}[Proof of Theorem~\ref{thm:ZeroCellMeanVolume}]
Throughout the proof, membership in an elliptic chamber refers
to its support. Tonelli's theorem gives
\begin{equation}\label{eq:ZeroCellVolumeByHitting}
\EE\cH_\kappa^d(C_{d,t,\kappa}^p)
=
\int_{\MM_\kappa^d}
\PP\bigl(x\in|C_{d,t,\kappa}^p|\bigr)
\,\cH_\kappa^d(\dint x),
\end{equation}
as an equality in $[0,\infty]$, without a preliminary
boundedness or integrability assumption.

\medskip
\emph{Nonpositive curvature.}
Suppose first that $\kappa\leq0$. For fixed
$x\in\MM_\kappa^d$, let $\cA(p\,|\,x)$ be the set of hyperplanes
separating $p$ and $x$. Hyperplanes containing either endpoint
have $\widetilde\mu_{d-1}^\kappa$-measure zero. Apart from these
hyperplanes, separation is equivalent to meeting the geodesic
segment $[p,x]$. Hence \eqref{eq:HyperplaneSegmentCrofton} gives
\[
\widetilde\mu_{d-1}^\kappa\bigl(\cA(p\,|\,x)\bigr)
=
\frac{d_\kappa(p,x)}{\pi}.
\]
By Definition~\ref{def:ZeroCell}, $x$ belongs to the zero cell
exactly when no hyperplane of the process separates $p$ and $x$.
The void probability formula for Poisson processes therefore yields
\begin{equation}\label{eq:ZeroCellHittingNonpositive}
\PP(x\in C_{d,t,\kappa}^p)
=
\PP\bigl(\eta_t^\kappa(\cA(p\,|\,x))=0\bigr)
=
e^{-t d_\kappa(p,x)/\pi}.
\end{equation}
Substituting \eqref{eq:ZeroCellHittingNonpositive} into
\eqref{eq:ZeroCellVolumeByHitting} and using geodesic polar
coordinates around $p$ gives
\begin{equation}\label{eq:ZeroCellPolarIntegral}
\EE\cH_\kappa^d(C_{d,t,\kappa}^p)
=
\omega_d\int_0^\infty
e^{-tr/\pi}\sn_\kappa(r)^{d-1}\,\dint r,
\qquad \kappa\leq0.
\end{equation}
Here $\sn_\kappa(r)^{d-1}$ is the radial volume density, and
$\omega_d$ is the area of the unit sphere in the tangent space $T_p\MM_\kappa^d$ at $p$.

Since $\sn_\kappa(r)\sim r(1+o(1))$ as $r\downarrow0$, the integral in
\eqref{eq:ZeroCellPolarIntegral} is finite near zero. For
$\kappa=0$, it is also finite at infinity for every $t>0$.
For $\kappa<0$, we have
\[
e^{-tr/\pi}\sn_\kappa(r)^{d-1}
=
(2\sqrt{-\kappa})^{-(d-1)}
\exp\!\left(\left((d-1)\sqrt{-\kappa}-\frac{t}{\pi}\right)r\right)(1+o(1)),\qquad r\to\infty.
\]
Thus \eqref{eq:ZeroCellPolarIntegral} is finite exactly when
$t>\pi(d-1)\sqrt{-\kappa}$. At equality the integrand tends to a strictly
positive constant, and below this value it grows exponentially.
This proves \eqref{eq:ZeroCellMeanVolumeIntegral} and the stated
finiteness criterion for $\kappa\leq0$.

\medskip
\emph{Positive curvature.}
Suppose now that $\kappa>0$. Recall that
$q_\kappa:\SS_\kappa^d\to\MM_\kappa^d$
identifies antipodal points. By isometry invariance, we may assume that
$p=q_\kappa(\widehat p)$, where $\widehat p=(1,0,\ldots,0)$ is the north pole.
Let $\widehat Z$ be the spherical chamber containing
$\widehat p$ when $\eta_t^\kappa\neq\emptyset$, and put
$\widehat Z:=\SS_\kappa^d$ otherwise.

On the nonempty-process event, the volume of the
elliptic zero cell equals that of $\widehat Z$.
On the empty-process event, the volume of $\widehat Z$
is twice the volume of the whole elliptic space.
Consequently,
\begin{equation}\label{eq:volume_identity_expected_zero_cell_spherical}
\cH_\kappa^d(C_{d,t,\kappa}^p)
=
\cH_\kappa^d(\widehat Z)
-
{\bf1}\{\eta_t^\kappa=\emptyset\}
\cH_\kappa^d(\MM_\kappa^d).
\end{equation}

We now compute the expected volume of $\widehat Z$.
Fix $\widehat x\in\SS_\kappa^d$ and let $s$ be its
spherical distance from $\widehat p$, where
$0<s<\pi/\sqrt{\kappa}$.
The invariant measure of elliptic hyperplanes whose
spherical lifts separate $\widehat p$ and $\widehat x$
equals $s/\pi$.
Indeed, partition the shortest spherical geodesic arc
from $\widehat p$ to $\widehat x$ into segments of length
less than $\pi/(2\sqrt{\kappa})$. Each segment projects
to a minimizing geodesic segment in elliptic space,
so \eqref{eq:HyperplaneSegmentCrofton} applies.
Since the original spherical arc has length less than
$\pi/\sqrt{\kappa}$, almost every great hypersphere
meets it at most once. Thus the measures of hyperplanes
meeting the individual segments add up to $s/\pi$.

For each fixed $\widehat x$, almost surely
$\widehat x\in\widehat Z$ if and only if no lifted
great hypersphere separates $\widehat x$ from
$\widehat p$. 
 By the void probability formula for Poisson processes, $\PP(\widehat x\in\widehat Z)=e^{-ts/\pi}$.
Tonelli's theorem and spherical polar coordinates
around $\widehat p$ therefore give
\[
\EE\cH_\kappa^d(\widehat Z)
=
\omega_d\int_0^{\pi/\sqrt{\kappa}}
e^{-ts/\pi}\sn_\kappa(s)^{d-1}\,\dint s.
\]
Finally, by \eqref{eq:EllipticHyperplaneCount}, $\PP(\eta_t^\kappa=\emptyset)=e^{-t/\sqrt{\kappa}}$.
Taking expectations in~\eqref{eq:volume_identity_expected_zero_cell_spherical}
and using \eqref{eq:EllipticSpaceVolume} yields
\[
\EE\cH_\kappa^d(C_{d,t,\kappa}^p)
=
\omega_d\int_0^{\pi/\sqrt{\kappa}}
e^{-ts/\pi}\sn_\kappa(s)^{d-1}\,\dint s
-
\frac{\omega_{d+1}}{2\kappa^{d/2}}
e^{-t/\sqrt{\kappa}}.
\]
This is \eqref{eq:ZeroCellMeanVolumeIntegral}
for $\kappa>0$.
\end{proof}

\subsection{Consequences of the mean-volume formula}
\label{subsec:zero-cell-volume-consequences}

In this subsection we first evaluate the integral in
\eqref{eq:ZeroCellMeanVolumeIntegral} and observe that the mean volume of the zero cell smoothly depends
on the curvature. We then use volume weighting to obtain the second
moment of the typical-cell volume, give low-dimensional examples,
and derive the high-intensity expansion.

\paragraph{Closed forms for the mean volume of the zero cell.}
For $d\geq2$ and $t>0$, define the rational function
\begin{equation}\label{eq:ZeroCellRationalPart}
\mathcal R_{d,t}(\kappa)
:=
\omega_d\Gamma(d)\left(\frac{\pi}{t}\right)^d
\prod_{\substack{1\leq j\leq d-1\\
                 j\equiv d-1\ ({\rm mod}\,2)}}
\left(1+\frac{\pi^2j^2}{t^2}\kappa\right)^{-1},
\end{equation}
wherever the denominators are nonzero. In particular, it is
well defined for
$\kappa>-t^2/(\pi^2(d-1)^2)$. The index set in
\eqref{eq:ZeroCellRationalPart} consists of the positive odd
integers up to $d-1$ when $d$ is even, and the positive even
integers up to $d-1$ when $d$ is odd.

\begin{theorem}[Mean volume of the zero cell: Closed form]
\label{thm:ZeroCellVolumeClosedForm}
Let $d\geq2$ and $t>0$. For $\kappa=0$, and for $\kappa<0$
with $t>t_{\rm crit}(d,\kappa)$,
\begin{equation}\label{eq:ZeroCellVolumeRationalNonpositive}
\EE\cH_\kappa^d(C_{d,t,\kappa}^p)
=
\mathcal R_{d,t}(\kappa).
\end{equation}
For $\kappa>0$,
\begin{equation}\label{eq:ZeroCellVolumeRationalPositive}
\EE\cH_\kappa^d(C_{d,t,\kappa}^p)
=
\mathcal R_{d,t}(\kappa)
+
e^{-t/\sqrt{\kappa}}
\left(
(-1)^d\mathcal R_{d,t}(\kappa)
-
\frac{\omega_{d+1}}{2\kappa^{d/2}}
\right).
\end{equation}
In particular, the Euclidean value is
\begin{equation}\label{eq:ZeroCellVolumeEuclidean}
\EE\cH_0^d(C_{d,t,0}^p)
=
\mathcal R_{d,t}(0)
=
\omega_d\Gamma(d)\left(\frac{\pi}{t}\right)^d.
\end{equation}
\end{theorem}

\begin{proof}
We work in the parameter range where the expectation in
\eqref{eq:ZeroCellMeanVolumeIntegral} is finite. Set
$\lambda:=t/\pi$ and, for $0\leq n\leq d-1$, write
\[
I_n:=\int_0^{R_\kappa}
e^{-\lambda r}\sn_\kappa(r)^n\,\dint r,
\]
where $R_\kappa$ is as in Theorem~\ref{thm:ZeroCellMeanVolume}.
For $n\geq2$, direct differentiation gives
\[
\frac{\dint^2}{\dint r^2}\sn_\kappa(r)^n
=
n(n-1)\sn_\kappa(r)^{n-2}
-n^2\kappa\sn_\kappa(r)^n.
\]
Twice integrating by parts therefore yields
\begin{equation}\label{eq:ZeroCellIntegralRecursion}
(\lambda^2+n^2\kappa)I_n=n(n-1)I_{n-2},
\qquad 2\leq n\leq d-1.
\end{equation}
The boundary terms vanish. At a finite endpoint this follows
from the zeros of $\sn_\kappa^n$ and its first derivative, and
at infinity it follows from $\lambda>(d-1)\sqrt{-\kappa}$
when $\kappa<0$, or from exponential decay when $\kappa=0$.
The initial values of the recursion are
\[
\begin{array}{lll}
\kappa\leq0:
& I_0=\lambda^{-1},
& I_1=(\lambda^2+\kappa)^{-1},\\[0.5em]
\kappa>0:
& I_0=\dfrac{1-e^{-\lambda R_\kappa}}{\lambda},
& I_1=\dfrac{1+e^{-\lambda R_\kappa}}{\lambda^2+\kappa}.
\end{array}
\]
Iteration of \eqref{eq:ZeroCellIntegralRecursion}, including
these initial values when $d=2$ or $d=3$, gives
\begin{equation}\label{eq:ZeroCellIntegralProduct}
\omega_d\int_0^{R_\kappa}
e^{-tr/\pi}\sn_\kappa(r)^{d-1}\,\dint r
=
\begin{cases}
\mathcal R_{d,t}(\kappa),&\kappa\leq0,\\[0.4em]
\mathcal R_{d,t}(\kappa)
\bigl(1+(-1)^d e^{-t/\sqrt{\kappa}}\bigr),&\kappa>0.
\end{cases}
\end{equation}
Substituting \eqref{eq:ZeroCellIntegralProduct} into
\eqref{eq:ZeroCellMeanVolumeIntegral} proves
\eqref{eq:ZeroCellVolumeRationalNonpositive} and
\eqref{eq:ZeroCellVolumeRationalPositive}. Setting $\kappa=0$
in \eqref{eq:ZeroCellRationalPart} gives
\eqref{eq:ZeroCellVolumeEuclidean}.
\end{proof}

Equation~\eqref{eq:ZeroCellVolumeEuclidean} is the classical
Euclidean mean-volume formula \cite[Theorem 10.4.9]{SW} with the
normalization conversion in \eqref{eq:FlatMeasureNormalization}.
For $\kappa<0$ and $t\leq t_{\rm crit}(d,\kappa)$, the
expectation is infinite by Theorem~\ref{thm:ZeroCellMeanVolume}.
Any finite value obtained by evaluating the rational function
\eqref{eq:ZeroCellRationalPart} outside the range specified in
\eqref{eq:ZeroCellVolumeRationalNonpositive} is not a mean volume.

\begin{remark}
\label{rem:ZeroCellGammaRepresentations}
The rational function $\mathcal R_{d,t}$ has a representation in terms of Gamma functions.
To present it, put $\nu:=(d+1)/2$. For $\kappa<0$ and
$t>t_{\rm crit}(d,\kappa)$, set
$z:=t/(2\pi\sqrt{-\kappa})$. Then
\begin{equation}\label{eq:ZeroCellGammaNegative}
\mathcal R_{d,t}(\kappa)
=
\pi^{\nu-1}\Gamma(\nu)(-\kappa)^{-d/2}
\frac{\Gamma(z+1-\nu)}{\Gamma(z+\nu)}.
\end{equation}
Indeed, after setting $s=\sqrt{-\kappa}\,r$, the substitution
$u=e^{-2s}$ gives
\[
\int_0^\infty e^{-2zs}\sinh^{d-1}s\,\dint s
=
2^{-d}\frac{\Gamma(d)\Gamma(z+1-\nu)}{\Gamma(z+\nu)}.
\]
The Gamma duplication formula now yields
\eqref{eq:ZeroCellGammaNegative}. For $\kappa=-1$, this is
the expression in \cite[Theorem~7.1]{BuehlerHugThaeleVisibility}
after the conversion described in
Remark~\ref{rem:ZeroCellVolumeLiterature}.

For $\kappa>0$, put $y:=t/(2\pi\sqrt{\kappa})$. The integral
in \eqref{eq:ZeroCellIntegralProduct} also has the representation
\begin{equation}\label{eq:ZeroCellGammaPositive}
\omega_d\int_0^{\pi/\sqrt{\kappa}}
e^{-tr/\pi}\sn_\kappa(r)^{d-1}\,\dint r
=
\frac{2\pi^\nu\Gamma(\nu)}{\kappa^{d/2}}
\frac{e^{-\pi y}}
     {\Gamma(\nu+\mathrm{i}y)\Gamma(\nu-\mathrm{i}y)}.
\end{equation}
Here the product in the denominator is
$|\Gamma(\nu+\mathrm{i}y)|^2$.
To relate \eqref{eq:ZeroCellGammaPositive} to
\eqref{eq:ZeroCellGammaNegative}, the Gamma recurrence formula,
with the factors paired about zero, gives
\begin{equation}\label{eq:ZeroCellGammaPositiveRatio}
\mathcal R_{d,t}(\kappa)
=
\pi^{\nu-1}\Gamma(\nu)\kappa^{-d/2}\,
\mathrm{i}^{\,d}
\frac{\Gamma(1-\nu+\mathrm{i}y)}
     {\Gamma(\nu+\mathrm{i}y)}.
\end{equation}
Thus the same Gamma ratio occurs at a positive real argument in
\eqref{eq:ZeroCellGammaNegative} and at an imaginary argument
in \eqref{eq:ZeroCellGammaPositiveRatio}. The reflection formula
gives
\[
\frac{1}{\Gamma(\nu+\mathrm{i}y)\Gamma(\nu-\mathrm{i}y)}
=
\frac{\sin\!\bigl(\pi(\nu-\mathrm{i}y)\bigr)}{\pi}
\frac{\Gamma(1-\nu+\mathrm{i}y)}{\Gamma(\nu+\mathrm{i}y)}.
\]
Since
\[
2e^{-\pi y}\sin\!\bigl(\pi(\nu-\mathrm{i}y)\bigr)
=
\mathrm{i}^{\,d}\bigl(1+(-1)^d e^{-2\pi y}\bigr),
\]
combining \eqref{eq:ZeroCellGammaPositiveRatio} with
\eqref{eq:ZeroCellIntegralProduct} proves
\eqref{eq:ZeroCellGammaPositive}. This also identifies explicitly
the parity factor in \eqref{eq:ZeroCellIntegralProduct}.
To obtain the elliptic mean volume from
\eqref{eq:ZeroCellGammaPositive}, one still subtracts
$\omega_{d+1}e^{-t/\sqrt{\kappa}}/(2\kappa^{d/2})$,
as prescribed by \eqref{eq:ZeroCellMeanVolumeIntegral}.
\end{remark}

Recall that Proposition~\ref{prop:CurvatureSmoothness} established smooth
dependence on curvature for the cell intensity and several
mean geometric characteristics of the typical cell.
We now show that the mean volume of the zero cell
has the same property.

\begin{proposition}
\label{prop:ZeroCellVolumeSmoothness}
For fixed $d\geq2$ and $t>0$, the map
\[
\kappa\longmapsto\EE\cH_\kappa^d(C_{d,t,\kappa}^p)
\]
is finite and $C^\infty$ in a neighbourhood of zero.
\end{proposition}

\begin{proof}
Choose $0<\delta<t^2/(\pi^2(d-1)^2)$.
The rational function in \eqref{eq:ZeroCellRationalPart} has no
poles on $(-\delta,\delta)$, and
\eqref{eq:ZeroCellVolumeRationalNonpositive} applies on
$(-\delta,0]$. Every derivative of the correction term in
\eqref{eq:ZeroCellVolumeRationalPositive} tends to zero as
$\kappa\downarrow0$, as the exponential $e^{-t/\sqrt{\kappa}}$
dominates every negative power of $\kappa$. Extending that
correction by zero for $\kappa\leq0$ therefore gives a
$C^\infty$ function, which proves the assertion.
\end{proof}

\paragraph{The second moment of the typical-cell volume.}
For the following statements we assume
\eqref{eq:AdmissibleRegime}. Write
\begin{equation*}
V^{\rm typ}:=\cH_\kappa^d(C_{d,t,\kappa}^{\rm typ}),
\qquad
V^0:=\cH_\kappa^d(C_{d,t,\kappa}^p),
\end{equation*}
and recall from \eqref{eq:VolumeViaIntensity} that
$\EE V^{\rm typ}=p_\kappa(t)/\gamma_{d,t,\kappa}$.
In the next result we use the volume-weighting identity in
Proposition~\ref{prop:ZeroCellVolumeBias} to relate 
the mean volume of the zero cell to the second volume
moment of the typical cell. For results on the second moment of the zero-cell
volume, see \cite{BuehlerThaeleZeroCell}.

\begin{proposition}[Second volume moment of the typical cell]
	\label{prop:TypicalCellVolumeSecondMoment}
	Under \eqref{eq:AdmissibleRegime},
	\begin{equation}\label{eq:TypicalCellVolumeSecondMoment}
		\gamma_{d,t,\kappa}\EE[(V^{\rm typ})^2]
		=
		\EE\!\left[{\bf1}\{\eta_t^\kappa\neq\emptyset\}V^0\right].
	\end{equation}
	In particular, for $\kappa\leq0$,
	\begin{equation}\label{eq:TypicalCellVolumeMomentProduct}
		\EE[(V^{\rm typ})^2]
		=
		\EE V^{\rm typ}\,\EE V^0
		=
		\frac{\EE V^0}{\gamma_{d,t,\kappa}}.
	\end{equation}
	For $\kappa<0$, the second moment is finite if and only if
	$t>t_{\rm crit}(d,\kappa)$. At the critical intensity,
	\begin{equation}\label{eq:CriticalTypicalCellVolumeMoments}
		0<\EE V^{\rm typ}<\infty,
		\qquad
		\EE V^0=\EE[(V^{\rm typ})^2]=+\infty.
	\end{equation}
	For $\kappa\geq0$, the second moment is finite for every $t>0$.
\end{proposition}
\begin{proof}
	Applying \eqref{eq:def_vol_typ} to the volume functional and using
	\eqref{eq:ZeroCellVolumeBias}, we obtain
	\[
	\frac{\EE[(V^{\rm typ})^2]}{\EE V^{\rm typ}}
	=
	\EE\!\left[V^0\mid\eta_t^\kappa\neq\emptyset\right].
	\]
	Since $\EE V^{\rm typ}=p_\kappa(t)/\gamma_{d,t,\kappa}$,
	multiplication by $p_\kappa(t)$ proves
	\eqref{eq:TypicalCellVolumeSecondMoment}.
	
	For $\kappa>0$, the zero cell equals the whole elliptic space
	on the empty-process event, which has probability
	$e^{-t/\sqrt{\kappa}}$. Hence
	\[
	\EE\!\left[{\bf1}\{\eta_t^\kappa\neq\emptyset\}V^0\right]
	=
	\EE V^0
	-
	\frac{\omega_{d+1}}{2\kappa^{d/2}}e^{-t/\sqrt{\kappa}}.
	\]
	For $\kappa\leq0$, the process is almost surely nonempty,
	so \eqref{eq:TypicalCellVolumeMomentProduct} follows from
	\eqref{eq:TypicalCellVolumeSecondMoment} and
	\eqref{eq:VolumeViaIntensity}.
	
	All identities remain valid when the second moment is infinite,
	since the first moment is finite. The remaining assertions
	follow from Theorem~\ref{thm:ZeroCellMeanVolume}.
\end{proof}

\paragraph{Low-dimensional hyperbolic examples.}
For $\kappa<0$, combining
\eqref{eq:CellIntensityAllCurvatures},
\eqref{eq:ZeroCellVolumeRationalNonpositive},
\eqref{eq:VolumeViaIntensity} and
\eqref{eq:TypicalCellVolumeMomentProduct} gives the low-dimensional values collected in
Table~\ref{tab:HyperbolicCellVolumeMoments}.
The denominator of each expression for $\EE V^0$ contains
$t^2-t_{\rm crit}(d,\kappa)^2$. The same factor appears in
$\EE[(V^{\rm typ})^2]$ by
\eqref{eq:TypicalCellVolumeMomentProduct}, whereas
$\EE V^{\rm typ}$ remains finite at the critical intensity.

\begin{table}[htbp]
	\centering
	\small
	\setlength{\tabcolsep}{4pt}
	\renewcommand{\arraystretch}{1.35}
	\begin{tabular}{@{}ccccc@{}}
		\hline
		$d$ & $t_{\rm crit}$ & $\EE V^{\rm typ}$ & $\EE V^0$
		& $\EE[(V^{\rm typ})^2]$\\
		\hline
		$2$ & $\pi\sqrt{-\kappa}$
		& $\displaystyle\frac{4\pi}{t^2+2\kappa}$
		& $\displaystyle\frac{2\pi^3}{t^2+\pi^2\kappa}$
		& $\displaystyle\frac{8\pi^4}
		{(t^2+2\kappa)(t^2+\pi^2\kappa)}$\\[1.4em]
		$3$ & $2\pi\sqrt{-\kappa}$
		& $\displaystyle\frac{6\pi^2}{t(t^2+6\kappa)}$
		& $\displaystyle\frac{8\pi^4}{t(t^2+4\pi^2\kappa)}$
		& $\displaystyle\frac{48\pi^6}
		{t^2(t^2+6\kappa)(t^2+4\pi^2\kappa)}$\\[1.4em]
		$4$ & $3\pi\sqrt{-\kappa}$
		& $\displaystyle\frac{32\pi^2}
		{t^4+12\kappa t^2+24\kappa^2}$
		& $\displaystyle\frac{12\pi^6}
		{(t^2+\pi^2\kappa)(t^2+9\pi^2\kappa)}$
		& $\displaystyle\frac{384\pi^8}
		{
			(t^4+12\kappa t^2+24\kappa^2)
			(t^2+\pi^2\kappa)(t^2+9\pi^2\kappa)
		}$\\[1.3em]
		\hline
	\end{tabular}
	\caption{Cell-volume moments in hyperbolic space of curvature
		$\kappa<0$. The formulas for $\EE V^{\rm typ}$ hold for
		$t\geq t_{\rm crit}$, while those for $\EE V^0$ and
		$\EE[(V^{\rm typ})^2]$ hold for $t>t_{\rm crit}$.
		At equality the latter two moments are $+\infty$, as in
		\eqref{eq:CriticalTypicalCellVolumeMoments}.}
	\label{tab:HyperbolicCellVolumeMoments}
\end{table}

\paragraph{High-intensity asymptotics.}
As $t\to\infty$, the zero cell becomes increasingly small,
and its geometry approaches that of the Euclidean zero cell.
The finite product in \eqref{eq:ZeroCellRationalPart}
allows us to identify the first curvature-dependent
correction to the Euclidean mean volume.

\begin{corollary}
\label{cor:ZeroCellHighIntensity}
For fixed $d\geq2$ and $\kappa\in\RR$, as $t\to\infty$,
\begin{equation}\label{eq:ZeroCellHighIntensity}
\EE\cH_\kappa^d(C_{d,t,\kappa}^p)
=
\omega_d\Gamma(d)\left(\frac{\pi}{t}\right)^d
\left(
1-\frac{d(d-1)(d+1)}{6}\frac{\pi^2\kappa}{t^2}
+O_{d,\kappa}(t^{-4})
\right).
\end{equation}
\end{corollary}

The leading term in \eqref{eq:ZeroCellHighIntensity}
is the Euclidean value in \eqref{eq:ZeroCellVolumeEuclidean},
and the first relative correction is linear in the curvature.

\begin{proof}[Proof of Corollary~\ref{cor:ZeroCellHighIntensity}]
For fixed $d$ and $\kappa$, expand the finite product in
\eqref{eq:ZeroCellRationalPart}:
\[
\prod_{\substack{1\leq j\leq d-1\\
                 j\equiv d-1\ ({\rm mod}\,2)}}
\left(1+\frac{\pi^2j^2\kappa}{t^2}\right)^{-1}
=
1-\frac{\pi^2\kappa}{t^2}
\sum_{\substack{1\leq j\leq d-1\\
                j\equiv d-1\ ({\rm mod}\,2)}}j^2
+O_{d,\kappa}(t^{-4}).
\]
The finite sum is $d(d-1)(d+1)/6$.
For $\kappa\leq0$, \eqref{eq:ZeroCellVolumeRationalNonpositive}
therefore proves \eqref{eq:ZeroCellHighIntensity}, since its
assumptions hold for all sufficiently large $t$.
For $\kappa>0$, the additional term in
\eqref{eq:ZeroCellVolumeRationalPositive} is exponentially
small as $t\to\infty$ and remains smaller than every power of
$t^{-1}$ after division by the leading factor $t^{-d}$.
It does not change the expansion. 
\end{proof}

\subsection{Mean face numbers, Crofton functionals and intrinsic volumes of the zero cell}
\label{subsec:zero-cell-face-numbers}
We now turn to the face structure of the zero cell $C_{d,t,\kappa}^p$.
By definition, $f_j(c)=\#\cF_j(c)$ counts the
$j$-faces of a cell $c$, each counted once even
if it has two local incidences with $c$.
In positive curvature, when $\eta_t^\kappa=\emptyset$,
the zero cell $C_{d,t,\kappa}^p$ is defined to be
$\MM_\kappa^d$, regarded as a boundaryless object with
no proper faces. Thus $f_j(C_{d,t,\kappa}^p)=0$
for $j<d$, while $f_d(C_{d,t,\kappa}^p)=1$. Recall also that the Crofton functionals $U_k^\kappa$
were defined in \eqref{eq:CroftonFunctionalDefinition}.

\paragraph{Crofton functionals from face numbers.}
The following theorem relates the expected Crofton
functionals of the zero cell to its expected face numbers.
The identity remains valid even when the zero cell
is unbounded.

\begin{theorem}[Efron identity]
\label{thm:ZeroCellEfron}
For every $t>0$ and $0\leq k\leq d$,
\begin{equation}\label{eq:EfronIdentity}
\EE U_k^\kappa(C_{d,t,\kappa}^p)
=
\frac{k!}{2t^k}\,
\EE f_{d-k}(C_{d,t,\kappa}^p).
\end{equation}
The identity is understood in $[0,\infty]$.
\end{theorem}

\begin{remark}[Related literature]
	\label{rem:ZeroCellFaceNumberLiterature}
	For the Euclidean counterpart of \eqref{eq:EfronIdentity},
	see \cite[Equation~(6.12)]{HugSchneiderBook}.
	A general approach to Efron-type identities based on stopping
	sets was recently developed in
	\cite{LastMolchanovEfron}. In particular, Example~10.2 there
	recovers the identity relating the expected number of facets
	of the Euclidean Poisson zero cell to its hyperplane hitting
	measure.
	The complete expected $f$-vector of the Euclidean Poisson zero
	cell was determined in \cite[Theorem~1.1]{KabluchkoZeroFVector}.
	Related results for weighted spherical cells appear in
	\cite[Theorem~5.1]{KT21}, and the expected face numbers of
	the hyperbolic Poisson zero cell were obtained in
	\cite[Theorem~3.5]{BetaStarPolytopes}.
	The spherical results concern a fixed number of great
	hyperspheres and spherical faces, rather than the 
	Poissonized model in the elliptic space considered here.
\end{remark}

\begin{proof}[Proof of Theorem~\ref{thm:ZeroCellEfron}]
For $k=0$, \eqref{eq:EfronIdentity} follows from
\eqref{eq:CroftonEndpoints} and $f_d=1$. Fix $1\leq k\leq d$.
By local finiteness and general position, each
$(d-k)$-face of the zero cell has exactly $k$ supporting
hyperplanes, even when the face is unbounded.
A given unordered
$k$-tuple supports at most one such  face.

To apply the multivariate Mecke formula, consider first a
realization of the original process and insert hyperplanes
$H_1,\ldots,H_k$ whose intersection
$E:=H_1\cap\cdots\cap H_k$ has dimension $d-k$. Let $Z$ and
$Z^+$ denote the zero cells before and after insertion,
respectively. Then
\begin{equation}\label{eq:ZeroCellInsertionIntersection}
|Z^+|\cap E=|Z|\cap E.
\end{equation}
For $\kappa\leq0$, this follows because the inserted closed
halfspaces containing $p$ all contain $E$ in their boundaries.
For $\kappa>0$, choose a spherical lift of $p$ and the closed
hemispheres containing that lift. Their inequalities become
equalities on the spherical lift of $E$, so projection gives
\eqref{eq:ZeroCellInsertionIntersection} as well. This also
covers an empty original process: after insertion the whole
elliptic flat $E$ is one  $(d-k)$-face of $Z^+$.

Outside a set of hyperplane tuples of product measure zero,
a nonempty intersection in
\eqref{eq:ZeroCellInsertionIntersection} has nonempty relative
interior in $E$ and hence contributes one $(d-k)$-face. Tuples
with an empty or nontransverse intersection contribute zero.
Consequently, the multivariate Mecke formula
\cite[Theorem~4.4]{LPbook}, followed by
\eqref{eq:Formula_1}, gives
\begin{equation}\label{eq:ZeroCellFaceCroftonIntegral}
\begin{split}
\EE f_{d-k}(C_{d,t,\kappa}^p)
&=
\frac{t^k}{k!}
\int_{A_\kappa(d,d-1)^k}
\PP\bigl(|C_{d,t,\kappa}^p|\cap H_1\cap\cdots\cap H_k
          \neq\emptyset\bigr)\,
(\widetilde\mu_{d-1}^\kappa)^{\otimes k}
\bigl(\dint(H_1,\ldots,H_k)\bigr)\\
&=
\frac{t^k}{k!}\,
\EE\int_{A_\kappa(d,d-k)}
{\bf1}\{|C_{d,t,\kappa}^p|\cap E\neq\emptyset\}
\,\widetilde\mu_{d-k}^\kappa(\dint E).
\end{split}
\end{equation}
The first integrand is understood to vanish unless the
intersection of the $k$ hyperplanes has dimension $d-k$.
The factor $k!$ accounts for the orderings of the supporting
hyperplanes. By \eqref{eq:CroftonFunctionalDefinition}, the last
line of \eqref{eq:ZeroCellFaceCroftonIntegral} equals
$2t^k\EE U_k^\kappa(C_{d,t,\kappa}^p)/k!$, proving
\eqref{eq:EfronIdentity}. All sums and integrals are nonnegative, so Tonelli's theorem
applies without any finiteness assumption. In particular,
the argument remains valid when the zero cell is unbounded.
\end{proof}

\paragraph{$\kappa$-intrinsic volumes from face numbers.}
Our aim is to express the mean proper $\kappa$-intrinsic
volumes of the zero cell in terms of its expected face
numbers. To this end, we first invert the relation
\eqref{eq:CroftonIntrinsicTransform}. Subtracting the
expressions at consecutive indices of the same parity
gives, for every $c\in\cC_\kappa^d$,
\begin{equation}\label{eq:IntrinsicVolumesFromCrofton}
\upsilon_m^\kappa(c)
=
\begin{cases}
U_m^\kappa(c)-\kappa U_{m+2}^\kappa(c),
&0\leq m\leq d-2,\\[0.3em]
U_{d-1}^\kappa(c),&m=d-1.
\end{cases}
\end{equation}

For $\kappa>0$, recall that the zero cell equals
$\MM_\kappa^d$ when the hyperplane process is empty.
Since this object has no proper faces, we define
\[
\upsilon_m^\kappa(\MM_\kappa^d):=0,
\qquad \kappa>0,\quad 0\leq m<d.
\]
Thus the empty-process event contributes zero to
$\EE\upsilon_m^\kappa(C_{d,t,\kappa}^p)$.

\begin{corollary}
\label{cor:ZeroCellIntrinsicVolumesViaFaces}
Under \eqref{eq:AdmissibleRegime}, for $0\leq m\leq d-2$,
\begin{equation}\label{eq:ZeroCellIntrinsicVolumesViaFaces}
\EE\upsilon_m^\kappa(C_{d,t,\kappa}^p)
=
\frac{m!}{2t^m}\,
\EE f_{d-m}(C_{d,t,\kappa}^p)
-\kappa\frac{(m+2)!}{2t^{m+2}}\,
\EE f_{d-m-2}(C_{d,t,\kappa}^p).
\end{equation}
For $m=d-1$,
\begin{equation}\label{eq:ZeroCellPenultimateIntrinsicVolume}
\EE\upsilon_{d-1}^\kappa(C_{d,t,\kappa}^p)
=
\frac{(d-1)!}{2t^{d-1}}\,
\EE f_1(C_{d,t,\kappa}^p)
-
\begin{cases}
\displaystyle
\frac{e^{-t/\sqrt{\kappa}}}{2\kappa^{(d-1)/2}},&\kappa>0,\\[0.6em]
0,&\kappa\leq0.
\end{cases}
\end{equation}
For $\kappa\leq0$ these identities are understood in
$[0,\infty]$, with the second term in
\eqref{eq:ZeroCellIntrinsicVolumesViaFaces} omitted when
$\kappa=0$. For $\kappa>0$ all expectations displayed above
are finite.
\end{corollary}

\begin{proof}
For $\kappa\leq0$, combine
\eqref{eq:IntrinsicVolumesFromCrofton} with
\eqref{eq:EfronIdentity}. Since $-\kappa\geq0$, this involves
no subtraction of infinite expectations.
For $\kappa>0$, \eqref{eq:CroftonWholeEllipticSpace} gives
\[
U_m^\kappa(\MM_\kappa^d)
-\kappa U_{m+2}^\kappa(\MM_\kappa^d)=0,
\qquad 0\leq m\leq d-2.
\]
Thus the first line of \eqref{eq:IntrinsicVolumesFromCrofton}
remains valid on the empty-process event under our extension.
At $m=d-1$, however, the empty-event contribution to
$\EE U_{d-1}^\kappa(C_{d,t,\kappa}^p)$ is
$e^{-t/\sqrt{\kappa}}/(2\kappa^{(d-1)/2})$ and must be
subtracted. Applying \eqref{eq:EfronIdentity} proves
\eqref{eq:ZeroCellIntrinsicVolumesViaFaces} and
\eqref{eq:ZeroCellPenultimateIntrinsicVolume}. All Crofton
functionals are bounded by the corresponding values in
\eqref{eq:CroftonWholeEllipticSpace}, which also proves the
claimed finiteness.
\end{proof}

\paragraph{Explicit hyperbolic mean intrinsic volumes.}
We shall express the mean intrinsic volumes of the hyperbolic zero
cell in terms of the coefficients $A[n,\ell]$ introduced in
\cite{KabluchkoZeroFVector}. To define them, let
$Q_0(x)=Q_1(x)=1$ and, for $n\geq2$, put
\[
Q_n(x)
:=
\prod_{\substack{1\leq j\leq n-1\\
		j\equiv n-1\ ({\rm mod}\,2)}}
(1+j^2x^2).
\]
For $n\geq0$ and $\ell\in\mathbb Z$, define
\[
A[n,\ell]
:=
\begin{cases}
	[x^\ell]Q_n(x),&\ell\text{ even},\\
	[x^\ell]\bigl(\tanh(\pi/(2x))Q_n(x)\bigr),
	&\ell\text{ odd},\ n\text{ even},\\
	[x^\ell]\bigl(\coth(\pi/(2x))Q_n(x)\bigr),
	&\ell\text{ odd},\ n\text{ odd},
\end{cases}
\]
where $[x^\ell]$ denotes coefficient extraction
from the corresponding formal Laurent series in $x^{-1}$.
In particular, $A[n,0]=1$ and $A[n,\ell]=0$
for $\ell>n$.
These coefficients also appear in the expected face numbers
of hyperbolic Poisson zero cells, see
\cite[Equation~(3.15)]{BetaStarPolytopes}. For $\alpha>(d-1)\pi$, we also put
\begin{equation}\label{eq:HyperbolicIntrinsicVolumeFactor}
T_d(\alpha)
:=
\left(\frac{\alpha}{2\pi}\right)^d
\frac{
\Gamma\!\left(\frac{\alpha}{2\pi}-\frac{d-1}{2}\right)
}{
\Gamma\!\left(\frac{\alpha}{2\pi}+\frac{d+1}{2}\right)
}.
\end{equation}

\begin{theorem}[Mean intrinsic volumes of the hyperbolic zero cell]
\label{prop:HyperbolicZeroCellIntrinsicVolumes}
Let $\kappa<0$ and $t>t_{\rm crit}(d,\kappa)$. Then, for
$0\leq m\leq d-1$,
\begin{equation}\label{eq:HyperbolicZeroCellIntrinsicVolumes}
\EE\upsilon_m^\kappa(C_{d,t,\kappa}^p)
=
\frac{A[d,m]}{2}
\left(\frac{\pi}{t}\right)^m
T_d\!\left(\frac{t}{\sqrt{-\kappa}}\right).
\end{equation}
All these expectations are finite.
\end{theorem}

\begin{proof}
Put $a:=\sqrt{-\kappa}$ and $\alpha:=t/a$. By
\eqref{eq:HyperplaneMeasureScale}, rescaling to curvature $-1$
replaces the parameter $t$ by $\alpha$. The conversion in
\eqref{eq:HyperbolicParameterConversion} and
\cite[Equation~(3.14)]{BetaStarPolytopes} identifies this
$\alpha$ with the beta-star parameter in that reference.
For $0\leq q\leq d$, define $T_q(\alpha)$ by
\eqref{eq:HyperbolicIntrinsicVolumeFactor} with $d$ replaced
by $q$. All these quantities are finite and positive, and
$T_0(\alpha)=T_1(\alpha)=1$.

For $1\leq\ell\leq d$, write
$F_\ell:=\EE f_{d-\ell}(C_{d,t,\kappa}^p)$.
The duality in \cite[Equation~(3.14)]{BetaStarPolytopes} and
\cite[Theorem~3.5]{BetaStarPolytopes} give
\begin{equation}\label{eq:HyperbolicZeroCellFaceNumberSum}
F_\ell
=
\frac{\pi^\ell}{\ell!}
\sum_{\substack{\ell\leq q\leq d\\q\equiv d\ ({\rm mod}\,2)}}
T_q(\alpha)\bigl(A[q,\ell]-A[q-2,\ell]\bigr).
\end{equation}
Following that reference, entries $A[q,\ell]$ with $q<\ell$
are interpreted as zero. In particular, the possible term
$A[-1,1]$ is zero. The formula in
\eqref{eq:HyperbolicZeroCellFaceNumberSum} is used only for
$\ell\geq1$, the case $m=0$ will be treated separately below.

We need two recurrences. The definition of the coefficients
$A[q,\ell]$ implies
\begin{equation}\label{eq:ZeroCellCoefficientRecurrence}
	A[q,\ell]-A[q-2,\ell]
	=(q-1)^2A[q-2,\ell-2],
	\qquad q\geq2,\quad \ell\in\mathbb Z.
\end{equation}
Indeed, $Q_q(x)=(1+(q-1)^2x^2)Q_{q-2}(x)$, and the
additional hyperbolic function in the definition of
$A[q,\ell]$ is the same for $q$ and $q-2$.
Similarly, the functional equation $\Gamma(z+1)=z\Gamma(z)$
applied to \eqref{eq:HyperbolicIntrinsicVolumeFactor} gives
\begin{equation}\label{eq:ZeroCellFactorRecurrence}
	T_q(\alpha)-T_{q-2}(\alpha)
	=
	\frac{\pi^2(q-1)^2}{\alpha^2}T_q(\alpha),
	\qquad 2\leq q\leq d.
\end{equation}

Fix $1\leq m\leq d-1$. We combine the expected
face-number formula \eqref{eq:HyperbolicZeroCellFaceNumberSum}
with \eqref{eq:ZeroCellIntrinsicVolumesViaFaces}, or
\eqref{eq:ZeroCellPenultimateIntrinsicVolume} when $m=d-1$.
This gives
\begin{equation}\label{eq:HyperbolicIntrinsicVolumeTwoSums}
\begin{split}
2\left(\frac{t}{\pi}\right)^m
\EE\upsilon_m^\kappa(C_{d,t,\kappa}^p)
={}&
\sum_{\substack{m\leq q\leq d\\q\equiv d\ ({\rm mod}\,2)}}
T_q(\alpha)\bigl(A[q,m]-A[q-2,m]\bigr)\\
&+\frac{\pi^2}{\alpha^2}
\sum_{\substack{m+2\leq q\leq d\\q\equiv d\ ({\rm mod}\,2)}}
T_q(\alpha)\bigl(A[q,m+2]-A[q-2,m+2]\bigr).
\end{split}
\end{equation}
The second sum is empty when $m=d-1$. By \eqref{eq:ZeroCellCoefficientRecurrence} and
\eqref{eq:ZeroCellFactorRecurrence}, this sum equals
\[
\sum_{\substack{m+2\leq q\leq d\\q\equiv d\ ({\rm mod}\,2)}}
\bigl(T_q(\alpha)-T_{q-2}(\alpha)\bigr)A[q-2,m].
\]
The missing lower-index terms vanish because $A[q-2,m]=0$
when $q<m+2$. Hence the right-hand side of
\eqref{eq:HyperbolicIntrinsicVolumeTwoSums} reduces to
\[
\sum_{\substack{m\leq q\leq d\\q\equiv d\ ({\rm mod}\,2)}}
T_q(\alpha)A[q,m]
-
\sum_{\substack{m+2\leq q\leq d\\q\equiv d\ ({\rm mod}\,2)}}
T_{q-2}(\alpha)A[q-2,m]
=
T_d(\alpha)A[d,m].
\]
This proves \eqref{eq:HyperbolicZeroCellIntrinsicVolumes}
for $1\leq m<d$.
For $m=0$, we use $f_d(C_{d,t,\kappa}^p)=1$ and
\eqref{eq:ZeroCellIntrinsicVolumesViaFaces} to obtain
\[
2\EE\upsilon_0^\kappa(C_{d,t,\kappa}^p)
=1+\frac{2a^2}{t^2}F_2.
\]
Since $A[q,0]=1$ for $q\geq0$, equations
\eqref{eq:HyperbolicZeroCellFaceNumberSum},
\eqref{eq:ZeroCellCoefficientRecurrence} and
\eqref{eq:ZeroCellFactorRecurrence} yield
\[
\begin{split}
2\EE\upsilon_0^\kappa(C_{d,t,\kappa}^p)
&=
1+\frac{\pi^2}{\alpha^2}
\sum_{\substack{2\leq q\leq d\\q\equiv d\ ({\rm mod}\,2)}}
(q-1)^2T_q(\alpha)\\
&=
1+
\sum_{\substack{2\leq q\leq d\\q\equiv d\ ({\rm mod}\,2)}}
\bigl(T_q(\alpha)-T_{q-2}(\alpha)\bigr)
=T_d(\alpha).
\end{split}
\]
The last equality uses $T_0(\alpha)=T_1(\alpha)=1$.
Together with $A[d,0]=1$, this proves the remaining case.
Finiteness follows either from these finite sums or directly
from \eqref{eq:HyperbolicZeroCellIntrinsicVolumes}.
\end{proof}

\begin{remark}[Euclidean counterpart]
\label{rem:ZeroCellIntrinsicVolumeCounterparts}
For $\kappa=0$, the formula in
\cite[Section~1.3]{KabluchkoZeroFVector}, expressed in the
normalization \eqref{eq:FlatMeasureNormalization}, is
\begin{equation}\label{eq:EuclideanZeroCellIntrinsicVolumes}
\EE\upsilon_m^0(C_{d,t,0}^p)
=
\frac{A[d,m]}{2}\left(\frac{\pi}{t}\right)^m,
\qquad 0\leq m<d.
\end{equation}
Thus \eqref{eq:HyperbolicZeroCellIntrinsicVolumes} multiplies
each Euclidean mean in \eqref{eq:EuclideanZeroCellIntrinsicVolumes}
by the same factor $T_d(t/\sqrt{-\kappa})$, independently of
$m$. The Gamma recurrence also gives
\[
T_d(\alpha)
=
\prod_{\substack{1\leq j\leq d-1\\
                 j\equiv d-1\ ({\rm mod}\,2)}}
\left(1-\frac{\pi^2j^2}{\alpha^2}\right)^{-1}.
\]
In particular, $T_d(\alpha)\to1$ as $\alpha\to\infty$,
so \eqref{eq:HyperbolicZeroCellIntrinsicVolumes} tends to
\eqref{eq:EuclideanZeroCellIntrinsicVolumes} as $\kappa\uparrow0$
with $t$ fixed. This is also the factor relating
$\mathcal R_{d,t}(\kappa)$ in \eqref{eq:ZeroCellRationalPart}
to $\mathcal R_{d,t}(0)$ for $\kappa<0$ in the finite-mean regime.
\end{remark}

\begin{remark}[Spherical counterpart]
For positive curvature, a related formula is available
for tessellations of the unit sphere $\SS^d$ generated
by a fixed number $n$ of independent, uniformly
distributed great hyperspheres. Let $W_{n,d}$ denote
the spherical cell containing the north pole.
By \cite[Corollary~5.9]{KT21}, for $n\geq d+1$,
\[
\EE\upsilon_m^{\mathrm{sph}}(W_{n,d})
=
\frac{n!}{2\pi^n}B\{n+m,d\}A[d,m],
\qquad 0\leq m<d,
\]
where $B\{\cdot,\cdot\}$ is defined in that reference
and $\upsilon_m^{\mathrm{sph}}$ denotes the spherical
intrinsic volume given by the face-angle representation
\eqref{eq:ExteriorAngleFunctionals} with $\kappa=1$,
applied directly to spherical chambers.
Under rescaling to
curvature $\kappa>0$, the right-hand side is
multiplied by $\kappa^{-m/2}$.
The proof in \cite{KT21} likewise uses a difference
of two Crofton means and a telescoping sum.
Poissonizing the number of great hyperspheres
yields the corresponding mean intrinsic volumes
of the spherical cell containing the north pole.
The cases $n\leq d$ must be treated separately,
since they are not covered by the displayed formula. We omit the details.
\end{remark}

\subsection{Mean total face volumes of the zero cell}
\label{subsec:zero-cell-face-volumes}

For a compact chamber $c\in\cC_\kappa^d$, or a possibly
unbounded cell $c$ of $\PHT_t$, define
\begin{equation}\label{eq:TotalFaceVolumes}
L_k(c)
:=
\sum_{F\in\cF_k(c)}\cH_\kappa^k(F),
\qquad 0\leq k\leq d.
\end{equation}
In particular,
\[
L_0(c)=f_0(c),
\qquad
L_d(c)=\cH_\kappa^d(c).
\]
Our aim is to determine $\EE L_k(C_{d,t,\kappa}^p)$.
Unlike the preceding intrinsic-volume results, the identities
in this subsection hold for every $t>0$, without requiring
the zero cell to be bounded.
For unbounded cells, the sum in
\eqref{eq:TotalFaceVolumes} is understood in
$[0,\infty]$ and may contain infinitely many terms.

For $\kappa>0$, recall that the zero cell equals
$\MM_\kappa^d$ when $\eta_t^\kappa=\emptyset$.
Since this object has no proper faces, we set
\[
L_k(\MM_\kappa^d)=0\quad(0\leq k<d),
\qquad
L_d(\MM_\kappa^d)=\cH_\kappa^d(\MM_\kappa^d).
\]

\begin{corollary}
\label{cor:ZeroCellSkeletonVolume}
For every $t>0$ and every $0\leq k\leq d$,
\begin{equation}\label{eq:ZeroCellSkeletonVolume}
\EE L_k(C_{d,t,\kappa}^p)
=
\frac{t^{d-k}}{(d-k)!}
\frac{\omega_{k+1}}{\omega_{d+1}}
\EE\cH_\kappa^d(C_{d,t,\kappa}^p),
\end{equation}
as an identity in $[0,\infty]$.
\end{corollary}

Together with Theorems~\ref{thm:ZeroCellMeanVolume} and~\ref{thm:ZeroCellVolumeClosedForm},
Equation~\eqref{eq:ZeroCellSkeletonVolume} determines all the mean total
face volumes explicitly. For $k<d$, its coefficient is the
$k$-skeleton density from \eqref{eq:SkeletonIntensities}, so
\begin{equation*}
\EE L_k(C_{d,t,\kappa}^p)
=
\sigma_{d,k,t,\kappa}\,
\EE\cH_\kappa^d(C_{d,t,\kappa}^p),
\qquad 0\leq k<d.
\end{equation*}
This reformulation is restricted to $k<d$: in positive curvature
$\sigma_{d,d,t,\kappa}=p_\kappa(t)$, whereas the prefactor in
\eqref{eq:ZeroCellSkeletonVolume} is one for $k=d$.

\begin{proof}
For $k=d$, the assertion follows from \eqref{eq:TotalFaceVolumes}.
Fix $0\leq k<d$ and write $Z:=C_{d,t,\kappa}^p$. By local finiteness and general position, every
$k$-face has exactly $d-k$ supporting hyperplanes,
even when the face is unbounded. A given unordered tuple
supports at most one such face of the zero cell.

We use the insertion argument from the proof of
Theorem~\ref{thm:ZeroCellEfron}, with section volume in place of
the nonemptiness indicator. Insert $d-k$ hyperplanes with a
$k$-dimensional intersection $E$, and let $Z^+$ denote the new
zero cell. Equation~\eqref{eq:ZeroCellInsertionIntersection}
gives $|Z^+|\cap E=|Z|\cap E$.
Outside a set of hyperplane tuples of product measure zero, a
nonempty section is the support of one $k$-face of
$Z^+$, whose volume is therefore $\cH_\kappa^k(|Z|\cap E)$.
Tuples with an empty or nontransverse hyperplane intersection
contribute zero.

The multivariate Mecke formula \cite[Theorem~4.4]{LPbook},
followed by \eqref{eq:Formula_1}, consequently yields
\begin{equation}\label{eq:ZeroCellFaceVolumeIntegral}
\EE L_k(Z)
=
\frac{t^{d-k}}{(d-k)!}
\EE\int_{A_\kappa(d,k)}
\cH_\kappa^k(|Z|\cap E)\,
\widetilde\mu_k^\kappa(\dint E).
\end{equation}
The factorial removes the orderings of the supporting hyperplanes.
For $\kappa>0$, the same calculation includes an empty original
process: then $|Z|=\MM_\kappa^d$, and after insertion the uncut
flat $E$ contributes one  $k$-face, with its projective
$k$-volume counted once. Thus the expectation on the right-hand
side of \eqref{eq:ZeroCellFaceVolumeIntegral} is unconditioned.

For every realization, \eqref{eq:CroftonSpecialCase} with
$B=|Z|$ gives
\[
\int_{A_\kappa(d,k)}
\cH_\kappa^k(|Z|\cap E)\,
\widetilde\mu_k^\kappa(\dint E)
=
\frac{\omega_{k+1}}{\omega_{d+1}}\cH_\kappa^d(Z).
\]
This applies also when $Z$ is unbounded or $|Z|$ is
the whole elliptic space.
Substitution into \eqref{eq:ZeroCellFaceVolumeIntegral} proves
\eqref{eq:ZeroCellSkeletonVolume}. All sums and integrals are
nonnegative, so Tonelli's theorem justifies the interchanges
without any preliminary moment assumption.
\end{proof}

By Theorem~\ref{thm:ZeroCellMeanVolume} and
\eqref{eq:ZeroCellSkeletonVolume}, every
$\EE L_k(C_{d,t,\kappa}^p)$ is finite when
$\kappa\geq0$ and $t>0$. For $\kappa<0$ and
every $0\leq k\leq d$,
\[
\EE L_k(C_{d,t,\kappa}^p)<\infty
\quad\Longleftrightarrow\quad
t>t_{\rm crit}(d,\kappa).
\]
At $t=t_{\rm crit}(d,\kappa)$,
\begin{equation}\label{eq:ZeroCellFaceVolumesAtCriticality}
L_k(C_{d,t,\kappa}^p)<\infty\quad\text{almost surely},
\qquad
\EE L_k(C_{d,t,\kappa}^p)=+\infty,
\qquad 0\leq k\leq d.
\end{equation}
Indeed, the zero cell is a compact polytope with finitely
many faces by Section~\ref{subsec:cells-faces}, but its
mean volume is infinite. Since $L_0=f_0$, the divergence
of the expected number of vertices holds throughout
$0<t\leq t_{\rm crit}(d,\kappa)$.

\begin{remark}
\label{rem:ZeroCellFaceVolumeCurvatureFactor}
Recall the factor $T_d$ defined in
\eqref{eq:HyperbolicIntrinsicVolumeFactor}. For $\kappa<0$ and
$t>t_{\rm crit}(d,\kappa)$,
\eqref{eq:ZeroCellVolumeRationalNonpositive},
\eqref{eq:ZeroCellVolumeEuclidean} and
\eqref{eq:ZeroCellRationalPart} show that the hyperbolic mean
zero-cell volume is its Euclidean counterpart multiplied by
$T_d(t/\sqrt{-\kappa})$. Since the coefficient in
\eqref{eq:ZeroCellSkeletonVolume} is independent of $\kappa$,
\begin{equation*}
\EE L_k(C_{d,t,\kappa}^p)
=
T_d\!\left(\frac{t}{\sqrt{-\kappa}}\right)
\EE L_k(C_{d,t,0}^p),
\qquad 0\leq k\leq d.
\end{equation*}
Both tessellations use the same parameter $t$ and the normalization
\eqref{eq:FlatMeasureNormalization}. Thus the same factor that
relates the mean intrinsic volumes in
\eqref{eq:HyperbolicZeroCellIntrinsicVolumes} and
\eqref{eq:EuclideanZeroCellIntrinsicVolumes} also relates all
mean total face volumes, including the volume itself.
\end{remark}

\section{Typical faces and face intensities}
\label{sec:Kfaces}

\subsection{Typical \texorpdfstring{$k$}{k}-faces}
\label{subsec:typical-k-faces}

We now turn from full-dimensional cells to the lower-dimensional
faces of the tessellation. Fix $1\leq k\leq d-1$ and recall from
Section~\ref{subsec:cells-faces} that $\cF_k(\PHT_t)$ consists
of all distinct  $k$-faces, each counted once, regardless
of how many cells share it. We define the typical $k$-face by
sampling from this collection, rather than by selecting a face
of the typical full-dimensional cell, and introduce its
volume-weighted counterpart.
We then identify both distributions, up to isometry, with
suitable cell distributions of an intrinsic $k$-dimensional
Poisson hyperplane tessellation. This yields the face intensities,
from which we determine the mean face numbers and mean total
face volumes of the typical full-dimensional cell. Throughout
this section we assume \eqref{eq:AdmissibleRegime}.

\paragraph{The face space and centre kernels.}
Let $\cC_\kappa^{d,k}$ be the space of compact $k$-dimensional
polyhedral chambers embedded in $\MM_\kappa^d$, with their
supporting flats retained. For $\kappa\leq0$, these are compact
geodesically convex polytopes with nonempty relative interior in
some $E\in A_\kappa(d,k)$. For $\kappa>0$, they are the elliptic
chambers defined in Section~\ref{subsec:cells-faces}, now inside
such an $E$, together with a separate boundaryless whole-flat
object for each $E\in A_\kappa(d,k)$. We use the Borel structure
induced by the Hausdorff metric on embedded polytopes when
$\kappa\leq0$. For $\kappa>0$, we use the quotient Borel structure
on embedded spherical representatives under the antipodal map,
and adjoin a disjoint copy of $A_\kappa(d,k)$ for the boundaryless
whole-flat objects. As before, $|F|$ denotes the support of a
chamber $F$.

A \emph{covariant centre kernel} is a measurable probability
kernel $F\mapsto\zeta_F$ from $\cC_\kappa^{d,k}$ to
$\MM_\kappa^d$ such that
\begin{equation}\label{eq:FaceCentreKernelProperties}
\zeta_F(|F|)=1,
\qquad
\zeta_{\phi F}(\phi B)=\zeta_F(B),
\end{equation}
for every $F\in\cC_\kappa^{d,k}$, every Borel set
$B\subset\MM_\kappa^d$ and every $\phi\in I(\MM_\kappa^d)$.
Thus each face carries total mass one, located in its support.
This extends the centre maps of Definition~\ref{def:CentreMap},
which correspond to point-mass kernels.

Such a kernel can be chosen as follows. For $\kappa\leq0$, let
$z_k(F)$ be the circumcentre of $F$ in its supporting flat. For
$\kappa>0$, if $F$ is an ordinary elliptic chamber, let $z_k(F)$
be the projected normalized spherical barycentre, computed in its
supporting space as in \eqref{eq:SphericalCentre}. Put
\begin{equation}\label{eq:FaceCentreKernelConstruction}
\zeta_F(B)
:=
\begin{cases}
{\bf1}\{z_k(F)\in B\},
&F\text{ is not a boundaryless whole elliptic flat},\\[0.8em]
\displaystyle
\frac{\cH_\kappa^k(|F|\cap B)}{\cH_\kappa^k(F)},
&F\text{ is a boundaryless whole elliptic flat}.
\end{cases}
\end{equation}
These prescriptions are measurable and covariant by the argument
in Lemma~\ref{lem:CentreMapExists}, applied intrinsically in the
supporting flat, and by invariance of Hausdorff measure.
A whole elliptic flat has no distinguished covariant point,
because its isometry group acts transitively. The second line
of \eqref{eq:FaceCentreKernelConstruction} avoids making such a
choice. It applies only to the boundaryless object. A chamber
represented by a hemisphere retains its boundary structure even
when its support is the whole flat, and uses the first line.

\paragraph{Face intensity.}
Fix a covariant centre kernel $\zeta$. Its associated random
measure is
\begin{equation}\label{eq:KFaceCentreMeasure}
\xi_{k,t}^{\kappa,\zeta}(B)
:=
\sum_{F\in\cF_k(\PHT_t)}\zeta_F(B),
\qquad B\subset\MM_\kappa^d\ \text{Borel}.
\end{equation}
For point-mass kernels this is the counting measure of face
centres. In general it need not be a counting measure, but each
 face still contributes total mass one.

\begin{lemma}
\label{lem:KFaceIntensityFinite}
For every covariant centre kernel $\zeta$, there is a constant
$\gamma_{d,k,t,\kappa}^{\,\zeta}\in(0,\infty)$ such that
\[
\EE\xi_{k,t}^{\kappa,\zeta}(B)
=
\gamma_{d,k,t,\kappa}^{\,\zeta}\cH_\kappa^d(B)
\]
for every Borel set $B\subset\MM_\kappa^d$.
\end{lemma}

\begin{proof}
Suppose first that $\kappa\leq0$. Given a compact set $B$, choose
a compact geodesically convex window $W$ whose interior contains
$B$, and let $N_W$ be the number of hyperplanes meeting $W$.
A $k$-face meeting $B$ has exactly $d-k$ supporting hyperplanes,
all of which meet $W$. Once these hyperplanes have been chosen,
there are at most $2^{N_W}$ choices of sides of the remaining
hyperplanes in $W$. Since each face meeting $B$ has a relative
interior point in $W$, the number of such faces is at most
$\binom{N_W}{d-k}2^{N_W}$. By
\eqref{eq:FaceCentreKernelProperties}, this also bounds
$\xi_{k,t}^{\kappa,\zeta}(B)$. The bound has finite expectation
because $N_W$ is Poisson distributed with finite mean.
For $\kappa>0$, use the same bound with the total number $N$
of hyperplanes in \eqref{eq:EllipticHyperplaneCount}. It applies
also when a supporting flat is uncut. Thus the mean measure in
\eqref{eq:KFaceCentreMeasure} is locally finite in every case.
It is isometry invariant by \eqref{eq:FaceCentreKernelProperties},
so it is a finite nonnegative constant multiple of $\cH_\kappa^d$.
If that constant were zero, the measure
$\xi_{k,t}^{\kappa,\zeta}$ would almost surely vanish on every
member of a countable cover by bounded balls, and hence everywhere.
Since each face carries total mass one, there would be no
$k$-faces almost surely. This contradicts the strictly positive
$k$-skeleton density in \eqref{eq:SkeletonIntensities}.
\end{proof}

We call $\gamma_{d,k,t,\kappa}^{\,\zeta}$ the \emph{$k$-face
intensity}. Equivalently,
\begin{equation}\label{eq:KFaceIntensityDefinition}
\gamma_{d,k,t,\kappa}^{\,\zeta}
=
\frac{1}{\cH_\kappa^d(W)}
\EE\sum_{F\in\cF_k(\PHT_t)}\zeta_F(W),
\qquad 0<\cH_\kappa^d(W)<\infty.
\end{equation}
The quotient is independent of $W$. We retain the superscript
$\zeta$ until independence of the kernel has been established.

\paragraph{The typical face.}
The typical $k$-face samples faces with equal weight,
without weighting by their $k$-dimensional volumes. Its invariant
expectations are characterized by
\begin{equation}\label{eq:TypicalKFaceInvariantExpectations}
\gamma_{d,k,t,\kappa}^{\,\zeta}
\EE h(C_{d,t,\kappa}^{k,{\rm typ},\zeta})
=
\frac{1}{\cH_\kappa^d(W)}
\EE\sum_{F\in\cF_k(\PHT_t)}\zeta_F(W)h(F),
\end{equation}
for every nonnegative measurable isometry-invariant functional
$h:\cC_\kappa^{d,k}\to[0,\infty]$. Here isometry invariance means
$h(\phi F)=h(F)$ for all $\phi\in I(\MM_\kappa^d)$.

To specify a representative in the ambient space, recall the
probability kernel $\widetilde\lambda_\kappa(x,\cdot)$ from
\eqref{eq:RecenteringKernel}. It is supported on the isometries
sending $p$ to $x$, their inverses send $x$ to $p$.

\begin{definition}\label{def:TypicalKFace}
The \emph{typical $k$-face} associated with $\zeta$ is a
$\cC_\kappa^{d,k}$-valued random chamber
$C_{d,t,\kappa}^{k,{\rm typ},\zeta}$ with distribution
\begin{equation}\label{eq:TypicalKFaceDefinition}
\begin{split}
\PP(C_{d,t,\kappa}^{k,{\rm typ},\zeta}\in\mathcal A)
&:=
\frac{1}{\gamma_{d,k,t,\kappa}^{\,\zeta}\cH_\kappa^d(W)}
\EE\sum_{F\in\cF_k(\PHT_t)}
\int_{|F|}{\bf1}\{x\in W\}\\
&\hspace{2.6cm}\times
\int_{I(\MM_\kappa^d)}
{\bf1}\{\phi^{-1}F\in\mathcal A\}
\,\widetilde\lambda_\kappa(x,\dint\phi)\,
\zeta_F(\dint x),
\end{split}
\end{equation}
where $\mathcal A\subset\cC_\kappa^{d,k}$ is Borel and $W\subset\MM_\kappa^d$ is
any Borel set with $0<\cH_\kappa^d(W)<\infty$.
\end{definition}

As in \eqref{eq:TypicalCellDefinition}, the numerator of
\eqref{eq:TypicalKFaceDefinition}, viewed as a measure in $W$
for fixed $\mathcal A$, is isometry invariant and dominated by
$\EE\xi_{k,t}^{\kappa,\zeta}$. Thus the right-hand side is
independent of $W$, and \eqref{eq:KFaceIntensityDefinition}
shows that it has total mass one. For invariant $h$,
$h(\phi^{-1}F)=h(F)$ gives
\eqref{eq:TypicalKFaceInvariantExpectations}. Almost surely,
$p\in|C_{d,t,\kappa}^{k,{\rm typ},\zeta}|$, and its distribution
is invariant under the isometries fixing $p$. For the point-mass
choice in \eqref{eq:FaceCentreKernelConstruction}, its intrinsic
centre is $p$ whenever the face is not a boundaryless whole flat.

\begin{proposition}
\label{prop:FaceKernelIndependence}
The intensity $\gamma_{d,k,t,\kappa}^{\,\zeta}$ and the expectation
$\EE h(C_{d,t,\kappa}^{k,{\rm typ},\zeta})$ do not depend on the
covariant centre kernel $\zeta$, for every nonnegative measurable
isometry-invariant $h:\cC_\kappa^{d,k}\to[0,\infty]$.
\end{proposition}

\begin{proof}
This is the kernel version of
Proposition~\ref{prop:CentreIndependence}. For two kernels
$\zeta^{(1)},\zeta^{(2)}$, transport unit mass within each face
according to their product: the expected transport on rectangles is
\[
M(A\times B)
:=
\EE\sum_{F\in\cF_k(\PHT_t)}
\zeta_F^{(1)}(A)\zeta_F^{(2)}(B).
\]
It is diagonally invariant and has locally finite first marginal
by Lemma~\ref{lem:KFaceIntensityFinite}. Its marginals are the
mean measures associated with the two kernels. Equation
\eqref{eq:MassTransport} therefore proves equality of their
intensities. Replacing the unit mass for $F$ by
$h(F)\wedge n$ with $n\in\NN$, using
\eqref{eq:TypicalKFaceInvariantExpectations}, and letting
$n\to\infty$ proves equality of the expectations by monotone
convergence, exactly as in the proof of
Proposition~\ref{prop:CentreIndependence}.
\end{proof}

Henceforth we fix one such kernel and write
$\gamma_{d,k,t,\kappa}$ and $C_{d,t,\kappa}^{k,{\rm typ}}$.
Only the law up to isometry is independent of this choice, the
recentered representative in \eqref{eq:TypicalKFaceDefinition}
may depend on it.

\paragraph{Volume weighting.}
Before weighting the typical-face law by its volume, we verify
that the required normalizing expectation is finite and positive.
Recall that $\sigma_{d,k,t,\kappa}$ in
\eqref{eq:SkeletonIntensities} is the $k$-dimensional face volume
per unit ambient volume, not the face-count intensity.

\begin{lemma}\label{lem:KFaceVolumeNormalization}
The typical $k$-face satisfies
\begin{equation}\label{eq:KFaceVolumeNormalization}
\gamma_{d,k,t,\kappa}
\EE\cH_\kappa^k(C_{d,t,\kappa}^{k,{\rm typ}})
=
\sigma_{d,k,t,\kappa}.
\end{equation}
In particular,
$0<\EE\cH_\kappa^k(C_{d,t,\kappa}^{k,{\rm typ}})<\infty$.
\end{lemma}

\begin{proof}
Distribute the volume of each face according to its centre kernel.
The expected transport is
\[
M(A\times B)
:=
\EE\sum_{F\in\cF_k(\PHT_t)}
\cH_\kappa^k(|F|\cap A)\zeta_F(B).
\]
It is diagonally invariant. Since $\zeta_F$ has total mass one,
its first marginal is $\EE S_k(A)=
\sigma_{d,k,t,\kappa}\cH_\kappa^d(A)$ by
\eqref{eq:SkeletonMeasures} and \eqref{eq:SkeletonIntensities}.
Its second marginal is
$\EE\sum_F\cH_\kappa^k(F)\zeta_F(B)$.
Equations \eqref{eq:MassTransport} and
\eqref{eq:TypicalKFaceInvariantExpectations}, applied to
$h(F)=\cH_\kappa^k(F)$, therefore give
\eqref{eq:KFaceVolumeNormalization}. All integrands are
nonnegative. The asserted positivity and finiteness follow from
Lemma~\ref{lem:KFaceIntensityFinite} and
\eqref{eq:SkeletonIntensities}.
\end{proof}

\begin{definition}
\label{def:VolumeWeightedTypicalKFace}
The \emph{volume-weighted typical $k$-face}
$C_{d,t,\kappa}^{k,{\rm vol}}$ has distribution
\begin{equation}\label{eq:VolumeWeightedTypicalKFaceDefinition}
\PP(C_{d,t,\kappa}^{k,{\rm vol}}\in A)
:=
\frac{
\EE\!\left[
\cH_\kappa^k(C_{d,t,\kappa}^{k,{\rm typ}})
{\bf1}\{C_{d,t,\kappa}^{k,{\rm typ}}\in A\}
\right]
}{
\EE\cH_\kappa^k(C_{d,t,\kappa}^{k,{\rm typ}})
},
\qquad  A\subset\cC_\kappa^{d,k}\ \text{Borel}.
\end{equation}
\end{definition}

For a nonnegative measurable isometry-invariant functional
$h:\cC_\kappa^{d,k}\to[0,\infty]$, combining
\eqref{eq:VolumeWeightedTypicalKFaceDefinition},
\eqref{eq:TypicalKFaceInvariantExpectations} and
\eqref{eq:KFaceVolumeNormalization} gives
\begin{equation*}
\sigma_{d,k,t,\kappa}
\EE h(C_{d,t,\kappa}^{k,{\rm vol}})
=
\frac{1}{\cH_\kappa^d(W)}
\EE\sum_{F\in\cF_k(\PHT_t)}
\cH_\kappa^k(F)\zeta_F(W)h(F).
\end{equation*}
Proposition~\ref{prop:FaceKernelIndependence}, applied also to
$F\mapsto\cH_\kappa^k(F)h(F)$, shows that these expectations
are independent of the centre kernel. In the next subsection, we identify these distributions up to isometry
with intrinsic $k$-dimensional cell distributions and evaluate
$\gamma_{d,k,t,\kappa}$ explicitly.

\subsection{Sectional representation and face intensities}
\label{subsec:sectional-representation}

Fix $1\leq k\leq d-1$ and assume \eqref{eq:AdmissibleRegime}.
We identify the typical $k$-face and its volume-weighted version
with cells of an intrinsic $k$-dimensional Poisson hyperplane
tessellation. The comparison is up to isometry: the supporting
flat of the ambient typical face is random, whereas the
lower-dimensional model space is fixed.

\paragraph{The lower-dimensional comparison cells.}
Consider the Poisson hyperplane process in $\MM_\kappa^k$ with
parameter $t$ and invariant hyperplane measure normalized
intrinsically in dimension $k$, as in
Section~\ref{subsec:space-hyperplanes}. For $k=1$, this is the
Poisson point process with intensity measure
$t\cH_\kappa^1/\pi$ on a geodesic line when $\kappa\leq0$, or
on the elliptic circle of length $\pi/\sqrt{\kappa}$ when
$\kappa>0$.

Write $\PHT_t^{*,k}$ for the induced tessellation with the
following convention: in positive curvature an empty process
leaves one cell, namely the boundaryless whole space
$\MM_\kappa^k$. The star distinguishes this convention from
$\PHT_t$, which has no cells on an empty realization. The local
convention is needed because an uncut supporting $k$-flat is one
 face of the ambient tessellation.
Let $C_{k,t,\kappa}^{*,\rm typ}$ and
$C_{k,t,\kappa}^{*,\rm vol}$ be the typical and volume-weighted
typical cells of $\PHT_t^{*,k}$. They are defined by the Palm
construction and volume weighting of
Section~\ref{subsec:typical-k-faces}, applied to the cells in
$\MM_\kappa^k$ and using the kernels in
\eqref{eq:FaceCentreKernelConstruction} intrinsically in that
space. Thus a boundaryless whole elliptic space carries its
normalized volume measure, whereas an ordinary chamber uses its
intrinsic centre. In contrast, $C_{k,t,\kappa}^p$ denotes the
zero cell of the intrinsic process, as in
Definition~\ref{def:ZeroCell}. No star is needed on this object, as
its definition already uses the whole elliptic space on an empty
realization.

All these cells are compact almost surely. Indeed, in negative
curvature,
\[
t\geq\pi(d-1)\sqrt{-\kappa}
  >\pi(k-1)\sqrt{-\kappa},
\]
so the intrinsic process is strictly above its boundedness
threshold when $k\geq2$. For $k=1$, its cells are bounded
intervals for every $t>0$. Compactness in the other two geometries
holds for every $t>0$.

Denote the intensity of $\PHT_t^{*,k}$ by
$\gamma_{k,t,\kappa}^*$. Recalling the ordinary intrinsic
cell intensity $\gamma_{k,t,\kappa}$, we have
\begin{equation}\label{eq:LocalCellIntensity}
\gamma_{k,t,\kappa}^*
=
\begin{cases}
\displaystyle
\gamma_{k,t,\kappa}
+\frac{2\kappa^{k/2}}{\omega_{k+1}}e^{-t/\sqrt{\kappa}},
&\kappa>0,\\[0.7em]
\gamma_{k,t,\kappa},&\kappa\leq0.
\end{cases}
\end{equation}
For $k\geq2$, $\gamma_{k,t,\kappa}$ is given by
\eqref{eq:CellIntensityAllCurvatures} in dimension $k$. For
$k=1$, it equals $t/\pi$ by
\eqref{eq:OneDimensionalCellIntensity}. The added term in
\eqref{eq:LocalCellIntensity} is the probability
$e^{-t/\sqrt{\kappa}}$ of the empty intrinsic process divided by
$\cH_\kappa^k(\MM_\kappa^k)=\omega_{k+1}/(2\kappa^{k/2})$.
In particular, $\gamma_{k,t,\kappa}^*$ is finite and positive.

Let $h:\cC_\kappa^{d,k}\to[0,\infty]$ be a measurable
isometry-invariant functional, that is,
$h(\phi F)=h(F)$ for every $F\in\cC_\kappa^{d,k}$
and $\phi\in I(\MM_\kappa^d)$.
When evaluating $h$ on a cell of an intrinsic
$k$-dimensional tessellation in $\MM_\kappa^k$, we identify
that space isometrically with a $k$-flat of $\MM_\kappa^d$.
The value of $h$ does not depend on this identification,
since every isometry between two $k$-flats extends to an
ambient isometry. This convention applies to all expectation
identities below.

For a nonnegative measurable isometry-invariant functional
$h:\cC_\kappa^{d,k}\to[0,\infty]$ and a Borel set $B\subset\MM_\kappa^k$ of positive finite
$k$-volume, the intrinsic typical-cell construction gives
\begin{equation}\label{eq:LocalCellCampbell}
\EE\sum_{c\in\PHT_t^{*,k}}\zeta_c(B)h(c)
=
\gamma_{k,t,\kappa}^*\cH_\kappa^k(B)
\EE h(C_{k,t,\kappa}^{*,\rm typ}).
\end{equation}
The same formula holds with both sides zero when $B$ has
$k$-volume zero. Since the local cell supports cover the whole
space up to boundaries on every realization, the volume transport
used to prove \eqref{eq:KFaceVolumeNormalization} gives
\begin{equation}\label{eq:LocalVolumeIntensity}
\gamma_{k,t,\kappa}^*
\EE\cH_\kappa^k(C_{k,t,\kappa}^{*,\rm typ})=1.
\end{equation}
Here the outgoing measure is $\cH_\kappa^k$, including on the
empty realization in positive curvature. This also verifies the
finite positive normalization required to define
$C_{k,t,\kappa}^{*,\rm vol}$.

\begin{theorem}[Sectional representation of typical faces]
\label{thm:TypicalKFaceSection}
Under \eqref{eq:AdmissibleRegime}, for every nonnegative measurable
isometry-invariant functional
$h:\cC_\kappa^{d,k}\to[0,\infty]$,
\begin{equation}\label{eq:TypicalKFaceSection}
\EE h(C_{d,t,\kappa}^{k,\rm typ})
=
\EE h(C_{k,t,\kappa}^{*,\rm typ}),
\end{equation}
and
\begin{equation}\label{eq:WeightedTypicalKFaceSection}
\EE h(C_{d,t,\kappa}^{k,\rm vol})
=
\EE h(C_{k,t,\kappa}^{*,\rm vol})
=
\EE h(C_{k,t,\kappa}^p).
\end{equation}
Moreover, the $k$-face intensity satisfies
\begin{equation}\label{eq:gammadkt}
\gamma_{d,k,t,\kappa}
=
\sigma_{d,k,t,\kappa}\gamma_{k,t,\kappa}^*
=
\frac{t^{d-k}}{(d-k)!}
\frac{\omega_{k+1}}{\omega_{d+1}}
\gamma_{k,t,\kappa}^*.
\end{equation}
Here $\sigma_{d,k,t,\kappa}$ is the skeleton density in
\eqref{eq:SkeletonIntensities}, and
$\gamma_{k,t,\kappa}^*$ is given by \eqref{eq:LocalCellIntensity}.
\end{theorem}

\begin{remark}
For Euclidean volume-weighted faces, the representation by sections
of the zero cell is \cite[Theorem~1]{SchneiderWeightedFaces}.
In the isotropic case it identifies their shape law with that of
an intrinsic lower-dimensional zero cell. Related representations
for typical and volume-weighted spherical faces of fixed-number
great-hypersphere tessellations are given in
\cite[Theorems~3.1 and~3.3]{KT21}.
\end{remark}

\paragraph{The process induced on a flat.}
To prove Theorem~\ref{thm:TypicalKFaceSection}, we first verify
that intersecting the ambient process with a fixed flat preserves
the parameter $t$. For $L\in A_\kappa(d,k)$, let
$A_\kappa(L,k-1)$ be the space of hyperplanes of $L$, and let
$\widetilde\mu_{k-1}^{\kappa,L}$ denote its invariant measure
with the intrinsic $k$-dimensional normalization.

\begin{lemma}\label{lm:SectionProcess}
For each fixed $L\in A_\kappa(d,k)$, the random measure
\begin{equation}\label{eq:SectionalPoissonProcess}
\eta_t^{\kappa,L}
:=
\sum_{\substack{H\in\eta_t^\kappa\\
                 \dim(H\cap L)=k-1}}
\delta_{H\cap L}
\end{equation}
is a Poisson hyperplane process in $L$ with intensity measure
$t\widetilde\mu_{k-1}^{\kappa,L}$. More precisely, for every
nonnegative measurable $f:A_\kappa(L,k-1)\to[0,\infty]$,
\begin{equation}\label{eq:SectionHyperplaneMeasure}
\int_{A_\kappa(d,d-1)}
f(H\cap L){\bf1}\{\dim(H\cap L)=k-1\}
\,\widetilde\mu_{d-1}^\kappa(\dint H)
=
\int_{A_\kappa(L,k-1)}
f(J)\,\widetilde\mu_{k-1}^{\kappa,L}(\dint J).
\end{equation}
The first integrand is set to zero for the discarded intersections.
\end{lemma}

\begin{proof}
Let $\Theta_L$ be the image of $\widetilde\mu_{d-1}^\kappa$
under $H\mapsto H\cap L$, restricted to intersections of
dimension $k-1$. This measure is locally finite: hyperplanes of
$L$ meeting a compact set $B\subset L$ can only arise from
ambient hyperplanes meeting $B$, whose measure is finite. Every
isometry of $L$ extends to an ambient isometry preserving $L$,
so $\Theta_L$ is invariant. Consequently,
$\Theta_L=b\widetilde\mu_{k-1}^{\kappa,L}$ for some finite
$b\geq0$.

Choose a nondegenerate minimizing geodesic segment $s\subset L$,
shorter than $\pi/(2\sqrt{\kappa})$ when $\kappa>0$.
Outside a set of ambient hyperplanes of measure zero, an ambient
hyperplane meets $s$ if and only if its intersection with $L$
is a hyperplane of $L$ meeting $s$.
Recall that $\Theta_L=b\widetilde\mu_{k-1}^{\kappa,L}$
for some $b\in[0,\infty)$. Applying
\eqref{eq:HyperplaneSegmentCrofton} in the ambient space and
intrinsically in $L$, we obtain
\[
\begin{split}
	\frac{\cH_\kappa^1(s)}{\pi}
	&=
	\Theta_L\bigl(
	\{J\in A_\kappa(L,k-1):J\cap s\neq\emptyset\}
	\bigr)
	=
	b\frac{\cH_\kappa^1(s)}{\pi}.
\end{split}
\]
For $k=1$, the intrinsic hyperplanes are points of $L$,
and their invariant measure is $\cH_\kappa^1/\pi$.
Since $\cH_\kappa^1(s)>0$, we conclude that $b=1$,
proving \eqref{eq:SectionHyperplaneMeasure}.

The mapping theorem for Poisson processes
\cite[Theorem~5.1]{LPbook} now shows that
$\eta_t^{\kappa,L}$ is a Poisson hyperplane process in $L$
with intensity measure $t\widetilde\mu_{k-1}^{\kappa,L}$.
\end{proof}

For a flat $E\in A_\kappa(d,k)$, write $\PHT_t^{*,E}$ for
the local tessellation induced by \eqref{eq:SectionalPoissonProcess},
counting the boundaryless whole flat once if the induced process
is empty in positive curvature. Lemma~\ref{lm:SectionProcess}
identifies its law, after an isometry $E\simeq\MM_\kappa^k$,
with that of $\PHT_t^{*,k}$. The notation involving $E$ will
be used only for this local tessellation in the proof.

\begin{proof}[Proof of Theorem~\ref{thm:TypicalKFaceSection}]
By Proposition~\ref{prop:FaceKernelIndependence}, we may use the
intrinsic kernels in \eqref{eq:FaceCentreKernelConstruction}.
Every  $k$-face has exactly $d-k$ supporting hyperplanes.
On their intersection flat, it is a cell of the tessellation
induced by the remaining hyperplanes. Conversely, each such
local cell gives one $k$-face after the supporting
hyperplanes are inserted. If there are no remaining cuts in
positive curvature, the uncut flat is counted once, as required
by the starred convention.

Let $W\subset\MM_\kappa^d$ be a compact window of positive
volume. The multivariate Mecke formula
\cite[Theorem~4.4]{LPbook} applied to the $d-k$ supporting
hyperplanes, followed by \eqref{eq:Formula_1}, gives
\begin{equation}\label{eq:FaceCampbellReduction}
\EE\sum_{F\in\cF_k(\PHT_t)}\zeta_F(W)h(F)
=
\frac{t^{d-k}}{(d-k)!}
\int_{A_\kappa(d,k)}
\EE\sum_{c\in\PHT_t^{*,E}}
\zeta_c(W\cap E)h(c)
\,\widetilde\mu_k^\kappa(\dint E).
\end{equation}
The factorial removes the orderings of the supporting hyperplanes,
it does not count incidences of a face with full-dimensional cells.
The remaining hyperplanes form a Poisson process with the original
intensity measure. Lemma~\ref{lm:SectionProcess} and
\eqref{eq:LocalCellCampbell} therefore identify the inner
expectation in \eqref{eq:FaceCampbellReduction} as
\[
\gamma_{k,t,\kappa}^*\cH_\kappa^k(W\cap E)
\EE h(C_{k,t,\kappa}^{*,\rm typ}).
\]
For bounded $h$, substitution and
\eqref{eq:CroftonSpecialCase} yield
\begin{equation}\label{eq:FaceCampbellEvaluation}
\EE\sum_{F\in\cF_k(\PHT_t)}\zeta_F(W)h(F)
=
\frac{t^{d-k}}{(d-k)!}
\frac{\omega_{k+1}}{\omega_{d+1}}
\gamma_{k,t,\kappa}^*\cH_\kappa^d(W)
\EE h(C_{k,t,\kappa}^{*,\rm typ}).
\end{equation}
The identity extends to every nonnegative measurable invariant
$h$ by monotone convergence. Taking $h\equiv1$ and using
\eqref{eq:KFaceIntensityDefinition} proves \eqref{eq:gammadkt}.
Dividing \eqref{eq:FaceCampbellEvaluation} by
$\gamma_{d,k,t,\kappa}\cH_\kappa^d(W)$ and applying
\eqref{eq:TypicalKFaceInvariantExpectations} proves
\eqref{eq:TypicalKFaceSection}.

Taking $h$ to be $k$-volume in \eqref{eq:TypicalKFaceSection}
shows that the two typical cells have the same finite positive
mean volume. Applying the same identity to
$F\mapsto\cH_\kappa^k(F)h(F)$ and dividing by this common mean
therefore identifies their volume-weighted laws. This proves the
first equality in \eqref{eq:WeightedTypicalKFaceSection}.

It remains to identify the starred volume-weighted cell with the
intrinsic zero cell. For a bounded nonnegative invariant $h$,
transport $k$-volume weighted by $h(c)$ from each local cell to
its centre kernel. The first marginal has density
$\EE h(C_{k,t,\kappa}^p)$, because the starred cells cover
$\MM_\kappa^k$ on every realization. By
\eqref{eq:LocalCellCampbell}, the second marginal has density
$\gamma_{k,t,\kappa}^*
\EE[\cH_\kappa^k(C_{k,t,\kappa}^{*,\rm typ})
h(C_{k,t,\kappa}^{*,\rm typ})]$.
Lemma~\ref{lm:MassTransport}, applied intrinsically in dimension
$k$, gives
\begin{equation}\label{eq:LocalZeroCellVolumeBias}
\EE h(C_{k,t,\kappa}^p)
=
\gamma_{k,t,\kappa}^*
\EE\!\left[
\cH_\kappa^k(C_{k,t,\kappa}^{*,\rm typ})
h(C_{k,t,\kappa}^{*,\rm typ})
\right].
\end{equation}
For $k=1$ the same principle applies to the isometry group of a
line or a circle. Together with \eqref{eq:LocalVolumeIntensity},
\eqref{eq:LocalZeroCellVolumeBias} proves the remaining equality
in \eqref{eq:WeightedTypicalKFaceSection}. Truncation and
monotone convergence extend it to every $h$ in the statement.
Unlike \eqref{eq:ZeroCellVolumeBias} for the ordinary typical
cell, no conditioning on a nonempty process is needed here.
\end{proof}

\begin{corollary}[Mean $k$-volume of the typical $k$-face]
\label{cor:MeanTypicalKFaceVolume}
Under \eqref{eq:AdmissibleRegime}, for $1\leq k\leq d-1$,
\begin{equation}\label{eq:MeanVolumeTypicalKFace}
\begin{split}
\EE\cH_\kappa^k(C_{d,t,\kappa}^{k,\rm typ})
&=\frac{1}{\gamma_{k,t,\kappa}^*}=\frac{\omega_{k+1}}{2}
\left(
\sum_{j=0}^{\lfloor k/2\rfloor}
\kappa^j\frac{t^{k-2j}}{(k-2j)!}
+
\begin{cases}
\kappa^{k/2}e^{-t/\sqrt{\kappa}},
&\kappa>0\text{ and }k\text{ odd},\\[0.3em]
0,&\text{otherwise}
\end{cases}
\right)^{-1}.
\end{split}
\end{equation}
In particular, the mean depends on $k,t,\kappa$, but not on the
ambient dimension $d$, provided \eqref{eq:AdmissibleRegime} holds.
\end{corollary}

\begin{proof}
Combine \eqref{eq:gammadkt} with
\eqref{eq:KFaceVolumeNormalization} to obtain the first equality
in \eqref{eq:MeanVolumeTypicalKFace}. The second follows by
substituting \eqref{eq:CellIntensityAllCurvatures} in dimension
$k$ into \eqref{eq:LocalCellIntensity}, using
\eqref{eq:OneDimensionalCellIntensity} for $k=1$.
For positive curvature and even $k$, the added whole-space cell
cancels the negative exponential term in
$\gamma_{k,t,\kappa}$. For odd $k$, it supplies the positive
term displayed in \eqref{eq:MeanVolumeTypicalKFace}.
\end{proof}

For example, the mean length of the typical edge is
\[
\EE\cH_\kappa^1(C_{d,t,\kappa}^{1,\rm typ})
=
\begin{cases}
	\displaystyle \frac{\pi}{t},&\kappa\leq0,\\[0.7em]
	\displaystyle
	\frac{\pi}{t+\sqrt{\kappa}\,e^{-t/\sqrt{\kappa}}},
	&\kappa>0.
\end{cases}
\]
For comparison, the paper \cite{HugThaele2026}
determined the full length distribution of the typical edge
in a spherical Poisson hypersphere tessellation. On the unit
sphere, its mean is $2\pi/(2t+e^{-t})$, whereas our elliptic
formula with $\kappa=1$ gives $\pi/(t+e^{-t})$.
The difference arises when exactly $d-1$ hyperspheres are
present. In this case, the unique spherical edge is a full
great circle of length $2\pi$, while its projective image
has length $\pi$. With at least $d$ hyperspheres, spherical
edges occur in antipodal pairs, each corresponding to one
elliptic edge of the same length. For $d\geq3$, the mean area of the typical two-dimensional
face is $4\pi/(t^2+2\kappa)$, under the same assumptions.

\begin{corollary}
\label{cor:VolumeWeightedKFaceZeroSection}
Let $L\in G_\kappa(d,k)$ be fixed, so that $p\in L$.
Under \eqref{eq:AdmissibleRegime}, for every nonnegative measurable
isometry-invariant $h:\cC_\kappa^{d,k}\to[0,\infty]$,
\begin{equation}\label{eq:WeightedFaceZeroSection}
\EE h(C_{d,t,\kappa}^{k,\rm vol})
=
\EE h(C_{k,t,\kappa}^p)
=
\EE h(C_{d,t,\kappa}^p\cap L).
\end{equation}
The last intersection is equipped with the chamber structure
induced by the restricted hyperplanes. Under the natural coupling,
it is the zero cell of \eqref{eq:SectionalPoissonProcess}. After
an isometric identification of $L$ with $\MM_\kappa^k$ taking
$p$ to the intrinsic reference point, its distribution is that
of $C_{k,t,\kappa}^p$.
\end{corollary}

\begin{proof}
For $\kappa\leq0$, restrict the closed halfspaces defining
$C_{d,t,\kappa}^p$ to $L$. Each ambient hyperplane meeting $L$
gives the corresponding halfspace in $L$ containing $p$.  A
hyperplane disjoint from $L$ imposes no restriction on $L$.
Thus their intersection is exactly the zero cell of the induced
process. For $\kappa>0$, apply the same argument to the
hemispheres containing a fixed spherical lift of $p$, then
project to $L$. This retains the induced chamber structure,
not just the closed support. If the induced process is empty,
the local zero cell is the boundaryless whole $L$, consistently
with Definition~\ref{def:ZeroCell}.
Lemma~\ref{lm:SectionProcess} identifies its distribution with
that of the intrinsic zero cell. Combining this with
\eqref{eq:WeightedTypicalKFaceSection} proves
\eqref{eq:WeightedFaceZeroSection}.
\end{proof}

Finally, substitution into \eqref{eq:gammadkt} makes all face
intensities explicit. In particular, for $\kappa\leq0$,
\begin{equation*}
\gamma_{d,k,t,\kappa}
=
\frac{2}{\omega_{d+1}(d-k)!}
\sum_{j=0}^{\lfloor k/2\rfloor}
\kappa^j\frac{t^{d-2j}}{(k-2j)!},
\qquad 1\leq k\leq d-1.
\end{equation*}
For $\kappa>0$, equations \eqref{eq:gammadkt} and
\eqref{eq:LocalCellIntensity} give the corresponding expression.
For use in the next subsection, we also record the vertex intensity
separately. It is given by
\begin{equation}\label{eq:VertexIntensity}
\gamma_{d,0,t,\kappa}
=
\sigma_{d,0,t,\kappa}
=
\frac{2t^d}{\omega_{d+1}d!},
\end{equation}
as follows from \eqref{eq:SkeletonIntensities}, since $S_0$
counts vertices.

\subsection{Mean face vector of the typical cell}
\label{subsec:typical-cell-face-vector}

We now use the face intensities from the preceding subsection to
compute the expected face numbers of the full-dimensional typical
cell. Throughout we assume \eqref{eq:AdmissibleRegime}. Recall
that $f_k(c)=\#\cF_k(c)$ counts distinct faces of a cell $c$, and that $C_{d,t,\kappa}^{\rm typ}$ is the
Palm typical cell defined in \eqref{eq:TypicalCellDefinition}.
The cell intensity is $\gamma_{d,t,\kappa}$. For $k\geq1$, the
$k$-face intensity $\gamma_{d,k,t,\kappa}$ was evaluated in
\eqref{eq:gammadkt}.  For $k=0$, it is the vertex intensity in
\eqref{eq:VertexIntensity}.

\begin{proposition}
\label{thm:TypicalCellFaceVector}
Under \eqref{eq:AdmissibleRegime}, for every $0\leq k<d$,
\begin{equation}\label{eq:TypicalCellFaceVector}
\gamma_{d,t,\kappa}
\EE f_k(C_{d,t,\kappa}^{\rm typ})
=2^{d-k}\gamma_{d,k,t,\kappa}-\begin{cases}
\displaystyle
\frac{2^{d-k}\kappa^{k/2}}{\omega_{d+1}}
\frac{t^{d-k}}{(d-k)!}e^{-t/\sqrt{\kappa}},&\kappa>0,\\[0.8em]
0,&\kappa\leq0.
\end{cases}
\end{equation}
All these expectations are finite, including at the hyperbolic
critical intensity. The remaining face number is
$f_d(C_{d,t,\kappa}^{\rm typ})=1$ almost surely.
\end{proposition}

\begin{proof}
Fix $0\leq k<d$ and put $m=d-k$. At a relative interior point
of a $k$-face $F$, exactly $m$ hyperplanes meet.
Their linearly independent normals give $2^m$ local chamber
sectors. Let $b(F)$ be the number of distinct full-dimensional
cells incident with $F$.

For $\kappa\leq0$, all these sectors belong to different cells,
so $b(F)=2^m$. For $\kappa>0$, recall from
\eqref{eq:EllipticHyperplaneCount} that $N$ is the total number
of hyperplanes. Necessarily $N\geq m$ when a $k$-face exists.
If $N>m$, fix a spherical lift of a relative interior point of
$F$. At least one hyperplane is not among its supporting
hyperplanes, and all sectors at this lift have the same side
of that remaining great hypersphere. Hence no two of their
spherical chambers are antipodal, and they give $2^m$ distinct
elliptic cells. If $N=m$, the common supporting flat is uncut.
It is the unique $k$-face, and the $2^m$ sectors form
$2^{m-1}$ antipodal pairs. Thus, for every actual $k$-face,
\begin{equation}\label{eq:FaceIncidentCellCount}
b(F)
=
\begin{cases}
2^m,&\kappa\leq0,\\
2^m-2^{m-1}{\bf1}\{N=m\},&\kappa>0.
\end{cases}
\end{equation}

Choose a cell centre map $z$ and the covariant face kernels
$\zeta_F$ as in Section~\ref{subsec:typical-k-faces}. For $k=0$,
use the point mass at the vertex. Send one unit of mass from
$\zeta_F$ to the centre of each distinct cell incident with $F$.
The expected transport on Borel rectangles is
\begin{equation}\label{eq:FaceCellIncidenceTransport}
M_k(A\times B)
:=
\EE\sum_{F\in\cF_k(\PHT_t)}\zeta_F(A)
\sum_{\substack{c\in\PHT_t\\F\in\cF_k(c)}}
{\bf1}\{z(c)\in B\}.
\end{equation}
This measure is diagonally invariant. Its first marginal is
bounded by $2^m\gamma_{d,k,t,\kappa}\cH_\kappa^d$ and is
therefore locally finite, by Lemma~\ref{lem:KFaceIntensityFinite}
for $k\geq1$ and \eqref{eq:VertexIntensity} for $k=0$.

For $\kappa>0$, the measure
$A\mapsto\EE[{\bf1}\{N=m\}\sum_F\zeta_F(A)]$ is invariant
and has total mass $\PP(N=m)$, since there is exactly one
$k$-face on this event. Consequently, its density with respect
to ambient volume is
\begin{equation}\label{eq:UncutFaceIntensityContribution}
\frac{\PP(N=m)}{\cH_\kappa^d(\MM_\kappa^d)}
=
\frac{2\kappa^{k/2}}{\omega_{d+1}}
\frac{t^m}{m!}e^{-t/\sqrt{\kappa}},
\qquad \kappa>0.
\end{equation}
Here we used \eqref{eq:EllipticHyperplaneCount} and
\eqref{eq:EllipticSpaceVolume}. Equations
\eqref{eq:FaceIncidentCellCount} and
\eqref{eq:UncutFaceIntensityContribution} show that the first
marginal of \eqref{eq:FaceCellIncidenceTransport} has density
$2^m\gamma_{d,k,t,\kappa}$, minus the correction displayed in
\eqref{eq:TypicalCellFaceVector} when $\kappa>0$.

For a Borel set $B\subset\MM_\kappa^d$ of positive finite volume, its second
marginal is
\[
M_k(\MM_\kappa^d\times B)
=
\EE\sum_{c\in\PHT_t}{\bf1}\{z(c)\in B\}f_k(c)
=
\gamma_{d,t,\kappa}\cH_\kappa^d(B)
\EE f_k(C_{d,t,\kappa}^{\rm typ}),
\]
by \eqref{eq:TypicalCellInvariantExpectations}. Equating the
marginals by \eqref{eq:MassTransport} and dividing by
$\cH_\kappa^d(B)$ proves \eqref{eq:TypicalCellFaceVector}.
Finiteness follows from the bound on the first marginal and
$\gamma_{d,t,\kappa}>0$, without any prior moment assumption
on $f_k(C_{d,t,\kappa}^{\rm typ})$.
\end{proof}

The correction in \eqref{eq:TypicalCellFaceVector} is necessary
for elliptic faces. For example, with one hyperplane
there is one elliptic cell and one facet, although that
facet has two local incidences with the cell. This is the case
$m=1$ and $N=m$ in \eqref{eq:FaceIncidentCellCount}.

\paragraph{An explicit expression.}
For every integer $r\geq0$, recall that
\begin{equation}\label{eq:FaceVectorPolynomials}
P_r(t,\kappa)
:=
\sum_{j=0}^{\lfloor r/2\rfloor}
\kappa^j\frac{t^{r-2j}}{(r-2j)!}.
\end{equation}
In particular, $P_0(t,\kappa)=1$ and $P_1(t,\kappa)=t$.
As before, the coefficient with $j=0$ is one also for
$\kappa=0$.

\begin{theorem}[Mean face vector of the typical cell]
\label{thm:ExplicitTypicalCellFaceVector}
Under \eqref{eq:AdmissibleRegime}, for $0\leq k<d$ and
$\kappa\leq0$,
\begin{equation}\label{eq:NonpositiveTypicalCellFaceVector}
\EE f_k(C_{d,t,\kappa}^{\rm typ})
=
2^{d-k}\frac{t^{d-k}}{(d-k)!}
\frac{P_k(t,\kappa)}{P_d(t,\kappa)}.
\end{equation}
For $\kappa>0$ and $t>0$,
\begin{equation}\label{eq:PositiveTypicalCellFaceVector}
\EE f_k(C_{d,t,\kappa}^{\rm typ})
=
2^{d-k}\frac{t^{d-k}}{(d-k)!}
\frac{
P_k(t,\kappa)+\frac{(-1)^{k+1}}{2}
\kappa^{k/2}e^{-t/\sqrt{\kappa}}
}{
P_d(t,\kappa)-\frac{1+(-1)^d}{2}
\kappa^{d/2}e^{-t/\sqrt{\kappa}}
}.
\end{equation}
\end{theorem}

\begin{proof}
For $k\geq1$, substitute \eqref{eq:gammadkt} and
\eqref{eq:LocalCellIntensity} into
\eqref{eq:TypicalCellFaceVector}, and divide by the cell
intensity in \eqref{eq:CellIntensityAllCurvatures}.
When $\kappa>0$, the quantity
$\omega_{k+1}\gamma_{k,t,\kappa}^*/2$ equals $P_k$ plus
$\kappa^{k/2}e^{-t/\sqrt{\kappa}}$ exactly when $k$ is odd,
and equals $P_k$ when $k$ is even. The incidence correction
then subtracts
$\frac12\kappa^{k/2}e^{-t/\sqrt{\kappa}}$ in every dimension
$k$. The resulting coefficient is $(-1)^{k+1}/2$, as in
\eqref{eq:PositiveTypicalCellFaceVector}. The denominator is
$\omega_{d+1}\gamma_{d,t,\kappa}/2$ and is positive.
For $k=1$, use \eqref{eq:OneDimensionalCellIntensity}. For
$k=0$, substitute \eqref{eq:VertexIntensity} directly into
\eqref{eq:TypicalCellFaceVector}. The latter gives the same
expressions because $P_0=1$.
\end{proof}

\paragraph{Euclidean and spherical comparisons.}
Since $P_r(t,0)=t^r/r!$, equation
\eqref{eq:NonpositiveTypicalCellFaceVector} gives
\begin{equation}\label{eq:EuclideanTypicalCellFaceVector}
\EE f_k(C_{d,t,0}^{\rm typ})
=
2^{d-k}\binom{d}{k},
\qquad 0\leq k<d.
\end{equation}
These are precisely the face numbers of the $d$-dimensional
cube. Equation~\eqref{eq:EuclideanTypicalCellFaceVector} is
classical and, more generally, holds for stationary nondegenerate
Poisson hyperplane tessellations without an isotropy assumption,
see \cite[Theorems~10.3.1 and~10.3.2]{SW}.

\begin{remark}
\label{rem:TypicalCellFaceVectorLiterature}
The mean face numbers of typical spherical cells are closely
related to the random-cone formulas of Cover and Efron
\cite{CE67}, see also \cite[Corollary~4.1]{HS16} and
\cite[Theorem~5.2]{KT21}, the latter treating typical spherical
faces for nondegenerate directional distributions, without an
isotropy assumption. These results concern a fixed number of
great hyperspheres and conical or spherical faces, rather than
the Poisson process and elliptic faces used here.

For the compact elliptic space, our Palm definition gives
\[
\EE f_k(C_{d,t,\kappa}^{\rm typ})
=
\frac{\EE\sum_{c\in\PHT_t}f_k(c)}
     {\EE\#\PHT_t},
\qquad \kappa>0.
\]
Thus a comparison by Poissonization averages total incidence
counts and total cell counts before taking their ratio, rather
than simply averaging fixed-number typical-cell means with the
Poisson probabilities. The antipodal identifications for uncut
faces must also be retained, as in
\eqref{eq:TypicalCellFaceVector}.

These typical-cell formulas should be distinguished from the
expected face numbers of the zero cell. The Euclidean Poisson
zero cell was treated in \cite{KabluchkoZeroFVector}, the
fixed-number spherical weighted-cell analogue in \cite{KT21},
and the hyperbolic Poisson zero cell in \cite{BetaStarPolytopes}.
We discussed these results in
Remark~\ref{rem:ZeroCellFaceNumberLiterature}.
\end{remark}

\paragraph{The hyperbolic plane.}
For $d=2$ and $\kappa<0$, equation
\eqref{eq:NonpositiveTypicalCellFaceVector} becomes
\begin{equation}\label{eq:HyperbolicPlaneTypicalFaceNumbers}
\EE f_0(C_{2,t,\kappa}^{\rm typ})
=
\EE f_1(C_{2,t,\kappa}^{\rm typ})
=
\frac{4t^2}{t^2+2\kappa},
\qquad t\geq\pi\sqrt{-\kappa}.
\end{equation}
The equality of the two face numbers also holds for each
individual polygon. To compare with the planar literature,
write $A$ and $L$ for the area and perimeter of
$C_{2,t,\kappa}^{\rm typ}$. Equations
\eqref{eq:VolumeViaIntensity} and
\eqref{eq:MeanExteriorAngleFunctionals} give
\begin{equation}\label{eq:HyperbolicPlaneTypicalAreaPerimeter}
\EE A=\frac{4\pi}{t^2+2\kappa},
\qquad
\EE L=\frac{4\pi t}{t^2+2\kappa},
\qquad t\geq\pi\sqrt{-\kappa}.
\end{equation}
For the perimeter, each edge has exterior angle $1/2$, so
\eqref{eq:ExteriorAngleFunctionals} gives
$\upsilon_1^\kappa(c)=L(c)/(4\pi)$. This explains the second
expression in \eqref{eq:HyperbolicPlaneTypicalAreaPerimeter}.
All three means in \eqref{eq:HyperbolicPlaneTypicalFaceNumbers}
and \eqref{eq:HyperbolicPlaneTypicalAreaPerimeter} remain
finite at the critical intensity.

\begin{remark}
\label{rem:SantaloYanezTypicalCell}
The expressions in \eqref{eq:HyperbolicPlaneTypicalFaceNumbers}
and \eqref{eq:HyperbolicPlaneTypicalAreaPerimeter} agree with
\cite[Equation~(7.11)]{SantaloYanez}. Their parameter $\lambda$
is the expected number of lines meeting a unit geodesic segment
\cite[p.~144]{SantaloYanez}. Hence
\eqref{eq:HyperplaneSegmentCrofton} gives $\lambda=t/\pi$.
They also derive all six second-order moments of area, perimeter
and number of sides, see \cite[Equations~(7.12)--(7.17)]{SantaloYanez}.
In particular, their expression for the second area moment agrees,
when $t>\pi\sqrt{-\kappa}$, with the $d=2$ case of
\eqref{eq:TypicalCellVolumeSecondMoment}.

Their random polygon is defined through equal-weight sampling of
full polygons meeting an expanding observation domain. The
existence of the limiting distribution and convergence of the
empirical moments are assumed in the hyperbolic case
\cite[pp.~151 and~154]{SantaloYanez}. The derivation of their
unweighted polygon means additionally uses the boundary-effect
assumption \cite[Equation~(7.9), p.~155]{SantaloYanez}, which
they explicitly leave unproved in hyperbolic space. Our Palm
construction and mass-transport argument establish the three
first-moment formulas for the typical cell throughout
$t\geq\pi\sqrt{-\kappa}$, including its critical endpoint,
without requiring these assumptions. This does not assert
convergence of their window-based sampling procedure.
\end{remark}

\subsection{Mean total face volumes of the typical cell}
\label{subsec:typical-cell-face-volumes}

We conclude by replacing the face counts of the preceding
subsection by total face volumes. Recall from
\eqref{eq:TotalFaceVolumes} that
\[
L_k(c)=\sum_{F\in\cF_k(c)}\cH_\kappa^k(F),
\qquad 0\leq k\leq d.
\]
In particular,
$L_0(c)=f_0(c)$ and $L_d(c)=\cH_\kappa^d(c)$.

We also recall that $\gamma_{d,t,\kappa}$ denotes the cell
intensity, $\sigma_{d,k,t,\kappa}$ the $k$-skeleton density
from \eqref{eq:SkeletonIntensities}, and
$p_\kappa(t)=\PP(\eta_t^\kappa\neq\emptyset)$ the probability
defined in \eqref{eq:NonemptyProbability}. The factor
$(1+p_\kappa(t))/2$ appearing below equals $1$ for
$\kappa\leq0$ and $1-\frac12e^{-t/\sqrt{\kappa}}$ for
$\kappa>0$. It accounts for the antipodal identification
of cells incident with uncut faces in elliptic space.

\begin{proposition}
	\label{prop:TypicalCellTotalFaceVolumes}
	Under \eqref{eq:AdmissibleRegime}, for every $0\leq k<d$,
	\begin{equation}\label{eq:TypicalCellTotalFaceVolumes}
		\begin{split}
			\EE L_k(C_{d,t,\kappa}^{\rm typ})
			&=
			\frac{2^{d-k}\sigma_{d,k,t,\kappa}}{\gamma_{d,t,\kappa}}
			\frac{1+p_\kappa(t)}{2}
			=
			\frac{2^{d-k-1}\omega_{k+1}}{\omega_{d+1}}
			\frac{t^{d-k}}{(d-k)!}
			\frac{1+p_\kappa(t)}{\gamma_{d,t,\kappa}}.
		\end{split}
	\end{equation}
	For $k=d$,
	\begin{equation}\label{eq:TypicalCellTopFaceVolume}
		\EE L_d(C_{d,t,\kappa}^{\rm typ})
		=\frac{p_\kappa(t)}{\gamma_{d,t,\kappa}}.
	\end{equation}
	All these expectations are finite, including at the hyperbolic
	critical intensity.
\end{proposition}

\begin{proof}
Equation~\eqref{eq:TypicalCellTopFaceVolume} is
\eqref{eq:VolumeViaIntensity}. Fix $0\leq k<d$, put $m=d-k$,
and choose a covariant cell centre map $z$. Send one full copy
of the $k$-dimensional volume of each face to every distinct
incident cell. The expected transport on Borel rectangles is
\begin{equation}\label{eq:TypicalCellFaceVolumeTransport}
M_k(A\times B)
:=
\EE\sum_{F\in\cF_k(\PHT_t)}\cH_\kappa^k(F\cap A)
\sum_{\substack{c\in\PHT_t\\F\in\cF_k(c)}}
{\bf1}\{z(c)\in B\}.
\end{equation}
For $k=0$, $\cH_\kappa^0$ is counting measure, so the same
construction applies. The measure is diagonally invariant.
If $b(F)$ denotes the number of distinct cells incident with
$F$, then \eqref{eq:FaceIncidentCellCount} gives $b(F)\leq2^m$.
Consequently, its first marginal is bounded by
$2^m\EE S_k(A)=2^m\sigma_{d,k,t,\kappa}\cH_\kappa^d(A)$
and is locally finite. Lemma~\ref{lm:MassTransport} therefore
applies without a preliminary integrability assumption on
$L_k(C_{d,t,\kappa}^{\rm typ})$.

For $\kappa\leq0$, $b(F)=2^m$, so the first marginal equals
$2^m\sigma_{d,k,t,\kappa}\cH_\kappa^d(A)$. For $\kappa>0$,
\eqref{eq:FaceIncidentCellCount} gives instead
\[
M_k(A\times\MM_\kappa^d)
=
2^m\EE S_k(A)
-2^{m-1}\EE\bigl[{\bf1}\{N=m\}S_k(A)\bigr],
\]
where $N$ is the total number of hyperplanes. On $N=m$ there
is exactly one $k$-face, namely the uncut intersection
flat, whose projective volume is
$\omega_{k+1}/(2\kappa^{k/2})$. This formula also holds for
$k=0$, when the flat is a single point. The measure
$A\mapsto\EE[{\bf1}\{N=m\}S_k(A)]$ is isometry invariant.
Its density is its total mass divided by the ambient elliptic
volume, and hence, by \eqref{eq:EllipticHyperplaneCount} and
\eqref{eq:EllipticSpaceVolume}, equals
\begin{equation*}
\frac{\PP(N=m)\,\omega_{k+1}/(2\kappa^{k/2})}
     {\omega_{d+1}/(2\kappa^{d/2})}
=
e^{-t/\sqrt{\kappa}}
\frac{\omega_{k+1}}{\omega_{d+1}}\frac{t^m}{m!}
=
e^{-t/\sqrt{\kappa}}\sigma_{d,k,t,\kappa}.
\end{equation*}
Thus, in all three curvature regimes,
\[
M_k(A\times\MM_\kappa^d)
=
2^{m-1}\bigl(1+p_\kappa(t)\bigr)
\sigma_{d,k,t,\kappa}\cH_\kappa^d(A).
\]
On the other hand, regrouping the nonnegative sums in
\eqref{eq:TypicalCellFaceVolumeTransport} by cells and using
\eqref{eq:TypicalCellInvariantExpectations} gives
\[
M_k(\MM_\kappa^d\times B)
=
\EE\sum_{c\in\PHT_t}{\bf1}\{z(c)\in B\}L_k(c)
=
\gamma_{d,t,\kappa}\cH_\kappa^d(B)
\EE L_k(C_{d,t,\kappa}^{\rm typ}).
\]
Equating the marginals by \eqref{eq:MassTransport} and dividing
by $\gamma_{d,t,\kappa}\cH_\kappa^d(B)$, for a window of
positive finite volume, proves the first equality in
\eqref{eq:TypicalCellTotalFaceVolumes}. The second follows from
\eqref{eq:SkeletonIntensities}. The same calculation proves
finiteness, since the first marginal has finite density and
$\gamma_{d,t,\kappa}>0$.
\end{proof}

\begin{corollary}
\label{cor:TypicalCellFaceVolumesIntrinsic}
Under \eqref{eq:AdmissibleRegime}, for $0\leq k<d$,
\begin{equation*}
\EE L_k(C_{d,t,\kappa}^{\rm typ})
=
2^{d-k-1}\omega_{k+1}\bigl(1+p_\kappa(t)\bigr)
\EE\upsilon_k^\kappa(C_{d,t,\kappa}^{\rm typ}).
\end{equation*}
For $\kappa\leq0$, the relation simplifies and extends to
$k=d$:
\begin{equation}\label{eq:NonpositiveTypicalCellTotalFaceVolumes}
\EE L_k(C_{d,t,\kappa}^{\rm typ})
=
2^{d-k}\omega_{k+1}
\EE\upsilon_k^\kappa(C_{d,t,\kappa}^{\rm typ})
=
\frac{2^{d-k-1}\omega_{k+1}}{(d-k)!}
\frac{t^{d-k}}{P_d(t,\kappa)},
\qquad 0\leq k\leq d,
\end{equation}
where $P_d$ was defined in \eqref{eq:FaceVectorPolynomials}.
\end{corollary}

\begin{proof}
Compare \eqref{eq:TypicalCellTotalFaceVolumes} with
\eqref{eq:IntrinsicVolumesViaSkeletons}. For $\kappa\leq0$,
$p_\kappa(t)=1$. Substituting
\eqref{eq:MeanExteriorAngleFunctionals} and
$\gamma_{d,t,\kappa}=2P_d(t,\kappa)/\omega_{d+1}$ gives
\eqref{eq:NonpositiveTypicalCellTotalFaceVolumes} for $k<d$.
For $k=d$, use \eqref{eq:TopIntrinsicVolume} and
\eqref{eq:VolumeViaIntensity}.
\end{proof}

For $k=0$, $L_0=f_0$ and
$\sigma_{d,0,t,\kappa}=\gamma_{d,0,t,\kappa}$, so
\eqref{eq:TypicalCellTotalFaceVolumes} recovers the expected vertex number already determined in~\eqref{eq:TypicalCellFaceVector}.
 For $\kappa\leq0$ and $k=d-1$,
$L_{d-1}(c)=\cH_\kappa^{d-1}(\partial c)$, and
\eqref{eq:NonpositiveTypicalCellTotalFaceVolumes} gives
\[
\EE\cH_\kappa^{d-1}(\partial C_{d,t,\kappa}^{\rm typ})
=\frac{\omega_d t}{P_d(t,\kappa)},
\]
agreeing with \eqref{eq:MeanExteriorAngleFunctionals}, since
each facet has exterior angle $1/2$.

\begin{remark}
For $\kappa<0$ and at the critical intensity
$t=t_{\rm crit}(d,\kappa)$, Proposition~\ref{prop:TypicalCellTotalFaceVolumes} and
\eqref{eq:ZeroCellFaceVolumesAtCriticality} give
\begin{equation*}
\EE L_k(C_{d,t,\kappa}^{\rm typ})<\infty,
\qquad
\EE L_k(C_{d,t,\kappa}^{p})=+\infty,
\qquad 0\leq k\leq d.
\end{equation*}
Thus all total face volumes of the typical cell have finite means
at the critical intensity, while the corresponding zero-cell
means are infinite, although both cells are almost surely compact.
\end{remark}

\section*{Acknowledgements}
AG and ZK were  supported by the DFG under Germany's Excellence Strategy
EXC 2044/2--390685587, \emph{Mathematics M\"unster: Dynamics-Geometry-Structure}, by the DFG Priority Programme
SPP 2265, \emph{Random Geometric Systems}, and by the DFG Research
Training Group \emph{Rigorous Analysis of Complex Random Systems}
RTG 3027.  CT was supported by the DFG Priority Programmes SPP 2265,
\emph{Random Geometric Systems}, and SPP 2458,
\emph{Combinatorial Synergies}.

\section*{Declarations}

\subsection*{Statement on the use of generative AI}

The human authors conjectured Theorem~\ref{thm:cell-intensity} around
2022. At that time, they had proved that, in the hyperbolic bounded-cell
regime, the cell intensity is a polynomial in the intensity
parameter $t$, with coefficients that were not yet identified.  Guided by the corresponding spherical
formula and by the expectation that the cell intensity should depend
smoothly on the curvature, they conjectured the explicit coefficients
appearing in Theorem~\ref{thm:cell-intensity}. They also developed a
programme for proving the required smoothness property, but were
unable to carry it out because of substantial technical difficulties.
During this earlier period, the human authors obtained several other
results presented in this paper. These include the formula for the
expected volume of the zero cell, the sectional identities relating
lower-dimensional faces to intrinsic tessellations, and mean-value
formulas developed in a unified framework for all three curvature
regimes.

In 2026, ChatGPT was used extensively in the further development of
the project. During these interactions, ChatGPT suggested a proof of
Theorem~\ref{thm:cell-intensity}.  The proof presented in this paper is
based on that argument. ChatGPT also played a substantial role in the
subsequent mathematical exploration, in the derivation and
organization of additional results, and in the preparation and
editing of the manuscript.
All arguments and statements included in the paper were subsequently
checked, revised where necessary, and approved by the human authors,
who take full responsibility for the contents of the paper.

\subsection*{Conflict of interest statement}
The authors declare that they have no conflict of interest.

\subsection*{Data availability statement}
We do not analyze or generate any datasets.

\bibliographystyle{abbrv}
\phantomsection
\addcontentsline{toc}{section}{References}
\bibliography{references}

\end{document}